%% file: article.tex
\documentclass{siamart171218}
\usepackage{graphicx,grffile}
\usepackage{amsfonts, amssymb}
\usepackage{mathrsfs}
\usepackage{subfig}
\usepackage{xcolor}
\usepackage{grffile}
\usepackage{mathrsfs}
\usepackage{multirow}
\usepackage[normalem]{ulem}
\usepackage{bm}
\usepackage{indentfirst} 
\allowdisplaybreaks

\newsiamremark{remark}{Remark}

\newcommand\bbR{\mathbb{R}}
\newcommand\bbN{\mathbb{N}}
\newcommand\bbZ{\mathbb{Z}}

\newcommand\dd{\,\mathrm{d}}
\newcommand\Kn{{\rm Kn}}
\newcommand\Ma{\mathit{Ma}}

\newcommand\mM{\mathcal{M}}
\newcommand\mH{\mathcal{H}}

\newcommand\mT{\mathcal{T}}

\newcommand\brho{\boldsymbol{\omega}}
\newcommand{\bomega}{\boldsymbol{\omega}}

\newcommand\bG{\boldsymbol{G}}
\newcommand{\mG}{\mathcal{G}}
\newcommand\pd[2]{\dfrac{\partial {#1}}{\partial {#2}}}

\newcommand{\param}{\phi}

\numberwithin{equation}{section}

\graphicspath{{images/}}

\title{Learning invariance preserving moment closure model for Boltzmann-BGK equation}

\author{
Zhengyi Li \footnote{
Beijing International Center for Mathematical Research, Peking University, (lizhengyi@pku.edu.cn)}, Bin Dong \footnote{Beijing International Center for Mathematical Research \&
Center for Machine Learning Research, Peking University, (dongbin@math.pku.edu.cn)},
Yanli Wang\footnote{Beijing Computational Science Research Center, (ylwang@csrc.ac.cn)}
}

\usepackage{subfiles}
\begin{document}
\maketitle
\input{article_abstract.tex}
\input{article_Boltzmann.tex}
\input{article_scheme.tex}
\input{article_experiment.tex}
%\input{extra.tex}
\input{article_conclusion.tex}
\input{appendix.tex}
\bibliography{article,references}
\bibliographystyle{plain}

%\printbibliography

\end{document}

%% file: article_abstract.tex
\begin{abstract} 
\vspace*{4mm}

 As one of the main governing equations in kinetic theory, the Boltzmann equation is widely utilized in aerospace, microscopic flow, etc. Its high-resolution simulation is crucial in these related areas. However, due to the high dimensionality of the Boltzmann equation, high-resolution simulations are often difficult to achieve numerically. The moment method which was first proposed in Grad (Commun. Pure Appl. Math. 2(4):331-407,1949) is among the popular numerical methods to achieve efficient high-resolution simulations. We can derive the governing equations in the moment method by taking moments on both sides of the Boltzmann equation, which effectively reduces the dimensionality of the problem. However, one of the main challenges is that it leads to an unclosed moment system, and closure is needed to obtain a closed moment system. It is truly an art in designing closures for moment systems and has been a significant research field in kinetic theory. Other than the traditional human designs of closures, the machine learning-based approach has attracted much attention lately in Han et al. (Proc. Natl. Acad. Sci. U.S.A. 116(44):21983-21991, 2019) and Huang et al. (J. Non-Equilib. Thermodyn. 46(4):355-370, 2021). In this work, we propose a machine learning-based method to derive a moment closure model for the Boltzmann-BGK equation. In particular, the closure relation is approximated by a carefully designed deep neural network that possesses desirable physical invariances, i.e., the Galilean invariance, reflecting invariance, and scaling invariance, inherited from the original Boltzmann-BGK equation and playing an important role in the correct simulation of the Boltzmann equation. Numerical simulations on the 1D-1D examples including the smooth and discontinuous initial condition problems, Sod shock tube problem, the shock structure problems, and the 1D-3D examples including the smooth and discontinuous problems demonstrate satisfactory numerical performances of the proposed invariance preserving neural closure method. 

\medskip

\end{abstract}

\section{Introduction}
% introduce Boltzmann 
% introduce classical  method 
% introduce machine learning, recent work 
% conclusion 
Kinetic theory is widely used when modeling non-equilibrium dynamics in a variety of fields, such as rarefied gases, plasma, and semiconductors. As one of the most fundamental kinetic equations, the Boltzmann equation describes the behaviors of the dynamic system from a statistical standpoint. However, due to the high dimensionality of the Boltzmann equation, its efficient numerical simulation is a long-standing challenge. The challenge is reflected in both the variables in the Boltzmann equation and the quadratic collision term. On the one hand, the Boltzmann equation itself is a seven-dimension integro-differential equation, including time, physical space, and microscopic velocity space. On the other hand, the quadratic collision model contains a high dimensional integral and a collision kernel with singularities. 
 
Various numerical methods have been proposed to solve the Boltzmann equation, characterized into the stochastic and deterministic methods. For the stochastic method, the direct simulation Monte Carlo method (DSMC) \cite{Bird} is widely used in the steady-state problem and highly rarefied flows. However, the usage of the DSMC method is limited due to its numerical noise and low efficiency. For the deterministic method, one of the classical methods is the discrete velocity method (DVM), which evaluates the numerical solution of the distribution function at specific points. DVM is easy to implement but has a low order of convergence and turns out to be inefficient. The Fourier spectral method \cite{Mouhot, Gamba}, where the trigonometric functions are used as the basis function to approximate the distribution function in the microscopic velocity space, has a high order of convergence and high efficiency. However, the microscopic velocity space is truncated into a finite space, which may lead to aliasing. Another important deterministic method is the moment method, where the Boltzmann equation is converted to the moment system by computing the moment coefficients of the Boltzmann equation. However, the moment system is not closed due to the convection term of the Boltzmann equation. 
 
The moment method for the Boltzmann equations was first proposed by Grad \cite{Grad}, where the distribution function is approximated by a series of basis functions whose weight function is the Maxwellian. In the framework of the Grad method, the closure is realized by simply truncating the expansion. However, such truncation may lead to inaccurate simulations and non-hyperbolicity even near the Maxwellian, which limits the application of the Grad moment system. The maximum entropy method is also adopted to derive the closed moment system \cite{Levermore, Mcdonald2013}. An optimization problem is solved to obtain the distribution function, and the related closed moment system is guaranteed to be hyperbolic \cite{Junk, Levermore1998BP}. However, the non-linear optimization problem can be expensive to solve \cite{Abramov, Schaerer, McDonald}. In \cite{Fan_new}, a globally hyperbolic closure method is proposed by modifying the governing equations of the moment coefficients at the highest order. The quadrature-based moment closure which retains a quadrature approximation of the distribution function is proposed \cite{Robert1997}. Then several related methods such as HyQMOM are proposed and we refer \cite{yuan2012extended, fox2009high, PATEL2019100006} and the references therein for more details of HyQMOM and the extended methods. In \cite{Koellermeier, KoellermeierRGD2014}, a hyperbolic quadrature-based closure method is proposed by modifying the governing equations of moment coefficients at the last two orders. However, these closure methods rely on mathematical derivations, such as the Boltzmann equation and the distribution function itself. These mathematical approaches are able to provide theoretical guarantees, such as hyperbolicity, on the closed moment systems. However, the accuracy of simulations still has plenty of room for improvement, which is the main objective of the latest machine learning-based methods. This paper introduces a new method to close the Grad moment system by learning the moment coefficients of the highest order from the simulation data. 

In recent years, machine learning methods have been introduced to derive the surrogate models for the kinetic equations, including the Boltzmann equation, Vlasov equation, radiation transfer equations, etc. In the pioneering work \cite{han2019uniformly}, a machine learning framework for moment closure of arbitrary order through neural networks is introduced, where a generalized moment system is learned by an encoder-decoder network and then closed by another neural network. A convolutional neural network is utilized for the closure of the Euler-Poisson equation in \cite{bois2020neural}. In \cite{Zhang2020DatadrivenDO}, the hydrodynamic equations are learned using the training data generated by DSMC. From the standpoint of the maximum entropy method, the entropy function is approximated by a neural network \cite{Porteous2021maxent1, Schotthofer2021maxent2}, which reduces the computational cost of the maximum entropy closure. Deep learning-based methods have been adopted to solve the kinetic system directly \cite{lou2021PINN-BGK, Xiao2021UsingNN}. In \cite{lou2021PINN-BGK}, PINN is used to solve the Boltzmann equation system directly, and in \cite{Xiao2021UsingNN}, the neural network is utilized to approximate the complex quadratic collision term. 

The closed moment system derived through the neural network is also expected to preserve some properties of the Boltzmann equation, such as hyperbolicity and physical invariance. Preservation of hyperbolicity plays an important role in the highly efficient simulation of the Boltzmann equation. 
For example, a closure of the Euler equation through a neural network in \cite{huang2021learning} ensures the hyperbolicity of the closed moment system. A moment closure for the RTE, which also maintains hyperbolicity, is studied in \cite{huang2021machine1, huang2021machine2, huang2021machine3}. 
Another point we are interested in is preserving physical invariances. Since the Boltzmann equation maintains several physical invariances such as the Galilean invariance, reflecting invariance, and scaling invariance, which is also natural for the classical closed moment systems, we could expect that the closed moment system derived through the neural network could also preserve these invariances. However, they are not automatically satisfied, and several works have been done through specially designing neural networks to achieve this. We note that the Galilean invariance is preserved in \cite{han2019uniformly} by defining the generalized moments properly. In \cite{huang2021learning}, the Galilean invariance is also preserved but only to close the Euler system. In \cite{Porteous2021maxent1, Schotthofer2021maxent2}, the reflecting and scaling invariances are incorporated through learning the closure function by maximum entropy method for the kinetic equations of scattering homogeneous medium. The neural networks in \cite{huang2021machine1} preserve the scaling invariance by normalizing the higher-order moments by zero-order moments to improve the performance of the closed model for the radiative transfer equation (RTE). Moreover, besides the application of invariance preserving closure by neural network for Boltzmann equations, in a work on the closure of the Reynolds-averaged Navier–Stokes equation, it has also shown that embedding these invariances into the model yields better performance than learning the invariance directly from the training data \cite{ling2016machine}. These works showed that incorporating invariances can be beneficial. However, there is no work that has considered all three invariances at once for the moment system of the Boltzmann equation. 

In this paper, we will introduce a new neural network-based moment closure method by considering as many physical invariances as possible. Such invariance includes the Galilean invariance, reflecting invariance, and scaling invariance. The neural closure network is specifically designed such that the aforementioned invariances are strictly enforced. Instead of learning the invariance directly from the training data, special techniques are utilized to impose these invariances in the neural closure network. For example, the derivative of the macroscopic velocity instead of itself is utilized as the input data to preserve the Galilean invariance. The original moment coefficients and their flips are used as input to enforce reflecting invariance. Moreover, the moment coefficients are all dimensionless by the density and temperature as input and then scaling back the output to preserve the scaling invariance. Thus, we will refer to the proposed method as the invariance preserving neural closure method (IPNC). Numerical results show that these embedded invariances help IPNC produce more accurate results and significantly improve the generalization of the model. On the other hand, we will prepare a benchmark data set based on the setting of \cite{han2019uniformly}, on which a relative comprehensive comparison between traditional moment closure methods and machine learning-based methods is presented. Such benchmark data sets may benefit the community by providing a common test-bed for new machine learning methods.

The rest of this paper is organized as follows. In Sec. \ref{sec:pre}, we will review some fundamental properties of Boltzmann's equation, especially the physical invariances. The globally hyperbolic moment method (HME) \cite{Fan_new} is also reviewed in this section. In Sec. \ref{sec:sym_net}, the particular design of the network to preserve the physical invariances of the Boltzmann equation is introduced. The numerical experiments to validate this invariance preserving closure method are shown in Sec. \ref{sec:num}. We conclude the paper in Sec. \ref{sec:conclusion}. Several supplementary materials, including the details of the random parameters in the numerical experiments, are provided in App. \ref{app:supp}.

%%% Local Variables: 
%%% mode: latex
%%% TeX-master: "article"
%%% End: 

%% file: article_Boltzmann.tex
\section{Boltzmann-BGK equation and moment method}
\label{sec:pre}
In this work, an invariance preserving moment closure method by neural networks will be proposed. This approach is motivated by the regularized moment method \cite{NRxx}. The Boltzmann-BGK equation, especially the invariance properties of the BGK equation, will be briefly reviewed in this section. Then the regularized moment method based on which we will derive the closed moment system is also presented. 

\subsection{Boltzmann-BGK equation}
 Boltzmann equation as one of the most fundamental kinetic equations, describes the movement of the particles from the statistical physical viewpoint. Here, we consider the simplified 1D BGK model \cite{BGK}, which assumes a simple relaxation to equilibrium, and takes the form 
\begin{equation}
    \label{eq:BGK}
    \pd{f(t, x, v)}{t} + v  \pd{f(t,x ,v)}{x} = \frac{1}{\Kn} [\mM(f) -f], \qquad x \in D \subset \bbR, \quad v \in \bbR, \quad t> 0,
\end{equation}
where $f(t, x, v)$ is the distribution function of the particles, $x$ is the physical space, $v$ is the microscopic velocity.  The Knudsen number $\Kn$, which is defined as the ratio of the mean free path and the typical length, is always utilized to describe the regime of the fluid dynamics. Here $\mM(f)$ is the Maxwellian equilibrium as 
\begin{equation}
    \label{eq:maxwellian}
    \mM(f) = \frac{\rho(t, x)}{\sqrt{2 \pi \theta(t, x)}}\exp\left(-\frac{(v - u(t, x))^2}{2\theta(t, x)} \right),
\end{equation}
where $\rho(t, x)$, $u(t, x)$ and $\theta(t, x)$ is the density, macroscopic velocity, and temperature respectively. The relationship between these macroscopic variables and the distribution function is 
\begin{equation}
    \label{eq:macro}
    \rho = \int_\bbR f(t, x, v) \dd v, \qquad 
    m := \rho u = \int_\bbR v f(t, x, v) \dd v, \qquad 
E := \frac{1}{2}   \rho u^2 +\frac{1}{2} \rho \theta =\frac{1}{2} \int_{\bbR} v^2 f(t, x, v) \dd v,
\end{equation}
where $m$ is the momentum, and $E$ is the total energy. 
Moreover, the relationship between the macroscopic variables, the distribution $f$ and the Maxwellian equilibrium \eqref{eq:maxwellian} is 
\begin{equation}
    \label{eq:f_M}
    \int_\bbR \left(\begin{array}{c}
         1  \\
         v \\
        \frac{1}{2} v^2
    \end{array}\right) f \dd v = \int_\bbR \left(\begin{array}{c}
         1  \\
         v \\
        \frac{1}{2}  v^2
    \end{array}\right) \mM(f) \dd v = \left(\begin{array}{c}
         \rho  \\
         m \\
        E
    \end{array}\right).
\end{equation}
We can also get the compressible Euler equation by multiplying \eqref{eq:BGK} with  $(1, v, \frac{1}{2}v^2)^T$ and integrating over $v\in\bbR$ as 
\begin{equation}
    \label{eq:Euler}
    \left\{
    \begin{aligned}
&\pd{\rho}{t} + \pd{(\rho u)}{x} = 0, \\
&\pd{\rho u}{t} + \pd{(\rho u^2 + p)}{x} = 0,  \\
&\pd{E}{t} + \pd{((E +p) u + q)}{x} = 0,
    \end{aligned}
    \right.
\end{equation}
where $p$ and  $q$ are the pressure and  heatflux respectively, which are defined as 
\begin{equation}
    \label{eq:heatflux}
   p(t,x) = \rho(t, x)\theta(t, x) = \int_{\bbR}(v-u)^2 f(t,x,v) \dd v, \qquad  q(t, x) =\frac{1}{2} \int_\bbR (v -u)^3 f(t,x,v)\dd v. 
\end{equation}

The BGK equation has several physical invariances, such as the Galilean invariance, scaling invariance, and reflecting invariance. Preserving these physical invariances is essential when solving the BGK equation numerically. For the classical numerical methods, these properties are naturally preserved in most cases. However, they are not readily satisfied with neural network-based closures. In the next section, we will briefly introduce the important physical invariances that we are interested in.

\subsection{Physical invariances of the Boltzmann equation}
\label{sec:invar}
The Boltzmann equation itself has some physical invariances, that is, the property that the form of the equation remains invariant under certain transformations. The invariances we mainly consider are the Galilean invariance, reflecting invariance, and scaling invariance. These invariances are the manifestations of some fundamental physical laws, which will be explained successively. These properties are also verified to be important when designing the closed moment system by machine learning methods \cite{huang2021learning, han2019uniformly}. 

%\subsubsection{Galilean Invariance}
We start with the Galilean invariance. The Galilean invariance means that the forms of equations are invariant under Galilean transformation, which consists of rotation, translation and Galilean boost. As is stated in \cite{han2019uniformly}, the Boltzmann equation having the Galilean invariance means that if  $f(t, x, v)$ is the solution to the Boltzmann equation, then $\tilde{f}(t, x, v)$ is also the solution to the system, where $\tilde{f}(t, x, v)$ is defined as 
\begin{equation}
    \label{eq:Galilean}
    \tilde{f}(t, x, v) = f(t, x - t u', v - u')
\end{equation}
 with any $u' \in \bbR$. 
% Galilean invariance for Boltzmann equation is that if $f(tm$

% Specifically, for every $t_0,x_0,v_0 \in \mathbb{R}$, define $f^*(t,x,v)=f(t-t_0,x-x_0-v_0 t,v-v_0)$. 
% If $f$ is a solution of the Boltzmann equation, then so is $f^*$.

% Rotation.
% \begin{equation}
% f^*=f , t^*=t , x^*=Gx, v^*=Gv, \Kn^*=\Kn
% \end{equation}

% Translation. $x_0,t_0$ are arbitrary values.
% \begin{equation}
% f^*=f ,  t^*=t+t_0, x^*=x+x_0, v^*=v ,\Kn^*=\Kn
% \end{equation}

% Uniform motion. $v_0$ is an arbitrary value.
% \begin{equation}
% f^*=f ,  t^*=t,x^*=x+v_0 t, v^*=v+v_0 ,\Kn^*=\Kn
% \end{equation}

%\subsubsection{Reflection Invariance}
The second invariance we are concerned is the reflecting invariance. The reflecting, also called parity transformation, implies the flip in the sign of one spatial coordinate. Since the classical mechanics is invariant under parity transformation, so is the Boltzmann equation. In this case, the reflecting invariance means that if $f(t, x, v)$ is a solution to the Boltzmann equation, then $\tilde{f}(t, x, v)$ is also the solution to the Boltzmann equation, with 
\begin{equation}
    \label{eq:rot}
\tilde{f}(t, x, v) = f(t, -x, -v).    
\end{equation}

% That is, if define $f'(t,x,v)=f(t,-x,-v)$. If $f$ is a solution of the Boltzmann equation, then so is $f'$.

%\subsubsection{Scaling Invariance}

Finally, the scaling invariance corresponds to the dimensional correctness of the equation, which means that the physical phenomenon should be independent of the unit of measurement chosen. Precisely, supposing $f(t,x,v)$ is a solution to the BGK equation \eqref{eq:BGK} with Knudsen number $\Kn$, consider $\tilde{f}(t^*,v^*,x^*)$ with the following scaling relationship
\begin{equation}
    \label{eq:scaling_inv}
\tilde{f}=f/f_0,\quad  
t^*=t/t_0,\quad  
x^*=x/x_0,\quad  
v^*=v t_0/x_0,\quad 
\Kn^*=\Kn/t_0,
\end{equation}
for any $f_0, t_0, x_0 \in \mathbb{R}^+$. Then, $\tilde{f}$ is also a solution to Boltzmann-BGK equation with Knudsen number ${\Kn}^*$.

These invariances are fundamental for the Boltzmann equation. In the simulation, we would also expect that the numerical scheme would conserve these invariances. For most of the classical numerical methods, such as the Fourier spectral method \cite{Mouhot, Gamba}, the direct simulation Monte Carlo method \cite{Bird}, and the moment method \cite{Struchtrup2005, Fan_new}, these invariances are naturally preserved. In this paper, we focus on the moment method. We will first introduce the moment method briefly, and then the invariances of the moment system will be discussed.

\subsection{Moment method}
\label{sec:moment_method}
The moment method was first proposed by Grad \cite{Grad} in 1949. One of the primary concerns of the moment method is how to obtain the closed moment system. In the framework of Grad's method, the closed moment system is derived simply by truncating the higher order of moment coefficients. However, this will lead to a non-hyperbolic system. In \cite{Fan_new}, a globally hyperbolic moment method  (HME) is proposed, and we follow the same expansion framework here. In HME, the distribution function is approximated as
\begin{equation}
    \label{eq:expansion} 
    f(t, x, v) \approx\sum\limits_{\alpha \leqslant M} f_\alpha(t,x) \mH_{\alpha}(\xi), \qquad \xi = \frac{v-u}{\sqrt{\theta}}, \qquad \alpha \in \bbN,
\end{equation}
where $\mH_{\alpha}(\cdot)$ is the basis function defined as 
\begin{equation}
    \label{eq:basis} 
    \mH_{\alpha}(\xi) = \frac{1}{\sqrt{2\pi}}\theta^{-\frac{\alpha+1}{2}} He_{\alpha}(\xi)\exp\left(-\frac{\xi^2}{2}\right),
\end{equation}
with $He_{\alpha}(\cdot)$ the Hermite polynomial as 
\begin{equation}
    \label{eq:Hermite}
    He_{\alpha}(v) = (-1)^{\alpha}\exp\left(\frac{v^2}{2}\right) 
    \dfrac{\dd^{\alpha}}{\dd v^{\alpha}} \exp\left(-\frac{v^2}{2}\right). 
\end{equation}
With the orthogonality of the basic functions,  we can get $f_{\alpha}$ as 
\begin{equation}
    \label{eq:coef}
    f_{\alpha}(t, x) = \frac{\theta^{\alpha/2}}{\alpha!}\int_{\bbR} f(t, x, v) He_{\alpha}\left(\frac{v-u}{\sqrt{\theta}}\right) \dd v. 
\end{equation}
It holds for the moment coefficients that 
\begin{equation}
    \label{eq:coe}
    f_0 = \rho, \qquad f_{1} = 0, \qquad f_2 = 0, \qquad q = 3 f_3. 
\end{equation}
Substituting \eqref{eq:expansion} into \eqref{eq:BGK} and collecting coefficients for the same basis function, we can get the moment system with some rearrangement as \cite{Fan_new} 
\begin{equation}
    \label{eq:moment_Euler}
    \begin{aligned}
    & \pd{f_0}{t} + \left(u \pd{f}{x} + f_0 \pd{u}{x}\right) = 0,\\
    & f_0\left(\pd{u}{t} + u\pd{u}{x}\right) + f_0 \pd{\theta}{x} + \theta\pd{f_0}{x}=0,\\
    &\frac{f_0}{2} \left(\pd{\theta}{t} + u\pd{\theta}{x}\right) + \pd{q}{x} + p \pd{u}{x} = 0,
    \end{aligned}
\end{equation}
and 
\begin{equation}
    \label{eq:moment_system}
    \begin{aligned}    
& \frac{\partial f_{\alpha}}{\partial t} - \frac{1}{f_0} 
\frac{\partial p}{\partial x} f_{\alpha-1}
  - \frac{1}{f_0}  \left(
    \frac{\partial q}{\partial x} +
     p \frac{\partial u}{\partial x}
  \right)  f_{\alpha-2} \\
& \quad +  \left[
    \frac{\partial u}{\partial x} \left(
      \theta f_{\alpha-2} + (\alpha + 1) f_{\alpha}
    \right) + \frac{1}{2} \frac{\partial \theta} {\partial x} \left(
      \theta f_{\alpha-3} + (\alpha + 1) f_{\alpha-1}
    \right)
  \right] \\
& \quad + \left(
    \theta  \frac{\partial f_{\alpha - 1}}{\partial x} +
    u \frac{\partial f_{\alpha}}{\partial x}\right) +
     (\alpha + 1) \frac{\partial f_{\alpha+1}}{\partial x}
   =\nu(1 - \delta(\alpha))f_{\alpha}, \quad \forall M \geqslant |\alpha| > 2,
\end{aligned}
\end{equation}
where $\delta(\alpha)$ is defined as 
\begin{equation}
    \label{eq:delta}
    \delta(\alpha) = \left\{\begin{array}{ll}
    0,    & {\rm if}~ \alpha > 2,  \\
    1,     & {\rm otherwise.}
    \end{array}  \right.
\end{equation}
Here \eqref{eq:moment_Euler} is the same as the Euler system \eqref{eq:Euler}, and it is also verified that the Euler system is contained in the moment system. 
Let $\brho = (\rho, u, \theta, f_3, \cdots, f_M)^T \in \bbR^{M+1}$, then the moment system \eqref{eq:moment_Euler} and \eqref{eq:moment_system} could be written into the quasi-linear form as below 
\begin{equation}
    \label{eq:linear_system}
    \pd{\brho}{t} + {\bf A}_{M}(\brho) \pd{\brho}{x} + (M+1)\pd{f_{M+1}}{x} \boldsymbol{e}_{M+1} = G(\brho),
\end{equation}
where ${\bf A}_M$ is the coefficients matrix derived from \eqref{eq:moment_Euler} and \eqref{eq:moment_system}. $\boldsymbol{e}_{M+1} \in \bbR^{M+1}$ is a $M+1$ vector with the $(M+1)$-th entry equaling $1$ and others equaling $0$. $\bG(\brho)$ is the collision term. We refer \cite{Fan_new} for the detailed deduction of the moment system. 

It is obvious that \eqref{eq:linear_system} is not closed. In \cite{Grad}, $f_{M+1}$ is set directly as $0$, which lead to the seminal Grad-type moment system. However, the usage of the Grad-type moment system is limited due to the loss of hyperbolicity, even in the region near the equilibrium \cite{FanDissertation}. In \cite{Fan_new}, a globally hyperbolic regularization method is proposed, which is shown in the Proposition below. 
\begin{proposition}
\label{thm:global} 
When $\theta > 0$, the moment system 
\begin{equation}
    \label{eq:global} 
    \pd{\brho}{t} + {\bf A}_M \pd{\brho}{x} -\mathcal{R}_m\boldsymbol{e}_{M+1} =\boldsymbol{G}(\brho), 
\end{equation}
is globally hyperbolic, where 
\begin{equation}
    \label{eq:R}
    \mathcal{R}_M = \frac{M+1}{2}\left(2 f_M \pd{u}{x} + f_{M-1} \pd{\theta}{x}\right).
\end{equation} 
The characteristic velocity of the moment system \eqref{eq:global} is \begin{equation}
    \label{eq:chara}
    s_j = u + c^j_{M+1}\sqrt{\theta}, \qquad j = 1, \cdots, M+1,
\end{equation} 
where $c^j_{M+1}$ is the $j$-th value of the Hermite polynomials $He_{M+1}(x)$. 
\end{proposition}

For the original Grad-type closure, as well as for its improved systems, such as HME, the invariances of Boltzmann equation, such as Galilean invariance, scaling invariance, reflecting invariance, etc., are naturally preserved, which will be discussed in detail in the next section.

\subsection{Invariances of moment closure models}
\label{sec:mom_inv}
As is introduced in Sec. \ref{sec:invar}, the Boltzmann equation itself has some invariances, such as Galilean invariance, reflecting invariance, scaling invariance. These invariances are the reflection of some basic physical laws. As a reduced model of the Boltzmann equation, the moment equation is expected to maintain these invariances. We will discuss the invariances of the moment systems and what constraints should be imposed on the closure relation to keep these invariances. 

Assuming that the closure relation can be represented by the following operator $\mG[\brho]$ as 
\begin{equation}
\label{eq:network}
   f_{M+1}(t,\cdot)  = \mG[\brho] = \mG[\rho(t,\cdot), u(t,\cdot), \theta(t,\cdot), f_3(t,\cdot), \cdots, f_M(t,\cdot), \Kn],
\end{equation}
then, the derivatives of the closure relation can be derived as the derivative of $\mG[\brho]$
\begin{equation}
\label{eq:network_der}
     \frac{\partial f_{M+1}(t,\cdot)}{\partial x} = \pd{\mG[\brho]}{x} \triangleq \mG_x [\rho(t,\cdot), u(t,\cdot), \theta(t,\cdot), f_3(t,\cdot), \cdots, f_M(t,\cdot), \Kn].
\end{equation}
In the absence of ambiguity, we abbreviate them as
\begin{equation}
\label{eq:short_net_work}
\begin{aligned}
    f_{M+1} & = \mG[\brho] = \mG[\rho, u, \theta, f_3, \cdots, f_M, \Kn], \\
        \frac{\partial f_{M+1}}{\partial x}  & = \mG_x[\brho] = \mG_x[\rho, u, \theta, f_3, \cdots, f_M, \Kn].
\end{aligned}
\end{equation}
% We abbreviate the Boltzmann equation as $\mathcal{B}(f)=0$, where $f=f(t,x,v)$ is the distribution function, and we note that the space in which $f$ lies is $\mathbb{V}$. Call the Boltzmann equation invariant under some transformation $\mathcal{F}: \mathbb{V}\rightarrow \mathbb{V}$, if for any $f$ satisfying $\mathcal{B}(f)=0$, we have $\mathcal{B}(\mathcal{F}(f))=0$.

\paragraph{Galilean invariance}
The Galilean invariance says that classical mechanics has the same form in all inertial systems. For the Boltzmann equation, Galilean invariance means that the form of the equation remains invariant under Galilean transformations. We can split the Galilean transformation into translation invariance in time, rotation and translation invariance in space, and invariance under Galilean boost. Specifically, for the original Boltzmann equation, Galilean invariance means that it satisfies the property \eqref{eq:Galilean}.  

% for every $t_0,x_0,u_0 \in \mathbb{R}$, define $f^*(t,x,v)= \mathcal{F}_{\rm Gal}(f;u_0) = f(t-t_0,x-x_0-u_0 t, v-u_0)$. 
% If $f$ is a solution of the Boltzmann equation, then so is $f^*$.

For the moment system \eqref{eq:linear_system}, several constraints should be added to the closure term to preserve the Galilean invariance. As is stated in \cite{han2019uniformly}, the moment coefficients defined in \eqref{eq:coef}  preserves the Galilean invariance. So we need to make the closure relation maintain Galilean invariance to make the closed moment system Galilean invariant. This means that the closure relation should not vary with the choice of the inertial reference systems, which could be achieved by letting the closure relation satisfy the following constraints for arbitrary constant $u_0 \in \mathbb{R}$:
\begin{equation}
\begin{aligned}
    \label{eq:moment_gal_inv}
    & \mG[\rho, u, \theta, f_3, \cdots, f_M, \Kn]=\mG[\rho, u-u_0, \theta, f_3, \cdots, f_M, \Kn], \\
    & \mG_x[\rho, u, \theta, f_3, \cdots, f_M, \Kn]=\mG_x[\rho, u-u_0, \theta, f_3, \cdots, f_M, \Kn].
\end{aligned}
\end{equation}

\paragraph{Reflecting invariance}
For the reflecting invariance, it means that the physical law remains unchanged under the reflecting transformation. For the moment system, as long as the closure relations are mirror-symmetric, the closed moment system conserves the reflecting invariance. From the definition of the moment coefficients, we notice that when the direction of the coordinate axis is flipped, the moments of the even order remain unchanged, while the sign of the odd-order moments changes. Therefore, the closed moment system has reflecting invariance, if the closure relation remains unchanged under the transformation below \cite{Porteous2021maxent1}. Defining $\mathcal{P}[g](t, x) = g(t, -x)$, then $\mG[\cdot]$ satisfies 
{\small 
\begin{equation}
\label{eq:moment_gal_inv_ref}
\begin{aligned}
& \mG[\rho, u, \theta, f_3, \cdots, f_M, \Kn](t,x) =(-1)^{(M+1)} \mG[\mathcal{P}\rho, -\mathcal{P}u, \mathcal{P}\theta, (-1)^3\mathcal{P}f_3, \cdots, (-1)^{M}\mathcal{P}f_M, \Kn](t, -x), \\
& \mG_x[\rho, u, \theta, f_3, \cdots, f_M, \Kn](t, x) = (-1)^{M} \mG_x[\mathcal{P}\rho, -\mathcal{P}u, \mathcal{P}\theta, (-1)^3\mathcal{P}f_3, \cdots, (-1)^{M}\mathcal{P}f_M, \Kn](t, -x).
\end{aligned}
\end{equation}
}

\paragraph{Scaling invariance}
Scaling invariance corresponds to the dimensional correctness of physical laws. For the moment system, if the closure relation is dimensionless, then naturally the scaling invariance is guaranteed. Moreover, if the constraints below are satisfied, we can also say that the closed moment system conserves scaling invariance. 
\begin{align}
\lambda \mG[\rho, u, \theta, f_3, \cdots, f_M , \Kn]
&=\mG[\lambda \rho, u, \theta, \lambda f_3, \cdots, \lambda f_M , \Kn], \label{eq:scale-rho} \\
\mu^{M+1} \mG[\rho, u, \theta, f_3, \cdots, f_M , \Kn]
&=\mG[\rho, \mu u, \mu^2 \theta,\mu^3 f_3, \cdots, \mu^M f_M ,\mu^{-1} \Kn].  \label{eq:scale-veloc}
\end{align}
The constraints for $\mG_x[\cdot]$ can be deduced based on \eqref{eq:scale-rho} and \eqref{eq:scale-veloc}. 

% As a law describing a physical phenomenon, its form should be independent of the specific unit of measurement used. The property that the equations remain form invariant when we change the unit of measurement is known as scale invariance.

% However, it is very easy to neglect to maintain scale invariance when using neural networks for moment closure. In our work, the dimensionless treatment is done specifically to maintain the scale invariance.

For now, we have listed the constraints on the closure relation to preserve the invariances. For the Grad's moment system, $f_{M+1}$ is simply set as $ f_{M+1} = 0,$
% \begin{equation}
%     \label{eq:closure_grad}
%     f_{M+1} = 0,
% \end{equation}
which is easy to verify that this closure relation satisfies all the constraints mentioned above, and the Grad's moment system conserves the invariances. For the globally hyperbolic moment system, the regularized term for $f_{M+1}$ \cite{Fan_new} is 
\begin{equation}
    \label{eq:NRxx} 
    \pd{f_{M+1}}{x} = f_M \pd{u}{x} + f_{M-1}\pd{\theta}{x}.
\end{equation} 
We can also easily check that the closure relation \eqref{eq:NRxx} satisfies the constraints imposed on the moment system. Precisely, it holds that 
\begin{equation}
\label{eq:NRxx-inv1} 
\begin{aligned}
\pd{f_{M+1}}{x} &= f_M \pd{u}{x} + f_{M-1}\pd{\theta}{x} = f_M \pd{(u-u_0)}{x} + f_{M-1}\pd{\theta}{x}, \qquad \forall u_0 \in \bbR.
\end{aligned}
\end{equation}
Then the closure equation \eqref{eq:NRxx} satisfies the Galilean invariance. Moreover, we can derive that 
\begin{equation}
\label{eq:NRxx-inv2} 
\begin{aligned}
\pd{f_{M+1}}{x} &
= (-1)^M \pd{(-1)^{M+1}f_{M+1}}{(-x)} 
= (-1)^M \left((-1)^{M}f_M \pd{(-u)}{(-x)} + (-1)^{M-1} f_{M-1}\pd{\theta}{(-x)}\right), \\
\pd{f_{M+1}}{x}
& = \frac{1}{\lambda \mu^{M+1} }\pd{ \lambda \mu^{M+1} f_{M+1}}{x} = \frac{1}{\lambda \mu^{M+1} } \left(\lambda \mu^M f_M \pd{\mu u}{x} + \lambda \mu^{M-1} f_{M-1}\pd{\mu^2 \theta}{x}\right),
\end{aligned}
\end{equation}
which means that the closure equation \eqref{eq:NRxx} preserves the reflection and scaling invariance. All together, we can deduce that the globally hyperbolic moment method proposed in \cite{Fan_new} conserves these invariances.

For the classical numerical methods to obtain the closed moment system, we can always derive the closing term's specific expression, making it easy to check whether the closed system keeps the invariances. However, when learning the moment closure model by a neural network, we can not get the explicit expressions for the closing term. Therefore, the neural network should be specially designed to ensure that the closed moment system presumes these invariances, which we will discuss in detail in the next section. 

%%% Local Variables: 
%%% mode: latex
%%% TeX-master: "article"
%%% End: 

%% file: article_scheme.tex
\section{Invariance preserving neural closure}
\label{sec:sym_net}

Learning-based moment closure is a recently developed method in kinetic theory. In \cite{han2019uniformly, huang2021learning}, the machine learning methods are all successfully utilized in the moment closure of the kinetic equation. The core idea is to obtain the closure relation \eqref{eq:short_net_work} for $f_{M+1}$ using a neural network based on the known moments. 
%Learning-based moment closure is a recently developed method in kinetic theory. In \cite{han2019uniformly, huang2021learning}, the machine learning methods are all successfully utilized in the moment closure of the kinetic equation. The core idea is to get the closure function \eqref{eq:short_net_work} for $f_{M+1}$ using neural network based on the known moments. 

 Unlike the traditional moment closure method, it is always difficult for the closure relation derived by the machine learning method to express the closure term explicitly. Thus, it is not easy to verify whether it has physical invariances. Special designs should be made on the structure of the neural network to make this possible. For example, Galilean invariance is preserved in \cite{han2019uniformly, huang2021learning}. The scaling invariance and reflecting invariance are discussed in \cite{huang2021machine1, Porteous2021maxent1, Schotthofer2021maxent2} and \cite{Porteous2021maxent1} respectively. In this section, we will describe how to design our neural network to incorporate all three invariances. We will call the specially designed neural network the invariance preserving neural closure (IPNC) network. 

\begin{figure}
	\centering
	\subfloat[]{
		\begin{minipage}[c]{0.45\textwidth}
			\includegraphics[width=1\textwidth]{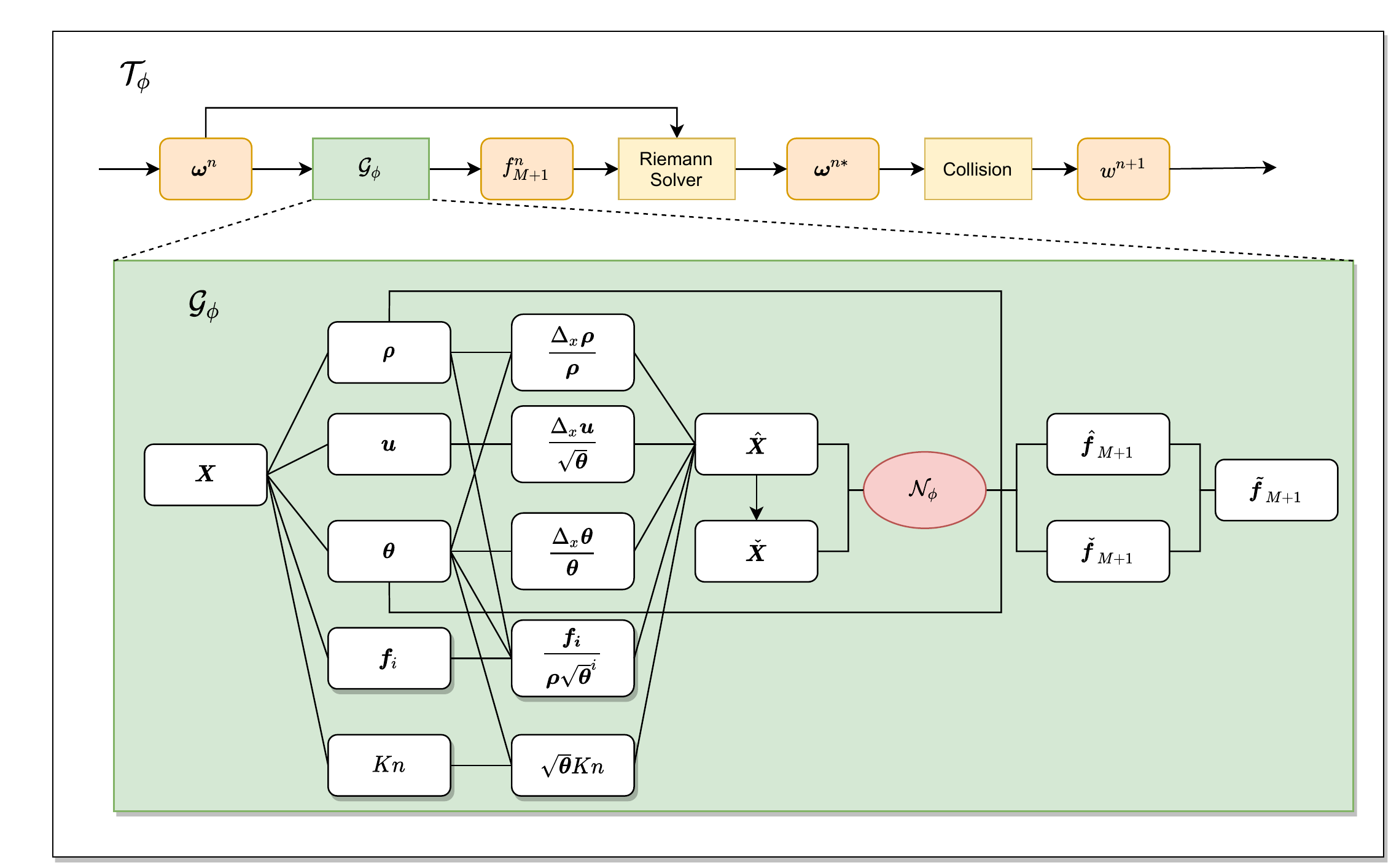}
		\end{minipage}
		\label{fig:Tau}
	}
    \subfloat[]{
    		\begin{minipage}[c]{0.45\textwidth}
   		 	\includegraphics[width=1\textwidth]{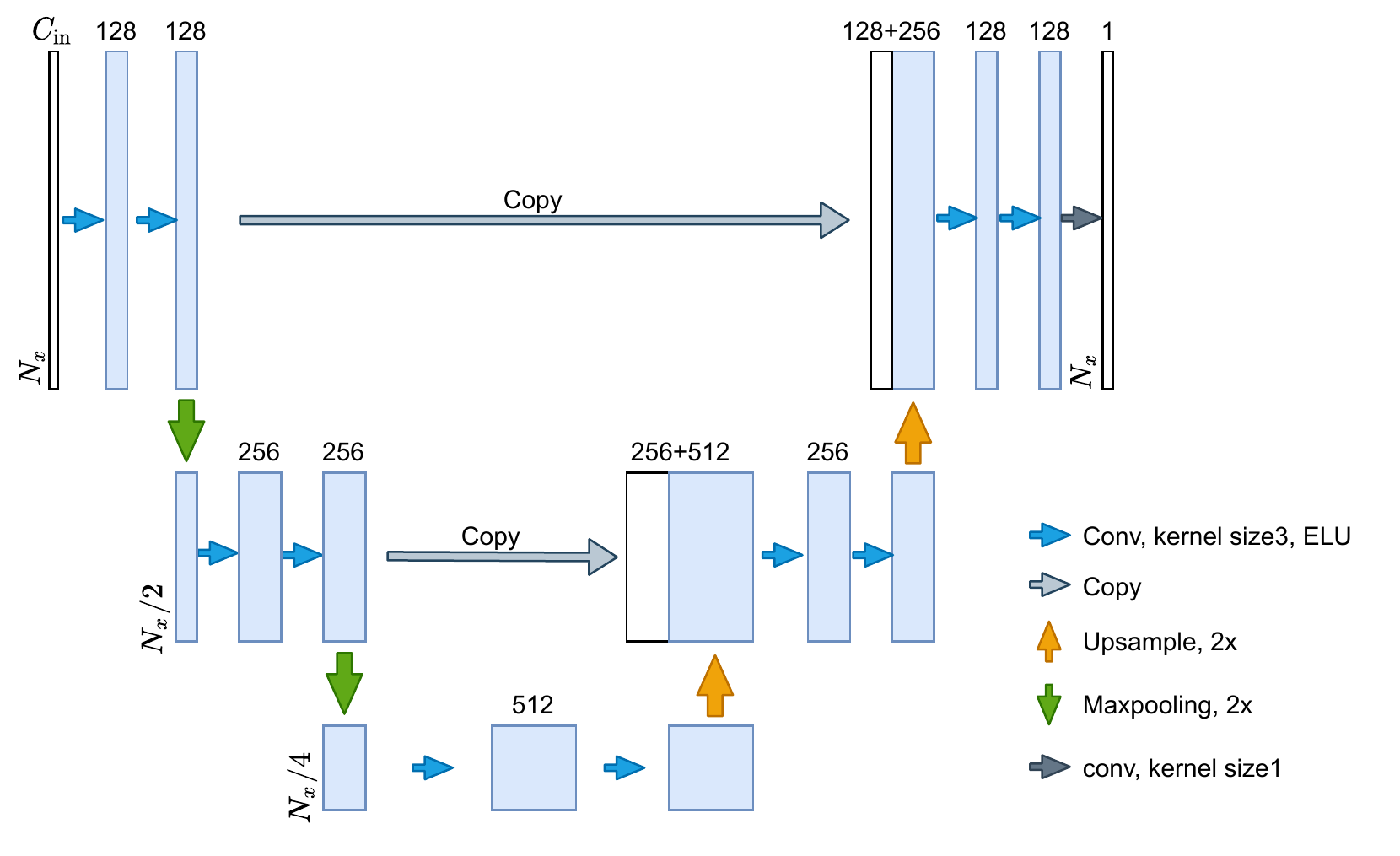}
    		\end{minipage}
		\label{fig:u-net}
    	}
\caption{(a) (Sec. \ref{sec:sym_net})  Schematic diagram of the time evolution operator $\mathcal{T}_{\param}[\cdot]$ and the main structure for the closure net.
(b) (Sec. \ref{sec:arc_neu}) The U-Net architecture. Each blue box corresponds to a multi-channel feature map while the white box represents the copied feature map. The number of channels is at the top of the box and the width of the feature map is at the bottom left corner of the box.}
\end{figure}

First, we will represent the time evolution operator in one time step as 
\begin{equation}
    \label{eq:one_step} 
\boldsymbol{\omega}^{i+1} :=\mathcal{T}_{\phi}[\boldsymbol{\omega}^i;\Kn,\Delta t]
% :=\left\{  
% \begin{array}{l}
% f_{M+1}^n = \mathcal{G}_{\phi}[\bomega^n;\Kn],     \\
% \mathcal{S}[\bomega^n,f_{M+1}^n;\Kn,\Delta t],
% \end{array}
% \right.
\end{equation}
where the operator $\mathcal{T}_{\param}[\cdot]$ contains two parts. The first part is the closure relation represented by a operator $\mathcal{G_{\param}}[\cdot]$ which is obtained by the neural network as
\begin{equation}
    \label{eq:T_first}
    f_{M+1}^i = \mathcal{G}_{\phi}[\bomega^i;\Kn].
\end{equation}
Here $\phi$ represents the parameters in the neural network, and the detailed property of $\mathcal{G_{\param}}[\cdot]$ will be introduced in Sec. \ref{sec:arc_neu}. The second part is the numerical scheme $\mathcal{S}[\bomega,f_{M+1};\Delta t]$ to update the numerical solution as
\begin{equation}
    \label{eq:T_second}
   \boldsymbol{\omega}^{i+1} =  \mathcal{S}[\bomega^i,f_{M+1}^i;\Kn,\Delta t].
   \end{equation}
The exact structure of it is shown in Figure \ref{fig:Tau}. For the classical numerical method, such as HME, there is also a similar time evolution operator. However, the main difference is that in HME the closure relation is given by \eqref{eq:NRxx} while the closure operator $\mathcal{G}_{\phi}[\cdot]$ in $\mT_{\phi}[\cdot]$ is obtained from learning. Here, we will focus on how to preserve the invariances. Special design should be applied to the neural network, and the methodology to choose the input and output variables should also be carefully studied.

% \begin{figure}[!htb]
% \centering
% \includegraphics[width=0.4\textwidth]{images/U-Net.pdf}
% \caption{(Sec. \ref{sec:arc_neu}) The U-Net architecture. Each blue box corresponds to a multi-channel feature map while the white box represents the copied feature map. The number of channels is at the top of the box and the width of the feature map is at the bottom left corner of the box.}
% \label{fig:u-net}
% \end{figure}

% \begin{figure}[!htb]
% \centering
% \includegraphics[width=0.4\textwidth]{images/Tau.pdf}
% \caption{(Sec. \ref{sec:sym_net})  Schematic diagram of the time evolution operator $\mathcal{T}_{\param}[\cdot]$ and the main structure for the closure net.} 
% \label{fig:Tau}
% \end{figure}

% \begin{figure}[!htb]
% \centering
% \includegraphics[width=0.4\textwidth]{images/Tau.pdf}
% \caption{(Sec. \ref{sec:sym_net})  Schematic diagram of the time evolution operator $\mathcal{T}_{\param}[\cdot]$ and the main structure for the closure net.} 
% \label{fig:Tau}
% \end{figure}

\subsection{Invariance preservation}
\label{sec:net_inv}

\paragraph{Galilean Invariance}
As is stated in Sec. \ref{sec:mom_inv}, the classical Grad moment equation \eqref{eq:linear_system} conserves Galilean invariance since the moment coefficients are defined based on the local velocity. Therefore, as long as the closure relation for $f_{M+1}$ satisfies Galilean invariance, the closed moment model preserves Galilean invariance.

The Galilean invariance includes translational, rotation, and motion invariance. To maintain translational invariance, we will use the neural network which naturally maintains translational invariance, such as a fully convolutional neural network. As this is an initial work for us to derive the closure relation by the neural network, we only consider the problem with 1D spatial space and 1D microscopic velocity space. Thus, we will not consider the rotation invariance for the 
closure. 

As motion invariance, if the closure relation is independent of the macroscopic velocity $u$, the Galilean invariance will be satisfied \cite{han2019uniformly, huang2021learning}. However, as one of the most important physical variables, we should include information on the macroscopic velocity when learning the closure relation by a neural network. Noticing that the spatial derivatives of the macroscopic variables such as $u_x$ and $u_{xx}$ all satisfy the motion invariance, these derivatives are included in the neural network instead of the macroscopic velocity $u$. 

\paragraph{Reflecting invariance}
As is discussed in Sec. \ref{sec:mom_inv}, as long as the closure relation for $f_{M+1}$ preserves the reflecting invariance \eqref{eq:moment_gal_inv}, the closed moment system will inherit the reflecting invariance.  
 
However, it is difficult to restrict the parameters of the neural network directly so that it gives a closure relation with reflecting invariance. Instead, we can obtain a closure relation that satisfies reflecting invariance by adjusting the input and output of the neural network. In our simulation, the testing data and its flip in the spatial direction are all entered into the neural network. The output is derived by combining the original output and the output flipped in the spatial direction.
This process is done by defining another neural network for the flipped moment coefficients. Noting that for the moment with odd order, flipping the spatial coordinates will lead to a change in the sign of the moment coefficients. Therefore, the closure relation for the flipping moment coefficients is defined as  
 \begin{equation}
     \label{eq:flip_network}
     \tilde{\mG}[\bomega](t, x) = (-1)^{M+1} \mG[\mathcal{P}\rho, -\mathcal{P}u, \mathcal{P}\theta, (-1)^3\mathcal{P}f_3, \cdots, (-1)^M\mathcal{P} f_M , \Kn](t, -x).
 \end{equation}
 Then the final closure relation to preserve reflecting invariance is derived by combing these two closure functions together 
\begin{equation}
\label{eq:reflect}     
\mG^{\rm sym}[\bomega](t, x)
= \frac{\mG[\bomega](t, x)+
\tilde{\mG}[\bomega](t, x)}{2}. 
\end{equation}
It is easy to verify that the closure relation  \eqref{eq:reflect}  with the specially designed input and output will conserve the reflecting invariance. 

\paragraph{Scaling invariance}
As we mentioned in Sec. \ref{sec:mom_inv}, if the closure relation for $f_{M+1}$ satisfies \eqref{eq:scale-rho} and \eqref{eq:scale-veloc}, then the closed moment system will preserve the scaling invariance.

To satisfy this, the closure network is designed in the following steps. First, the density $\rho$ and the temperature $\theta$ is scaled to $1$ by setting 
\begin{equation}
    \label{eq:scale_net} 
    \lambda = \rho^{-1}, \qquad \mu = \theta^{-1/2},
\end{equation}
in \eqref{eq:scale-rho} and \eqref{eq:scale-veloc}. As is stated that to preserve Galilean invariance, the macroscopic velocity $u$ is replaced by the derivatives such as $u_x$, the closure relation \eqref{eq:short_net_work} is changed into  
\begin{equation}
\label{eq:scale-form-1}     
\mG[\rho, u, \theta, f_3, \cdots, f_M , \Kn]
=
\rho\theta^{(M+1)/2} \mG \left[1, \frac{u_x}{\sqrt{\theta}}, 1, \frac{f_3}{\rho\theta^{3/2}}, \cdots,\frac{f_M}{\rho\theta^{M/2}},\sqrt{\theta} \Kn\right].
\end{equation}
However, as two of the most important physical variables, the information of the density $\rho$ and temperature $\theta$  should be included in the neural network. A similar method as preserving Galilean invariance is utilized here. The information of the derivatives of the density $\rho_x$ and temperature $\theta_x$ are included in the neural network as 
\begin{equation}
\label{eq:scale-form-2}     
\mG[\rho, u, \theta, f_3, \cdots, f_M , \Kn]
=
\rho\theta^{(M+1)/2} \mG \left[\frac{\rho_x}{\rho}, \frac{u_x}{\sqrt{\theta}}, \frac{\theta_x}{\theta}, \frac{f_3}{\rho\theta^{3/2}}, \cdots,\frac{f_M}{\rho\theta^{M/2}},\sqrt{\theta} \Kn\right].
\end{equation}

Finally, since we only have the information of the discretized difference of these macroscopic variables, the difference instead of the continuous functions of the derivatives is included as the input of the neural network. Consequently, the final form of the closure relation for the neural network is 
\begin{equation}
\label{eq:scale-form}     
\mG[\rho, u, \theta, f_3, \cdots, f_M , \Kn]
=
\rho\theta^{(M+1)/2} \mathcal{N}\left[\frac{\Delta_x \rho}{\rho}, \frac{\Delta_x u}{ \sqrt{\theta}}, \frac{ \Delta_x \theta}{\theta}, \frac{f_3}{\rho\theta^{3/2}}, \cdots,\frac{f_M}{\rho\theta^{M/2}},\sqrt{\theta} \Kn\right], 
\end{equation}
where $\mathcal{N}$ denotes an arbitrary convolutional neural network to maintain transitional invariance. 

\begin{remark}
With the special design of the neural network \eqref{eq:scale-form}, the input data is normalized by the density and temperature, with the density and temperature normalized to $1$. Therefore, if we encounter the problem where the magnitude of the physical variables for the testing data is much larger than those used in the training network, we can still derive a reasonable closure relation with this normalized neural network \eqref{eq:scale-form}. Moreover, with the scaled input data, the magnitude of the training data for different problems are similar, making the network more robust. 
\end{remark}

For now, we have finished the introduction of the special constraint on the neural network to preserve these invariances. We want to mention here that the invariance preserving neural network's design only requests the corresponding processing of the inputs and outputs and does not depend on the specific neural network structure, which grants vast flexibility in choosing neural networks. 

\subsection{Network architecture}
\label{sec:arc_neu}
% As is stated in Sec. \ref{sec:net_inv}, to preserve Galilean invariance, the neural network with transnational invariance is utilized. 
In this section, we present the network structure in detail. The closure operator $\mathcal{G}_{\phi}[\cdot]$ for the proposed IPNC takes the following form, 
\begin{equation}
    \label{eq:operator_Q}
    \tilde{\bm{f}}_{M+1} =  \mathcal{G}_{\phi}[{\bm X}]= \frac{1}{2}\left[\widetilde{\mathcal{R}}\mathcal{N}_{\phi}(\mathcal{R} {\bm X})+
    \widetilde{\mathcal{P}}\widetilde{\mathcal{R}}\mathcal{N}_{\phi}(\mathcal{R}\mathcal{P} {\bm X})\right],
\end{equation}
where the input $\bm X$ is a second-order $N_x \times (M+2)$ tensor with $N_x$ corresponding to the spatial discretization and $M+2$ corresponding to the concatenation of the first $M-$th order moment coefficients and $\Kn$. For each row of $\bm X$, its entries are $\bm{X}_i = [\bomega_i; \Kn]$ with $\bomega_i = [\rho,u,\theta,f_3, ... . f_{M}]_i, i = 1, \cdots, N_x$. The operator $\mathcal{P}$ is the parity transformation proposed in \eqref{eq:flip_network} to preserve reflecting invariances. To be more precise, $\mathcal{P} {\bm X}$  is another second-order $N_x \times (M+2)$ tensor and its $i$-th row is 
\begin{equation}
	\label{eq:row_Px}
		( \mathcal{P}{\bm X})_i 
		 = \left[{\rho}, -u, {\theta}, (-1)^{3}{f_3}, \cdots,(-1)^{M}{f_M}, \Kn\right]_{N_x-i}.
\end{equation}
The operator $\mathcal{R}$ is the transformation proposed in \eqref{eq:scale-form} to preserve Galilean and scaling invariances. $\mathcal{R} {\bm X}$  is also a second-order $N_x \times (M+2)$ tensor whose $i$-th row is 
% When fed into the network, it first applies the transformation $\mathcal{R}_{\rho,u,\theta}$ to obtain $\tilde X$, which is still a $N_x\times(M+2)$ second-order tensor,
\begin{equation}
    \label{eq:row_Rx}
   ( \mathcal{R}{\bm X})_i 
    = \left[\frac{\Delta_x \rho}{\rho}, \frac{\Delta_x u}{ \sqrt{\theta}}, \frac{ \Delta_x \theta}{\theta}, \frac{f_3}{\rho\theta^{3/2}}, \cdots,\frac{f_M}{\rho\theta^{M/2}},\sqrt{\theta} \Kn\right]_i,
    \end{equation}
and $\mathcal{R}\mathcal{P}{\bm X} = \mathcal{R}(\mathcal{P}{\bm X})$ has a similar form.
% Then $\tilde X$ will be used as input to the neural network
% $$
% \tilde f_{M+1}=\mathcal{N}_{\phi}(\tilde X)
% $$
% where 

In \eqref{eq:operator_Q},  $\mathcal{N}_{\phi}$ is a generic neural network which can be a multilayer perceptron (MLP) or a convolutional neural network. We will describe the exact structure of $\mathcal{N}_{\phi}$ in the next sections. It takes the $N_x \times (M+2)$ tensor $\mathcal{R} {\bm X}$ or $\mathcal{R} \mathcal{P}{\bm X}$ as input and a vector with length $N_x$ as output. The final output $\tilde{\bm{f}}_{M+1}$ of the operator $\mathcal{G}_{\phi}[\cdot]$ is also a vector with length $N_x$, which is obtained by the transformation $\widetilde{\mathcal{R}}$ and the parity transformation$\widetilde{\mathcal{P}}$ as 
% \begin{equation}
%     \label{eq:out_G}
% \tilde{\bm{f}}_{M+1} =\frac{1}{2}\left(\widetilde{\mathcal{R}}\hat{\bm{f}}+\widetilde{\mathcal{P}}\widetilde{\mathcal{R}}\check{\bm{f}}\right)
% \end{equation}
% with 
\begin{equation}
	\begin{aligned}
    \label{eq:tilde_R_P}
& (\widetilde{\mathcal{R}} \hat{\bm{f}})_i = \rho_i\theta_i^{(M+1)/2}\hat{f}_i\quad\mbox{and} \quad
 (\widetilde{\mathcal{P} }\widetilde{R} \check{\bm{f}})_i =  (-1)^{M+1}(\widetilde{R} \check{\bm{f}})_{N_x-i}, \quad i = 1, \cdots, N_x,
	\end{aligned}
\end{equation}
where $\hat{\bm{f}}$ and $\check{\bm{f}}$ are the output of $\mathcal{N}_{\phi}(\mathcal{R}[\cdot])$ and $\mathcal{N}_{\phi}(\mathcal{R}\mathcal{P}[\cdot])$ respectively as
\begin{equation}
    \label{eq:out_Nx}
\hat{\bm{f}} = \mathcal{N}_{\phi}(\mathcal{R} {\bm X}) \quad\mbox{and}\quad \check{\bm{f}} = \mathcal{N}_{\phi}(\mathcal{R} \mathcal{P} {\bm X}).
\end{equation}

There are several methods to design the neural network $\mathcal{N}_{\phi}$, which will lead to different properties of the network. Here, to keep the translation invariance in the closed moment system, we require that the structure of the neural network preserves the translation invariance, such as the convolutional networks. In this work, the shared MLP \cite{qi2017pointnet} and the U-Net \cite{ronneberger2015u} which is a particular convolutional neural network are utilized. 

\paragraph{Shared MLP}
\label{sec:full_connect}
In this work, the MLP is shared in the sense that only one MLP is learned for all spatial points. In this network, the information at each spatial point is put into the same MLP.
The network takes the moments of each order at each spatial position as input and the predicted moment coefficient at $(M+1)$-th order at the corresponding spatial position as output. Since only local information is utilized to predict the moment coefficient at $(M+1)$-th order, it is spatially translation invariant. 
For the MLP, its structure is mainly a consecutive composition of linear transformations and non-linear activation functions. In the numerical test for the shock structure problem in Sec. \ref{sec:num_shock}, a MLP with $12$ hidden layers and $128$ neurons per layer is utilized. To prevent gradient vanishing, the skip connection technique in ResNet \cite{he2016resnet} is adopted, and ReLU \cite{nair2010relu} is chosen as the activation function $\sigma$. 

In a MLP, the linear transformation within a single layer can be expressed as
\begin{equation}
Y_{k}=\sum_{i=1}^{d}W_{i,k}X_{i}+b_{k},
\end{equation}
while in a shared MLP, the linear transformation within a single layer can be expressed as
\begin{equation}
\label{eq:shared_MLP}
Y_{j,k}=\sum_{i=1}^{d}W_{i,k}X_{j,i}+b_{k},
\end{equation}
which we abbreviate to 
\begin{equation}
R^{0}(\boldsymbol{X}) \rightarrow \boldsymbol{Y}:\boldsymbol{Y}= {\bf W}^{(0)} \boldsymbol{X} + \boldsymbol{b}^{(0)},
\end{equation}
where $\boldsymbol{X}\in \mathbb{R}^{N_x,C_{\rm in}}$ is the input data, $\boldsymbol{Y}\in \mathbb{R}^{N_x,C_{\rm out}}$ is the output data, ${\bf W}\in \mathbb{R}^{C_{\rm in},C_{\rm out}}$  and $\boldsymbol{b}\in \mathbb{R}^{C_{\rm out}}$ are the parameters trained by the network. This MLP is shared between different spatial points with $j$ as index.

For the shared MLP with $L$ layers utilized in Sec. \ref{sec:num_shock}, it could be represented in the following form
\begin{equation}
\boldsymbol{Y}= {\bf W}^{(L+1)} \sigma \circ ({\bf I}+R^{(L)}) \cdots \sigma \circ ({\bf I}+R^{(1)})(W^{(0)}\boldsymbol{X}_{\rm in}+b^{(0)}) + \boldsymbol{b}^{(L+1)}, \quad l=1,\cdots,L
\end{equation}
where $\sigma$ denotes the activation function, ${\bf I}$ denotes the identity mapping, $\boldsymbol{X}_{\rm in}$ represents the input and $\boldsymbol{Y}$ represents the output. In this work,  the input data is a second-order tensor $\mathcal{R} \bm X$ or $\mathcal{P} \mathcal{R} \bm X$ with shape $N_x \times (M+2)$ and $\bm Y$ is the vector $\hat{\bm{f}}$ or $\check{\bm{f}}$ with length $N_x$ as in \eqref{eq:out_Nx}. The hidden layer $R^{(l)}$ has $128$ neurons, which means that $C_{\rm out}^{(l)}$ in $R^{(l)}$ is equal to $128$.

% neurons per layer means the number of channels $C_{\rm out}$ in each hidden layer $\mathcal{R}$ denotes the number of features of the input, which in our network takes the value $C_{\rm in}=M+2$ in general, where $M+1$ represents the moments of each order in $\bomega$ together with the Knudsen number $\Kn$.

\paragraph{Convolutional neural network.}
\label{sec:u-net}
In a convolutional neural network, due to the convolutional operation, the $(M+1)$-th order moment coefficient prediction at a spatial position in the output of the network depends not only on the moment coefficients at the local position, but also on the moment coefficients in the neighborhood. The following expression shows the convolution operation in such a network
\begin{equation}
\label{eq:con-neural}
    Y_{j,k}=\sum_{i=1}^{d}\sum_{r=-p}^{p}W_{r,i,k}\hat X_{j+r,i}+b_{k},
\end{equation}
where ${\bf W}$ is called the convolutional kernel and its shape is $[C_{\rm out}, 2p+1, C_{\rm in}]$, where $C_{\rm in}$ is still the number of features of the input,  $C_{\rm out}$ is that of the output, and $(2p+1)$ is the width of the convolutional kernel. $C_{\rm in}$ and $C_{\rm out}$ will also be referred to as the number of channels in the network.
Its value is trained by the network. $\boldsymbol{\hat{X}}$ is the result after padding the input data. For example, for the problem of periodic boundary conditions, the circular padding approach \eqref{eq:padding-cir} is adopted. 
\begin{equation}
    \label{eq:padding-cir}
    \hat X_j=\left\{
    \begin{array}{ll}
    X_{n-1+j}, & -p \leq j<0,\\
	X_j,  &0 \leq j < n,\\
	X_{j-n+1},  &n\leq j < n+p. 
     \end{array}
    \right.
\end{equation}

In the convolution neural network we use, the input $\boldsymbol{X}_{\rm in}=\mathcal{R}\bm X$ or $\mathcal{P} \mathcal{R}\bm X$ in the first layer is the same input as in the shared MLP, but each column is padded first and then input to the network according to \eqref{eq:padding-cir}.

Then what follows is a series of linear transformations and non-linear activation functions similar to those of the shared MLP. However, in each linear transformation, unlike the output $Y_{j,k}$ in a shared MLP which depends only on $X_{j,k}$ with $k\in\{0,1,\cdots,M+1\}$, $Y_{j,k}$ in a convolutional layer depends on $X_{l,k}$ with $l\in\{i-p,\cdots,i,\cdots,i+p\},k\in\{0,1,\cdots,M+1\}$.

% For the problem of constant boundary conditions, we use the constant padding approach,  

% \begin{equation}
%     \label{eq:padding-cir}
%     \hat X_i=\left\{
%     \begin{aligned}
% 	&X_{L} \qquad &-p& \leq i<0,\\
% 	&X_i  \qquad &0&\leq i < n,\\
% 	&X_{R} &n&\leq i < n+p,
%     \end{aligned}
%     \right.
% \end{equation}

The advantage of utilizing a convolutional neural network instead of a shared MLP is that the convolutional neural network can exploit the information of the spatial derivatives, which can be helpful in finding a better closure relation. For example, in the Fourier law for the Navier-Stokes equation, the closure relation \cite{Struchtrup2005} is 
\begin{equation}
    \label{eq:fourier_law}
    \sigma_{ij} = -2 \mu \pd{v_{\langle i }}{x_k \rangle }, \qquad q_i = -\kappa \pd{\theta}{x_i},
\end{equation}
where $\mu$ is the viscosity, and $\kappa$ is thermal conductivity. This NS closure requires the information of the spatial derivative of the macroscopic velocity and temperature, which cannot be represented by a shared MLP network. 
 
\begin{remark}
The shared MLP can be viewed as a convolutional neural network with all convolutional kernels having $1\times1$ support, which can be seen by comparing \eqref{eq:shared_MLP} and \eqref{eq:con-neural}.
\end{remark}

In the numerical experiments, the U-Net \cite{ronneberger2015u}, which is a fully convolutional neural network, is utilized as the main body in closure. Because of its multi-scale architecture, the U-Net is widely used in many problems in scientific computing, such as Navier–Stokes simulations \cite{thuerey2020deep}, and PDE-solvers \cite{hsieh2019learning, um2020solver}. This kind of neural network contains two main components: an encoder and a decoder. The encoder is a series of convolutions and pooling operations with decreasing resolution, and the decoder is a series of convolution and upsample operations with increasing resolution. In the encoder, there are intermediate outputs that are also directly adopted by the decoder. The output of the U-Net ends up with the same resolution as the input. The specific structure of the network adopted in this work is a U-Net using 1D convolution, which is shown in Figure \ref{fig:u-net}. 

% \begin{figure}[!htb]
% \centering
% \includegraphics[width=0.4\textwidth]{images/U-Net.pdf}
% \caption{(Sec. \ref{sec:arc_neu}) The U-Net architecture. Each blue box corresponds to a multi-channel feature map while the white box represents the copied feature map. The number of channels is at the top of the box and the width of the feature map is at the bottom left corner of the box.}
% \label{fig:u-net}
% \end{figure}

The U-Net utilized here contains three max-pooling operations and three upsample operations in addition to the convolution. The input is first convolved into 128 channel data, and then the number of channels is doubled after each pooling. In the decoder, the intermediate results from the encoder and the upsample data are concatenated together, then convolved with the number of channels gradually decreasing. Finally, the data are passed through a convolution layer to output a channel. The ELU activation function \cite{clevert2015ELU} is used here as
\begin{equation}
\sigma(x)=\operatorname{ELU}(x, \alpha)= 
\begin{cases}
x, & \text { if } x \geq 0, \\ 
\exp \left({x}\right)-1, & \text { otherwise. }
\end{cases}
\end{equation}

In the convolution operation, to keep the vector length constant, we need to perform a padding operation. For periodic boundary conditions, we use the corresponding circular padding \eqref{eq:padding-cir}, and for the fixed boundary conditions, we use the corresponding replicate padding. Padding operations do not contain parameters, which means that we can use the same network for different boundary conditions only by changing the type of padding. 

% \section{Detail design of the neural network}

\subsection{Training}
\label{sec:spe_app}
This section presents details on the training of the IPNC network introduced in the previous sections. The training data is generated by solving the Boltzmann equation with many different initial conditions using DVM. The training data set contains several trajectories, each of which is a numerical solution to the Boltzmann equation at different time steps for a certain initial value. These trajectories will be considered as the ground truth. 

When training the network, supervised learning is adopted. We train the network by minimizing a loss function. This loss function measures the difference between the model predicted by IPNC and the ground-truth solutions by DVM. Depending on how different loss functions are designed, we classify the training approach used in this work into end-to-end learning and direct learning.

\subsubsection{End-to-end neural closure learning}
\label{sec:end_end}

If we consider the numerical algorithm containing the neural network as a whole, and then calculate the distance between the trajectory predicted by IPNC and the trajectory of the ground truth as the loss function and optimize it, such a training method is called end-to-end training.

\begin{figure}
	\centering
	\subfloat[End-to-end learning]{
		\begin{minipage}[c]{0.6\textwidth}
	\includegraphics[height=5cm]{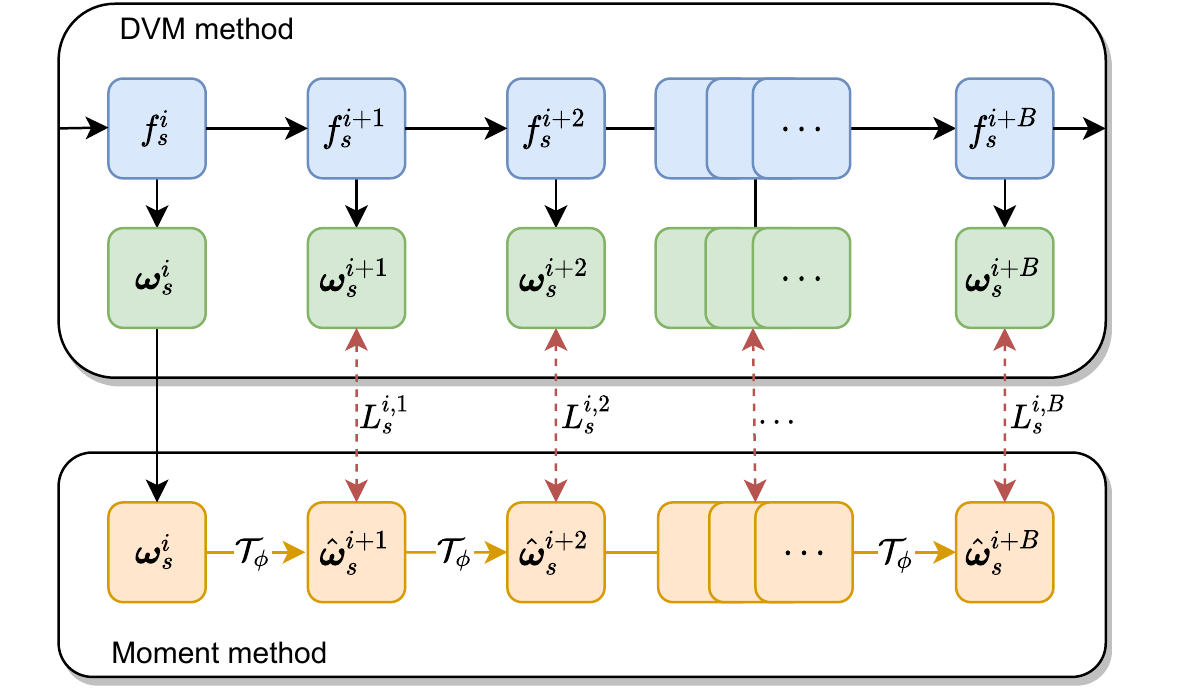}
		\end{minipage}
		\label{fig:pipeline}
	}
    \subfloat[Direct learning]{
	\begin{minipage}[c]{0.24\textwidth}
  	\includegraphics[height=5cm]{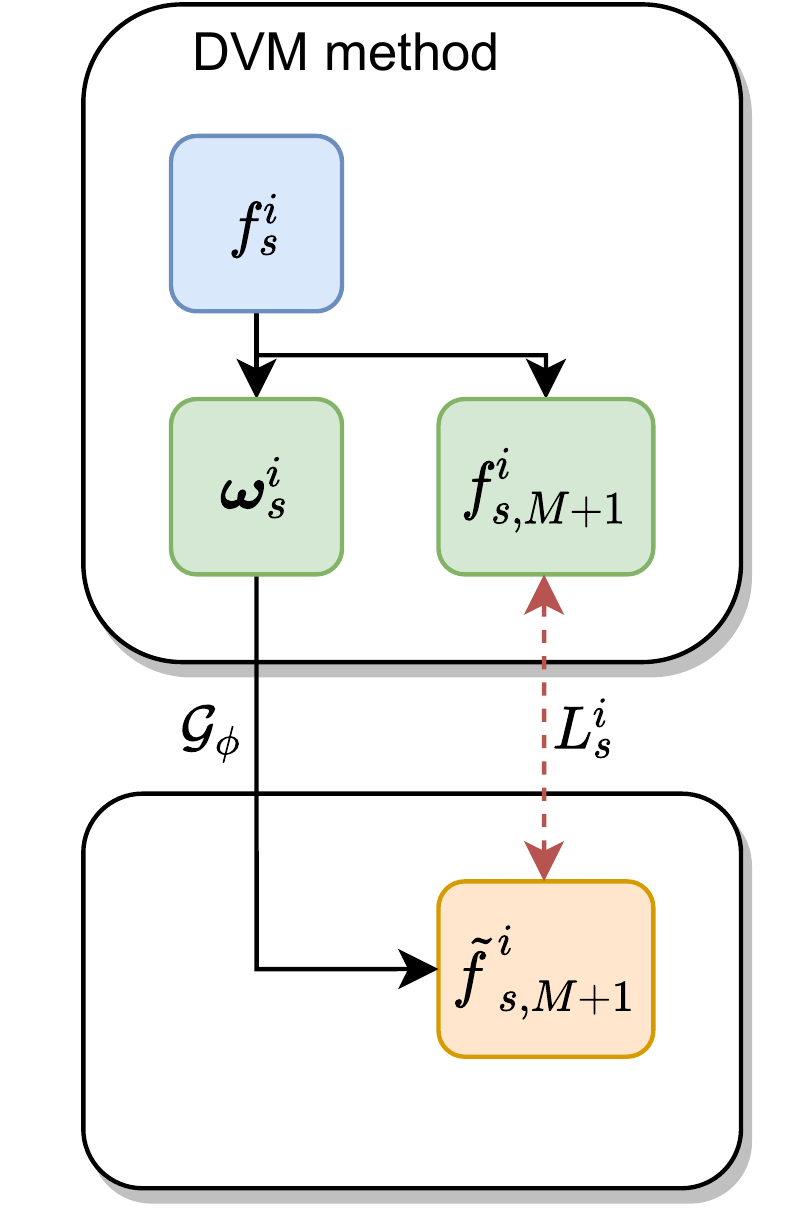}
    		\end{minipage}
		\label{fig:direct}
    	}
	\caption{(a) (Sec. \ref{sec:end_end}) The process of end-to-end closure learning. First, we select a segment of length $B$ from the DVM solution with the $s$-th initial condition and calculate the corresponding moments of each order to get $(\omega_s^{i+1}, \ldots, \omega_s^{i+B})$. Then the moment coefficients for the moment method are obtained by evolving the moment equation closed by IPNC with $\omega_s^i$ as the initial value, and finally, the closure network is trained by minimizing the distances $L_s^{i,b}, b=1,\ldots,B$, between the moment coefficients of two trajectories.
	(b) (Sec. \ref{sec:dir_learn}) The process of direct closure learning.  
    Here, $f_s^i$ is the solution selected at a certain time from the trajectories obtained from the DVM with the $s$-th initial condition and $\omega_s^i$ is calculated from $f_s^i$.
	Then a predicted $\tilde{f}^{i}_{s,M+1}$ is obtained by the closure network $\mathcal{G}_{\phi}$ using $\bomega_s^{i}$ in IPNC. Finally the closure network is trained by minimizing the distances $L_s^i$.}
	\label{fig:hor_2figs_1cap_2subcap}
\end{figure}

% \begin{figure}[!htb]
% \centering
% \includegraphics[width=0.4\textwidth]{images/E2E.pdf}
% \caption{(Sec. \ref{sec:end_end}) The process of end-to-end closure learning. First, a trajectory of the numerical solution is generated by solving the BGK equation by DVM, from which we obtain the moment coefficients $\bomega^{n}, n=0,1,\ldots,k$. Then another trajectory of moment coefficients is obtained by evolving the moment equation closed by IPNC, and finally, the closure network is trained by minimizing the distances $L^n, n=0,1,\ldots,k$, between the moment coefficients of two trajectories.}
% \label{fig:pipeline}
% \end{figure}

Figure \ref{fig:pipeline} shows the process of the end-to-end learning framework. Here, $\mathcal{T}_{\phi}$ is the PDE evolutionary operator, which is shown in Figure \ref{fig:Tau}. 
First, to construct our training data set, we will randomly generate $N_s$ initial values labeled as $\mathcal{S}=\{1,\ldots,N_s\}$ and evolve them to a certain time using DVM. 
Then, we will construct the training dataset in the form of $\{({\boldsymbol{w}}_s^{i+1},\ldots,{\boldsymbol{w}}^{i+B}_s)\}_{i \in \mathcal{I}, s\in \mathcal{S}}$, with $\boldsymbol{w}_s^i$. Here, $\boldsymbol{w}_s^i$ denotes the moment coefficients at $t^i$ solved using DVM with the $s$-th initial condition.  $\mathcal{I}$ is the set to choose the time index $i$, and $B$ is the length of the trajectory fragment. Each data point in our training dataset is a fragment containing $B$ time points taken from the trajectory obtained by the DVM. During the training, we use $\mathcal{T}_{\param}$ to evolve the numerical solution with ${\boldsymbol{w}}_s^i$ as the initial value to obtain the coefficients for the moment equations, which are symbolled as  $(\hat {\boldsymbol{w}}_s^{i+1},\ldots,\hat {\boldsymbol{w}}^{i+B}_s)$. The proceeding is as below, 
\begin{equation}
     \label{eq:hat_omega}
     \begin{aligned}
\hat{\boldsymbol{\omega}}_{s}^{i} = \boldsymbol{\omega}_{s}^i, \qquad 
     \hat{\boldsymbol{\omega}}_{s}^{i+1} = \mathcal{T}_{\param}[\hat{\boldsymbol{\omega}}_{s}^i], \qquad 
       \hat{\boldsymbol{\omega}}_{s}^{i+b} = \mathcal{T}_{\param}\cdots\mathcal{T}_{\param}  \left[\hat{\boldsymbol{\omega}}_{s}^i\right] \triangleq  \mathcal{T}^b_{\theta}[\hat{\boldsymbol{\omega}}_{s}^i] = \mathcal{T}^b_{\theta}[\boldsymbol{\omega}_{s}^i], \quad b = 1, \cdots, B.
     \end{aligned}
 \end{equation}
 Then, the distance between the short trajectory obtained by the DVM and the short trajectory obtained by $\mathcal{T}_{\phi}$ will be used as the loss function as
\begin{equation}
\label{eq:loss_new}
L(\phi)=\frac{1}{|\mathcal{S}||\mathcal{I}|}\sum_{s\in\mathcal{S}}\sum_{i\in\mathcal{I}}\sum_{b=1}^B\Vert 
\mathcal{T}_{\phi}^b(\boldsymbol{w}_s^{i} )-{\boldsymbol{w}}_s^{i+b} 
\Vert \triangleq \frac{1}{|\mathcal{S}||\mathcal{I}|}\sum_{s\in\mathcal{S}}\sum_{i\in\mathcal{I}}\sum_{b=1}^B L_s^{i,b},
\end{equation}
where $\Vert{\boldsymbol{w}} \Vert=\sum_{j=1}^{N_x}||\boldsymbol{w}(x_j)||_2^2$ and $L_s^{i,b} =\Vert 
\mathcal{T}_{\phi}^b(\boldsymbol{w}_s^{i} )-{\boldsymbol{w}}_s^{i+b} 
\Vert$.

Since IPNC is implemented using Pytorch \cite{pytorch}, where each step is differentiable, the whole process of the numerical solution can be derived using the automatic differentiation techniques. The gradient of the loss function for the network parameters is obtained and optimized using the stochastic gradient descent method AdamW \cite{AdamW}.

\begin{remark}
We also want to mention that the set $\mathcal{I}$ could be all the time steps we have simulated, but when there are a large number of time steps, only a subset of them is selected for training. For example, the index in $\mathcal{I}$ may be chosen randomly.
\end{remark}

\subsubsection{Direct neural closure learning}
\label{sec:dir_learn}

% \begin{figure}[!htb]
% \centering
% \includegraphics[width=0.4\textwidth]{images/Direct.pdf}
% \caption{(Sec. \ref{sec:dir_learn}) The process of direct closure learning. First, a trajectory of the numerical solution is generated by solving the BGK equation using DVM, from which we obtain the moment coefficients $\bomega^{n}$ and $f_{M+1}^{n}, n=0,1,\ldots, N_t$. Then a predicted $\tilde{f}^{n}_{M+1}$ is obtained by the closure network $\mathcal{G}_{\phi}$ using $\bomega^{n}$ in IPNC. Finally the closure network is trained by minimizing the distances between the $f_{M+1}^{n}$ and $\tilde f_{M+1}^{n}$.}
% \label{fig:direct}
% \end{figure}

Instead of end-to-end training, we can also use direct learning to obtain the closure relation. This means that we train this neural network to predict the moment coefficient at $(M+1)$-th order independently, rather than training it within a numerical solver for the moment system. 
Figure \ref{fig:direct} shows the process of the end-to-end learning framework. 
The loss function of direct learning is given as 
\begin{equation}
\label{eq:loss_dir}
L(\phi) = \sum_{s \in \mathcal{S}} \sum_{i \in \mathcal{I}} ||\tilde{f}_{ s,M+1}^i-f_{s,M+1}^{i}||_2^2= \sum_{s \in \mathcal{S}} \sum_{i \in \mathcal{I}} L^{i}_s,
\end{equation}
where $\tilde{f}_{s, M+1}^i=\mathcal{G}_{\phi}[\bomega^i_s]$ is the $(M+1)$-th order of moment coefficient predicted by the closure network at time $t = t^i$. $\mathcal{S}$ and $\mathcal{I}$ are the same as in end-to-end training. $\bomega^i_s$ which is derived by DVM, is the same as in the last section, and $f_{s, M+1}^i$ also obtained directly from DVM. 

%In this loss function, we define $\Vert{f} \Vert=\sum_{j=1}^{N_x}f(x_j)^2$.

\begin{remark}
From the loss function \eqref{eq:loss_dir}, we find that direct learning is simpler to implement and, in general, more time-efficient than end-to-end training. However, direct learning minimizes the average error of the closure itself, while end-to-end learning minimizes the average error of the solution. Since we want our numerical method to generate more accurate simulations, we would expect better performance from the end-to-end training than direct learning. We also note that the proposed IPNC method can derive the closed moment system for an arbitrary order of expansion. Moreover, the structure and training of the IPNC network can be flexibly adjusted according to the expansion order and can be incorporated with any numerical scheme solving the underlying moment system.
\end{remark}

Before ending this section, we present a technique we use to prepare the training data. We refer to this technique as channel-wise standardization. 

\paragraph{Channel-wise standardization} 
In the neural network, different channels correspond to different orders of moment coefficients. Since that the magnitude of the moments with different orders may have significant differences even after applying the scaling invariance technique \eqref{eq:scale-form}, inputting these moments directly into the neural network may be detrimental to the training of the neural network. Therefore, we adopt the channel-wise standardization on the training data \cite{bois2020neural}. For the neural network input of each channel, we calculate the mean and variance of each element of $\bm{X_{\rm in}}$ as 
\begin{equation}
    \label{eq:channel-std}
    \overline{\bm{X}}^{k} = \frac{1}{N_x N_t N_s} \sum_{j=1}^{N_x} \sum_{i=1}^{N_t}\sum_{s = 1}^{N_s} \bm{X}_{j,s}^{i,k},
     \quad  {\rm Var}(\bm{X}^k) = \frac{1}{N_x N_t N_s} \sum_{j=1}^{N_x} \sum_{i=1}^{N_t}\sum_{s = 1}^{N_s} (\bm{X}_{j,s}^{i,k} - \overline{\bm{X}}^k)^2, \quad k = 1, \cdots M+1,
\end{equation}
where $\bm{X}_{j,s}^{i,k}$ is the $k$th entry of $\bm{X}$  from the sample $s$ at the physical position $x_j$ and time step $i$. $N_x$ is the total grid number in the physical space, $N_t$ is the total time step and $N_s$ is the total number of samples. Then each $\bm{X}_{j,s}^{i,k}$ is normalized by the mean $\overline{\bm{X}}^k$ and the variance ${\rm Var}(\bm{X}^k)$. After the channel-wise standardization, the mean and variance for all the training data in different channels are close to $0$ and $1$, respectively. Moreover, the estimated mean and variance of each channel \eqref{eq:channel-std} are stored in the neural network along with the rest of the neural network parameters. For the testing data, they are normalized with the same mean and variance \eqref{eq:channel-std} calculated by the training data. With this channel-wise standardization, the training and testing data are normalized in the same way, and the order difference between different moment coefficients could be eliminated, which will make the network more robust.

%

%%% Local Variables: 
%%% mode: latex
%%% TeX-master: "article"
%%% End: 

%% file: article_experiment.tex
\section{Numerical experiments}
\label{sec:num}
In this section, we will adopt a similar numerical scheme as in HME \cite{NRxx} to solve the closed moment system by IPNC, which we will introduce briefly in Sec. \ref{sec:num_sch}. 
Then, several numerical examples are presented to validate the proposed IPNC method. First, in Sec. \ref{sec:num_smooth_dis} and \ref{sec:num_discon_dis}, we will test the performance of IPNC on the problems with smooth and discontinuous initial conditions, and several comparisons are also made with IPNC and other deep learning methods. In Sec. \ref{sec:num_sod}, the neural network trained by IPNC on the discontinuous initial condition in Sec. \ref{sec:num_discon_dis} is tested on the Sod shock tube problem without retraining. In Sec. \ref{sec:num_shock}, the shock structure problem is tested with IPNC. The ablation experiment to show the effect of the invariances in the numerical simulation is presented in Sec. \ref{sec:abl}.

Furthermore, we will also release a dataset based on the problem setting given in \cite{han2019uniformly} that can be used for benchmarking. We will also release the code used to generate the dataset. Our code supports 1D discrete velocity algorithms with several different temporal and spatial discrete schemes using GPU acceleration. The code of the proposed IPNC will be released after the paper is published.

\subsection{Numerical scheme for the moment equations} 
\label{sec:num_sch}
To solve the moment system, a numerical scheme similar to that of \cite{NRxx} is used. The standard finite-volume method is adopted in the $x$-direction. Suppose that $\Gamma_h$ is a uniform grid in $\bbR$ as
\begin{equation}
    \label{eq:cell}
    \Gamma_h = \left\{T_j = x_0 + \left(j \Delta x, (j+1)\Delta x\right): j \in \bbZ \right\},
\end{equation}
the numerical solution to approximate the distribution function $f$ at time $t=t_n$ and position $x_j$ is 
\begin{equation}
    \label{eq:dis_f}
    f_h^n(x,v) = f_j^n(v) = \sum_{0\leqslant \alpha \leqslant M} f_{j,\alpha}^n \mH_{\alpha}\left(\frac{v-u_j^n}{\sqrt{\theta_j^n}}\right), \qquad x \in T_j.
\end{equation}
The standard Strang's splitting method is utilized here, and solving the Boltzmann equation is split into the convection step and the collision step. For the convection step, the distribution function is updated by the standard finite volume method as 
\begin{equation}
    \label{eq:convection}
    f_j^{n+1, \ast}(v) = f_j^n(v) + K_{1,j}^n(v),
\end{equation}
where $K_{1,j}^n$ has the form 
\begin{equation}
    \label{eq:flux}
    K_{1,j}^n = -\frac{\Delta t^n}{\Delta x}\left[F_{j+1/2}^n(v) - F_{j-1/2}^n(v) \right].
\end{equation}
Here $F_{j+1/2}^n$ is the numerical flux between cell $T_j$ and $T_{j+1}$ at time $t^n$. The local Lax-Friedrichs scheme \cite{Toro} is utilized as 
\begin{equation}
    \label{eq:LF} 
    F^n_{j+1/2}(v) = \frac{1}{2}\left(v f_j^n + v f_{j+1}^n\right) -\frac{\lambda_{j+1/2}}{2}\left(f_{j+1}^n- f_{j}^n\right),
\end{equation}
where $\lambda_{j+1/2}=\max \left \{\left | \lambda_{j+1/2}^L \right |, \left | \lambda_{j+1/2}^R\right | \right \}$,
$\lambda_{j+1/2}^L$ and $\lambda_{j+1/2}^R$ are the fastest signal speeds in Proposition \ref{thm:global} as
\begin{equation}
    \label{eq:char_vel}
    \begin{aligned}
    \lambda_{j+1/2}^L & = \min\left\{u_j^n - C_{M+1}\sqrt{\theta_j^n}, u_{j+1}^n - C_{M+1}\sqrt{\theta_{j+1}^n} \right\}, \\
     \lambda_{j+1/2}^R& = \max\left\{u_j^n + C_{M+1}\sqrt{\theta_j^n}, u_{j+1}^n + C_{M+1}\sqrt{\theta_{j+1}^n} \right\}.
    \end{aligned}
\end{equation}
In the numerical flux \eqref{eq:LF}, the moment coefficient $f_{j, M+1}^n$ will be utilized when calculating $v f_{j}^n(v)$. In IPNC, $f_{j, M+1}^n$ is derived by the machine learning method, which will be adopted directly when computing the numerical flux \eqref{eq:LF}.  
%Moreover, since we do not have an analytic formula for $f_{j, M+1}^n$, it is difficult to derive the exact eigenvalues for the closed moment system obtained by the IPNC method. Therefore, the original eigenvalues in HME is adopted here. 

\begin{remark}
In the Grad method, $f_{j, M+1}^n$ is set directly as $0$,  while in HME, there is a non-conservative regularization term to obtain the closed system, and the numerical flux for the non-conservative regularization term is specially designed \cite{FanDissertation}. 
\end{remark}

For the collision step, it can be solved analytically for the BGK model as 
\begin{equation}
    \label{eq:collision}
    f^{n+1}_{j, \alpha} = f^{n+1, \ast}_{j, \alpha} \exp\left(-\frac{\Delta t^n}{\Kn}\right), \qquad 2 < \alpha \leqslant  M. 
\end{equation}
The time step is decided by the CFL condition as
\begin{equation}
    \label{eq:CFL}
    \Delta t^n = {\rm CFL} \frac{\Delta x}{\lambda_{\rm max}}, \qquad \lambda_{\rm max} = \max_j \left(|u_j^n| + C_{M+1}\theta_j^n\right).
\end{equation}
Throughout our numerical experiments, the linear reconstruction is utilized in the spatial space, and the CFL condition is set as $\rm CFL = 0.45$. 

\subsection{With smooth initial conditions}
\label{sec:num_smooth_dis}

We first validate IPNC on the problems with smooth initial conditions, which we will call it the wave problem for short. This is the same example from \cite{han2019uniformly}. The initial conditions are generated from the combination of two Maxwellian, who has the form below
\begin{equation}
    \label{eq:ex1_max}
    \mM^{\boldsymbol{U}}(x, v) = \frac{\rho(x)}{\sqrt{2 \pi \theta(x)}} \exp\left(-\frac{(v-u(x))^2}{2 \theta(x)}\right), \qquad \boldsymbol{U(x)} = (\rho(x), u(x), \theta(x)), 
\end{equation}      
with 
\begin{equation}
\label{eq:ex1_macro}
\left\{\begin{array}{l}
\rho(x)=a_{\rho} \sin \left(2 k_{\rho} \pi x / L+\phi_{\rho}\right)+b_{\rho}, \\
u(x)=0, \\
\theta(x)=a_{\theta} \sin \left(2 k_{\theta} \pi x / L+\phi_{\theta}\right)+b_{\theta},
\end{array}\right.
\end{equation}
where $a_{s}, k_{s}, \phi_{s}, b_{s}$ with $ s = \rho, \theta$ are all parameters to obtain different initial data. Here $a_{\rho}$ and $a_{\theta}$ are uniformly sampled from $[0.2,0.3]$, $b_{\rho}$ and $b_{\theta}$ are uniformly sampled from $[0.5,0.7]$, $\phi_{\rho}$ and $\phi_{\theta}$ are uniformly sampled from $[0, 2\pi]$, while $k_{\rho} \in \bbZ$ and $k_{\theta} \in \bbZ$ are randomly chosen from ${1,2,3,4}$. The continuous initial data is combined in the same way as \cite{han2019uniformly}  
\begin{equation}
\label{eq:ex1_ini}
f_{\text {smooth}}=\frac{\alpha_{1} \mM^{\boldsymbol{U}_1}(x, v) + \alpha_{2} \mM^{\boldsymbol{U}_2}(x, v)}{\alpha_{1}+\alpha_{2}+10^{-6}},
\end{equation}
where $\boldsymbol{U}_1$ and $\boldsymbol{U}_2$ are randomly generated from \eqref{eq:ex1_macro}, and $\alpha_1$ and $\alpha_2$ are uniformly sampled from $[0, 1]$. 

In order to compare with the HermMLC method in \cite{han2019uniformly}, DVM with the same IMEX numerical scheme \cite{Jin2010} as in \cite{han2019uniformly} is utilized to generate the training data. To be more precise, the computational region and the number of grids in the spatial space are set as $L = [-0.5, 0.5]$ and $N_x = 100$. The computational region in the microscopic velocity space is $[-10, 10]$ with $400$ grid points. The Knudsen number is randomly sampled from a log-uniform distribution on $[-3, 1]$ respect base $10$ as in \cite{han2019uniformly}. For IPNC, the numerical method in Sec. \ref{sec:num_sch} is utilized. The computational region and the number of grids in spatial space are the same as DVM, and the moment number is set as $M = 5$. We choose the same moment number $M=5$ in HermMLC and HME. 

In the numerical experiment, the end-to-end training approach in Sec. \ref{sec:end_end} is applied with the loss function \eqref{eq:loss_new}. Here, the training data set contains data in the time period of $[0,0.1]$, that is $\mathcal{I}=\{0\leqslant i \leqslant 96,i\in \mathbb{Z}\}$ with $t^i=0.001i$ and the length of the trajectory fragment $B$ is set as $B = 4$, taking both performance and efficiency into consideration. The sample number is set as $N_s = 100$. The testing data set contains $100$ samples which are generated from the initial condition \eqref{eq:ex1_ini} and Knudsen number as the training data. When generating the testing data, we fix the random seed to make sure that each time we can obtain the same $100$ samples so that the experiment can be reproduced. 

To validate the IPNC method, we have extended the simulation time to $t = 1.0$, where the generalization ability could also be tested. The density $\rho$, macroscopic velocity $u$, and temperature $\theta$ at $t = 0.1, 0.2$ and $1$ for a sample from the testing data is plotted in Figure \ref{fig:ex1_solution}. The initial condition for this specific sample is given in Appendix \ref{app:supp} . The numerical solutions by HermMLC, HME, and DVM are also plotted in Figure \ref{fig:ex1_solution}. We can see that, at all time steps, even at $ t = 1$ which is much longer than the final time of the training data, IPNC matches with the solution of DVM better than HME and HermMLC, which also validates the generalization ability and long time behavior of IPNC. 

\begin{figure}[!htb]
\centering
\subfloat[error quarterlies]{
\includegraphics[width=0.35\textwidth]{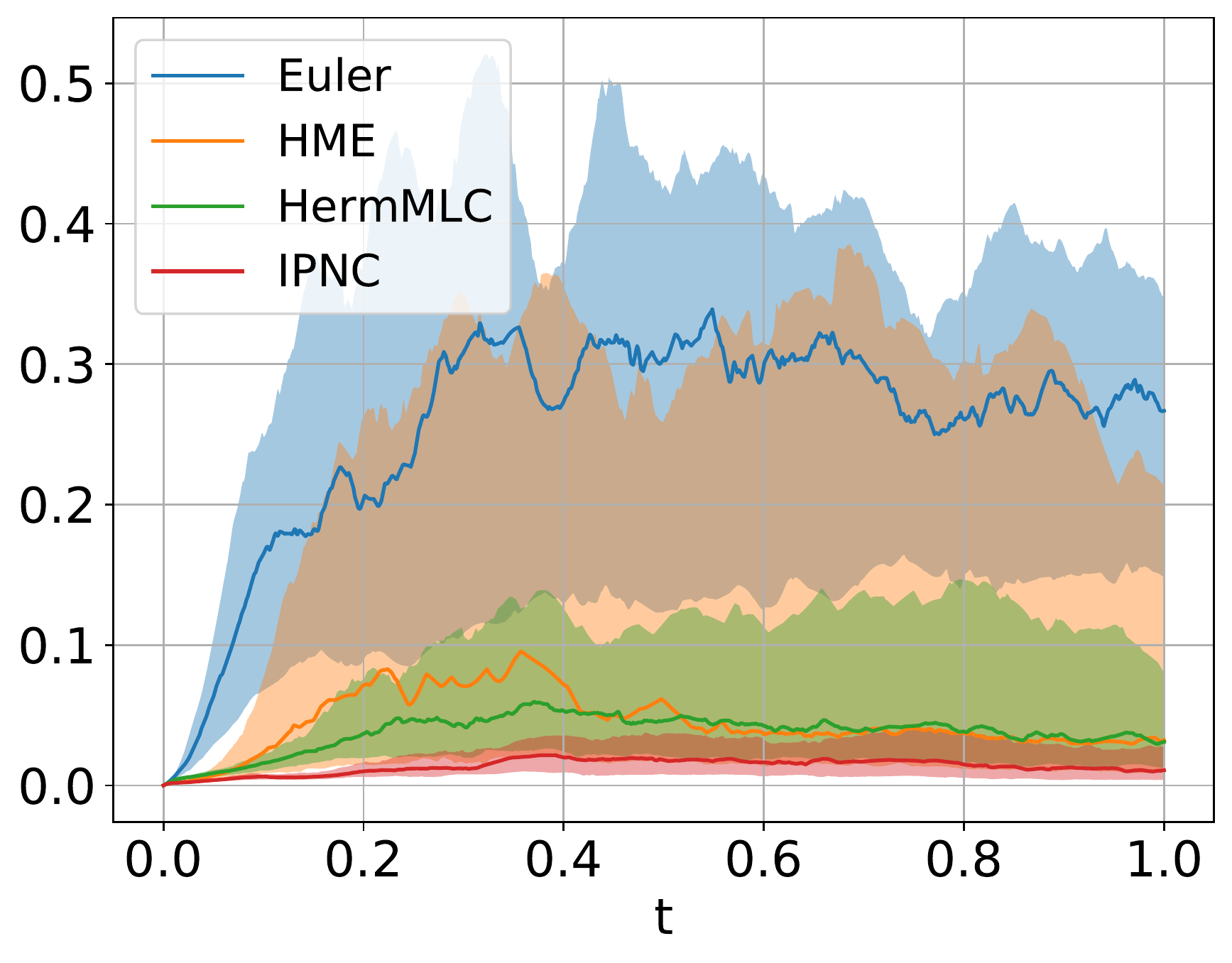}
 \label{fig:ex1_error_1_1000}} \qquad
\subfloat[different $\Kn$ at $t=0.1$]{
\includegraphics[width=0.36\textwidth]{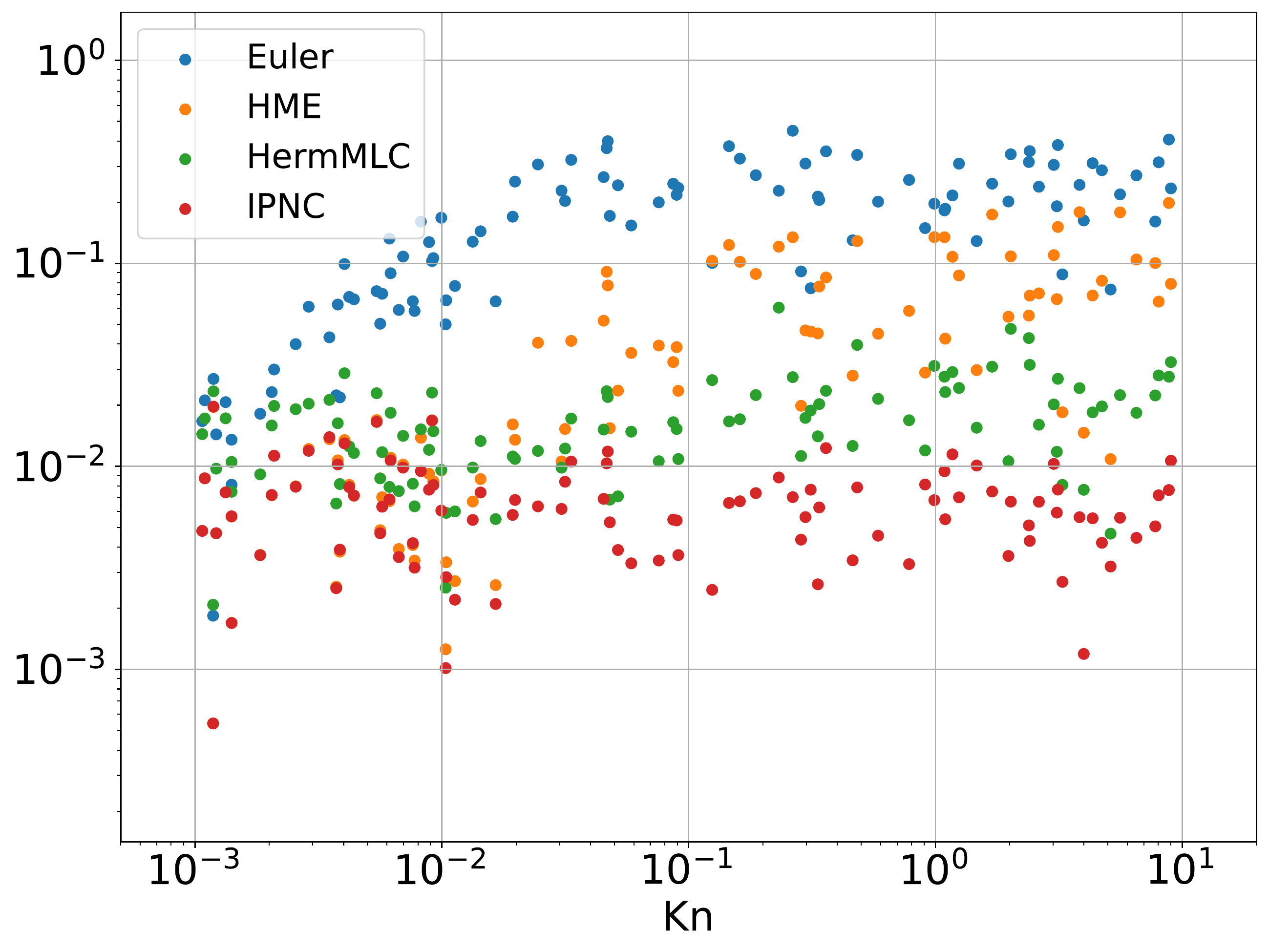}
\label{fig:ex1_error_2_100}} \\
\subfloat[different $\Kn$ at $t=0.2$]{
\includegraphics[width=0.35\textwidth]{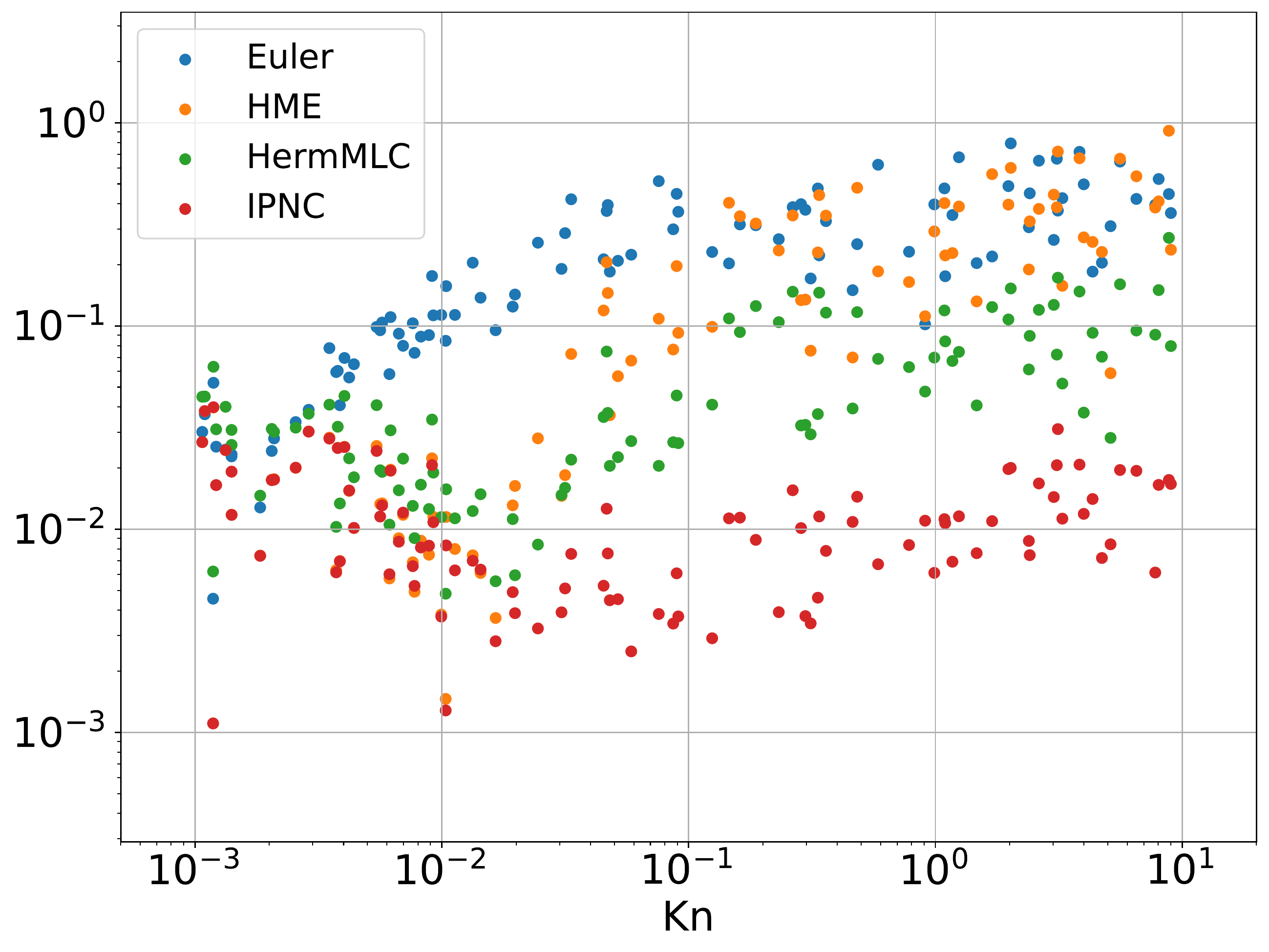}
\label{fig:ex1_error_2_200}}
\qquad
\subfloat[different $\Kn$ at $t=1.0$]{
\includegraphics[width=0.35\textwidth]{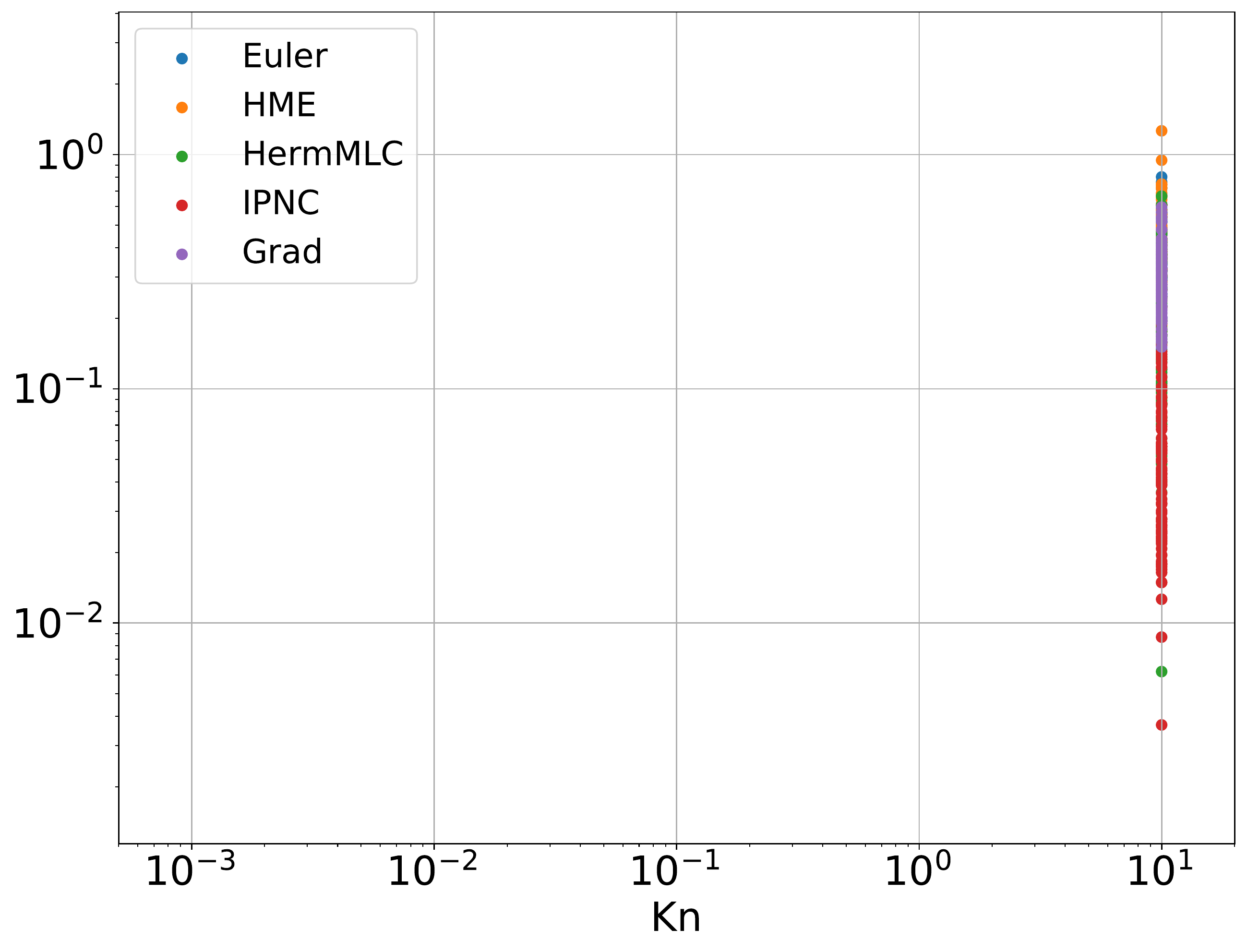}
\label{fig:ex1_error_2_1000}}
\caption{(Sec. \ref{sec:num_smooth_dis}) (a) The time evolution of the distribution of the relative error with $100$ samples of the initial condition \eqref{eq:ex1_ini} obtained by the different methods, where the $x$-axis is time $t$ and the $y$-axis is the relative error \eqref{eq:macro_error}. Here the solid line indicates the median and the translucent region is the area between the 1/4 and 3/4 quarterlies. (b-d) The relative error \eqref{eq:macro_error} at time $t = 0.1, 0.2$ and $1$, respectively, of the $100$ initial samples with different $\Kn$ and different initial conditions \eqref{eq:ex1_ini}. The $x$-axis is the Knudsen number and the $y$-axis is the relative error obtained by \eqref{eq:macro_error}.  }
\label{fig:waveH_solution_error_200}
\end{figure}

To quantitatively validate the performance of IPNC, we define the relative error of the macroscopic variables for the sample $s$ at the $i$-th time step as 
\begin{equation}
\label{eq:macro_error}
{\varepsilon}^{s,i} = \sqrt{\frac{1}{3}\left[ \frac{\| \overline{\rho}^{s, i} - \rho^{s,i}\|_x}{\|\rho^{s, i}\|_x} + \frac{\| \overline{u}^{s, i} - u^{s,i}\|_x}{1+\|u^{s, i}\|_x} + \frac{\| \overline{\theta}^{s, i} - \theta^{s,i}\|_x}{\|\theta^{s, i}\|_x}\right] },
\end{equation}
where the norm $\|\cdot\|_x$ is defined as 
$\|\omega\|_x = \sum\limits_{j = 1}^{N_x} \omega_j^2$. Here $\overline{\omega}^{s,i}$ with $\omega = \rho, u, \theta$ is the numerical solution obtained by IPNC with $s$-th initial condition at time step $i$, and $\omega^{s, i}$ is the reference solution obtained by DVM. The relative errors \eqref{eq:macro_error} of the solutions obtained by Euler, HME, HermMLC, and IPNC are calculated. Figure \ref{fig:ex1_error_1_1000} shows the distribution of the relative errors for different initial data with different methods, where the $1/4$, $3/4$ quarterlies and the median line of the distribution of the relative error for the $100$ samples with the time evolution is plotted. We can see that the median line of the error obtained by IPNC is much smaller than the other three methods, and the translucent region between $1/4$ and $3/4$ quarterlies is also much smaller. From this, we can see that the relative error obtained by IPNC is smaller and the distribution of this error is also much sharper.

To test the stability of IPNC with respect to the Knudsen number $\Kn$, the numerical tests with different $\Kn$ are conducted. Consider $\Kn\in [0.001, 10]$, where $100$ samples with different $\Kn$ and different initial data are generated uniformly at random. The distribution of the relative error of the solutions obtained by different numerical methods with these initial conditions at $t = 0.1, 0.2$ and $1$ is shown in Figure \ref{fig:ex1_error_2_100} to   \ref{fig:ex1_error_2_1000}, respectively. We can see that at all time steps, IPNC has a clear advantage over the other compared methods, indicating the strong stability of IPNC with respect to the Knudsen number. To show the behavior of IPNC with different Knudsen numbers quantitatively, we choose five different Knudsen numbers as $\Kn = 0.01, 0.1, 1.0$, and $10.0$ ranging from the continuous regime to rarefied gas regime. For each $\Kn$, $100$ samples with different initial data \eqref{eq:ex1_ini} are generated, and the average of the relative errors for different Knudsen number at time $t = 0.1, 0.2$ and $1$ is shown in Table \ref{tab:waveH}. The same relative errors of the Euler model, HME \cite{FanDissertation}, and HermMLC  \cite{han2019uniformly} are all listed in Table \ref{tab:waveH}, from which we can see that IPNC is the most accurate method.

\begin{figure}[!htb]
\centering
\subfloat[$\rho$]{
\begin{minipage}[c]{0.3\textwidth}
\includegraphics[width=\textwidth, height=0.75\textwidth]{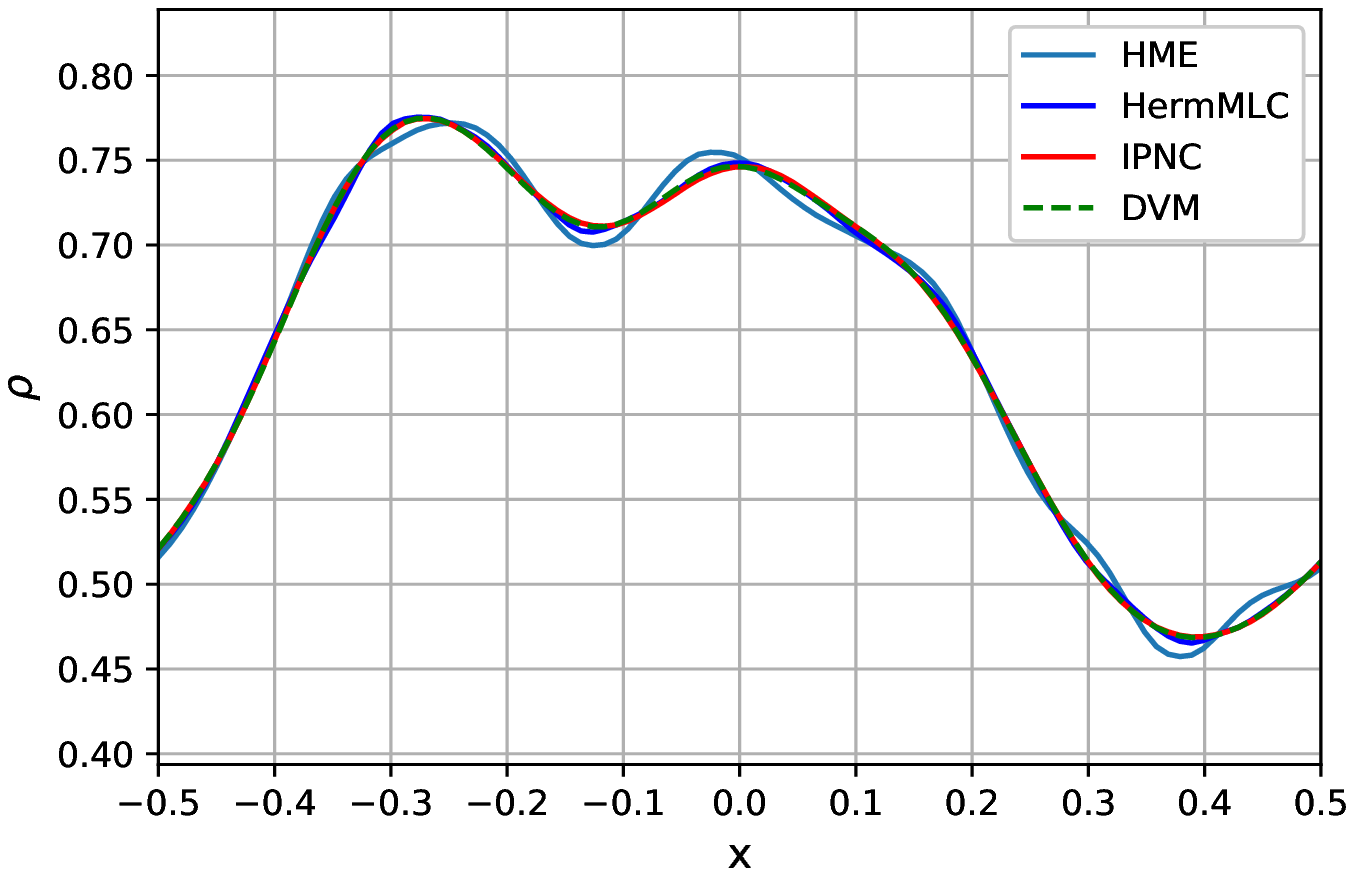}
\includegraphics[width=\textwidth, height=0.75\textwidth]{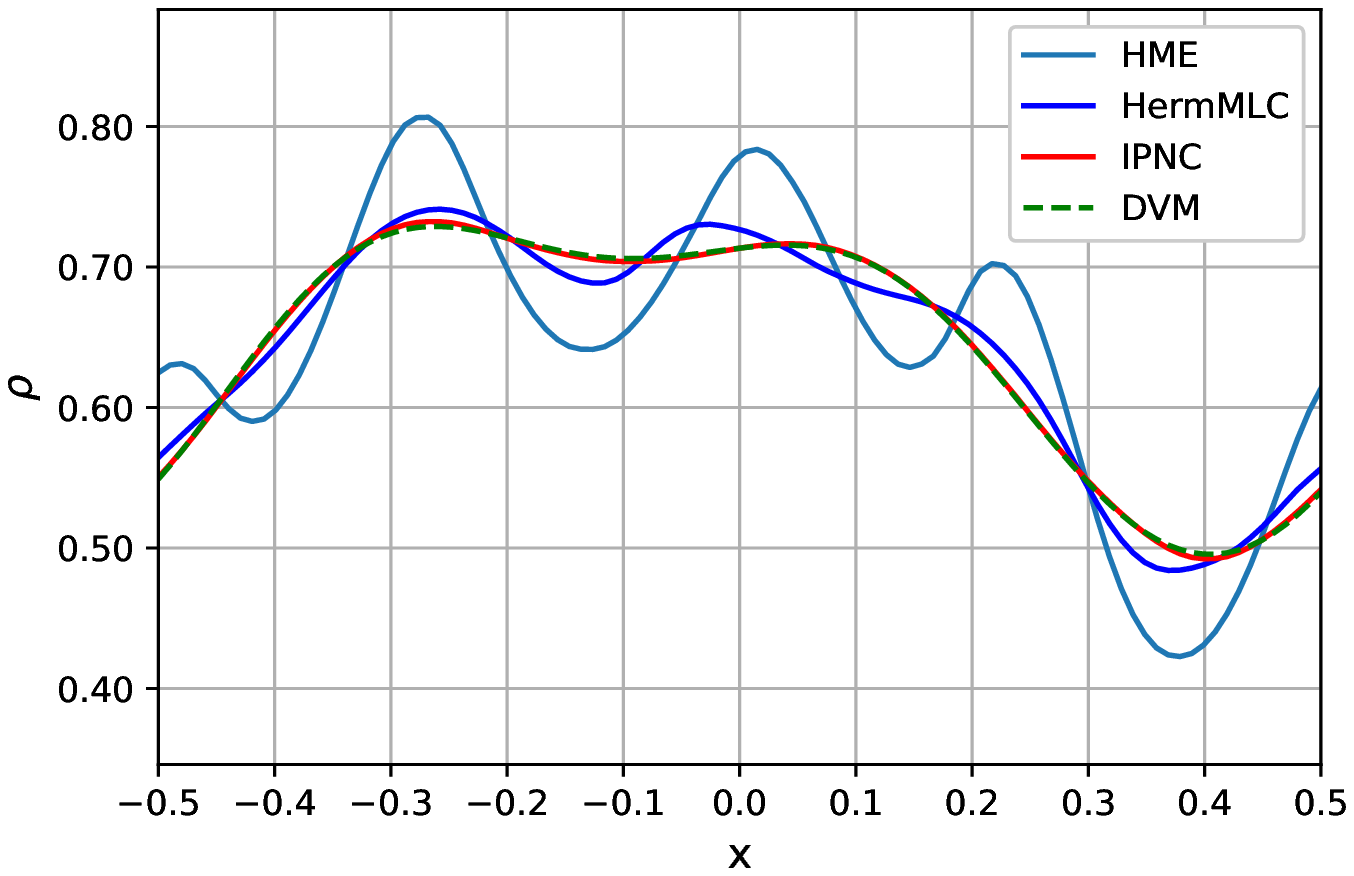}
\includegraphics[width=\textwidth, height=0.75\textwidth]{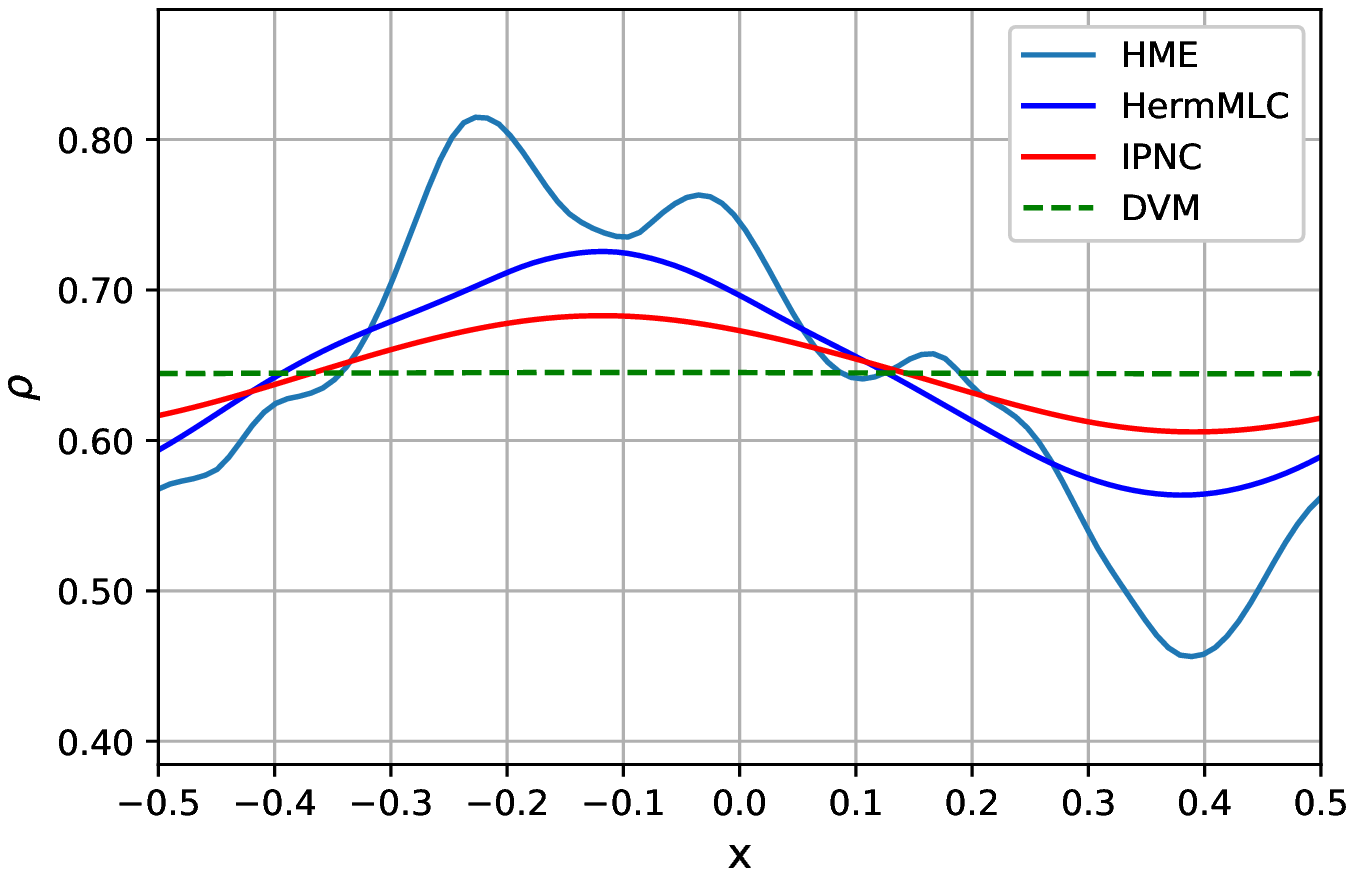}
\end{minipage}
}\quad
\subfloat[$u$]{
\begin{minipage}[c]{0.3\textwidth}
\includegraphics[width=\textwidth, height=0.75\textwidth]{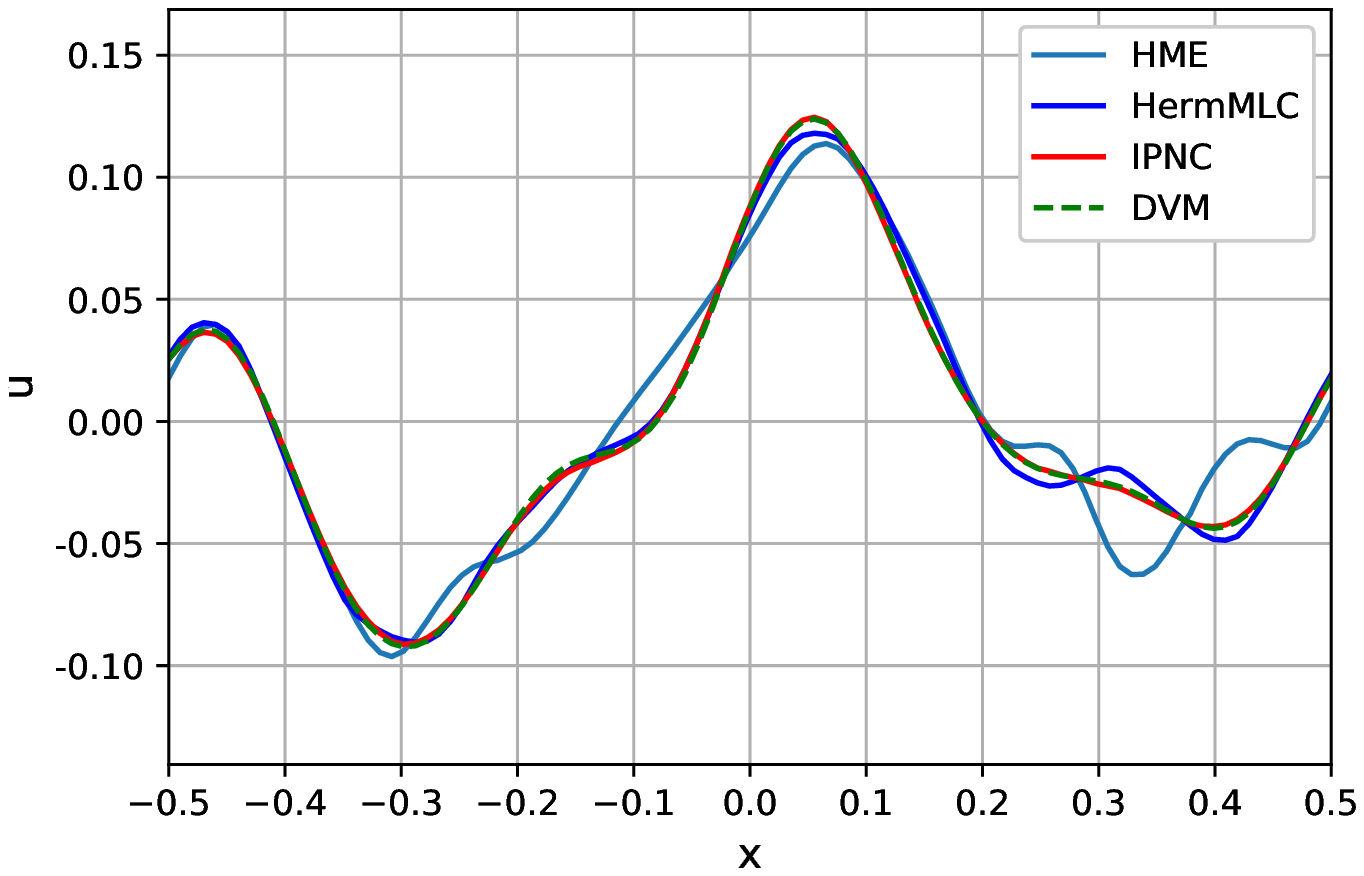}
\includegraphics[width=\textwidth, height=0.75\textwidth]{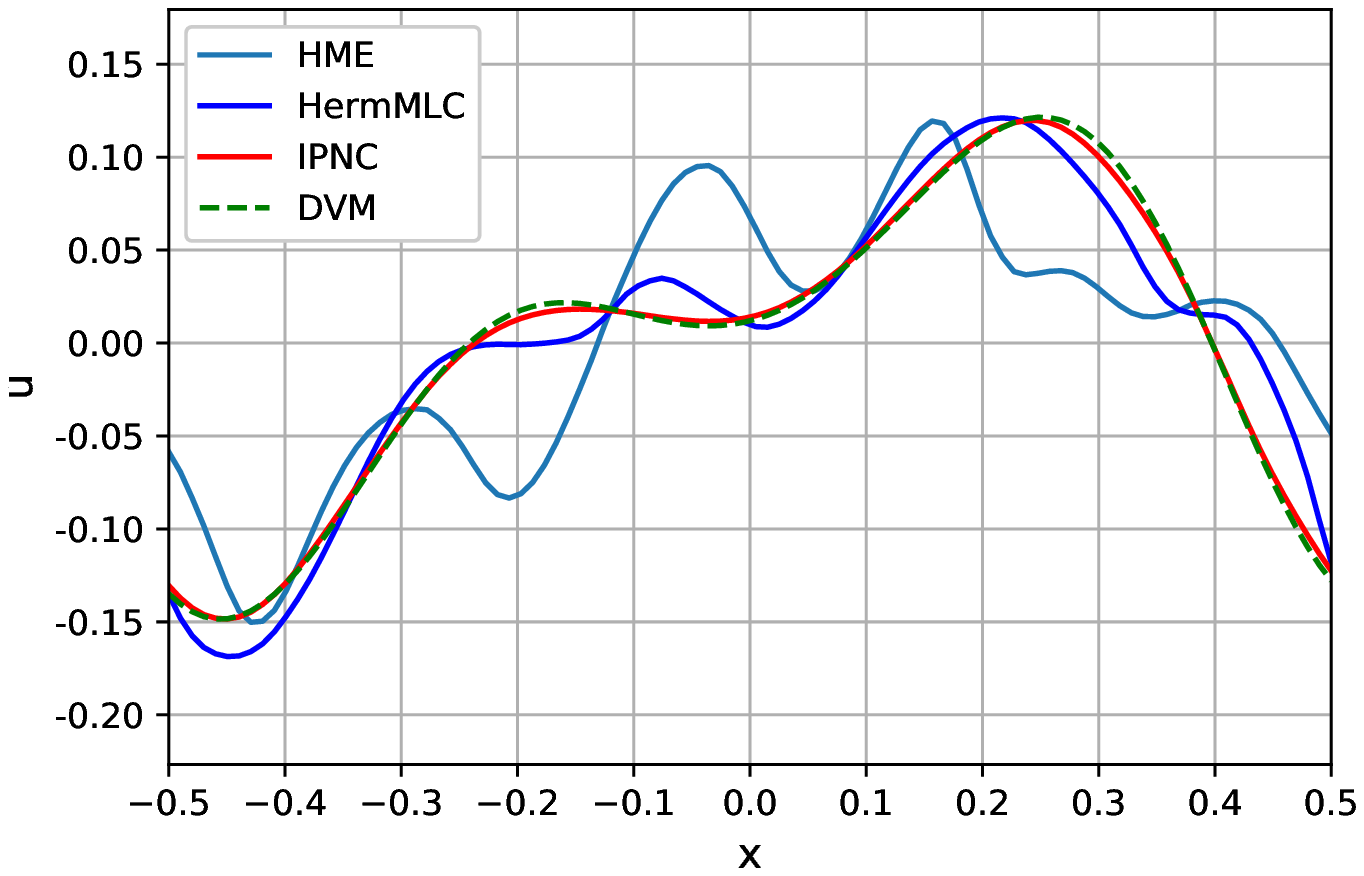}
\includegraphics[width=\textwidth, height=0.75\textwidth]{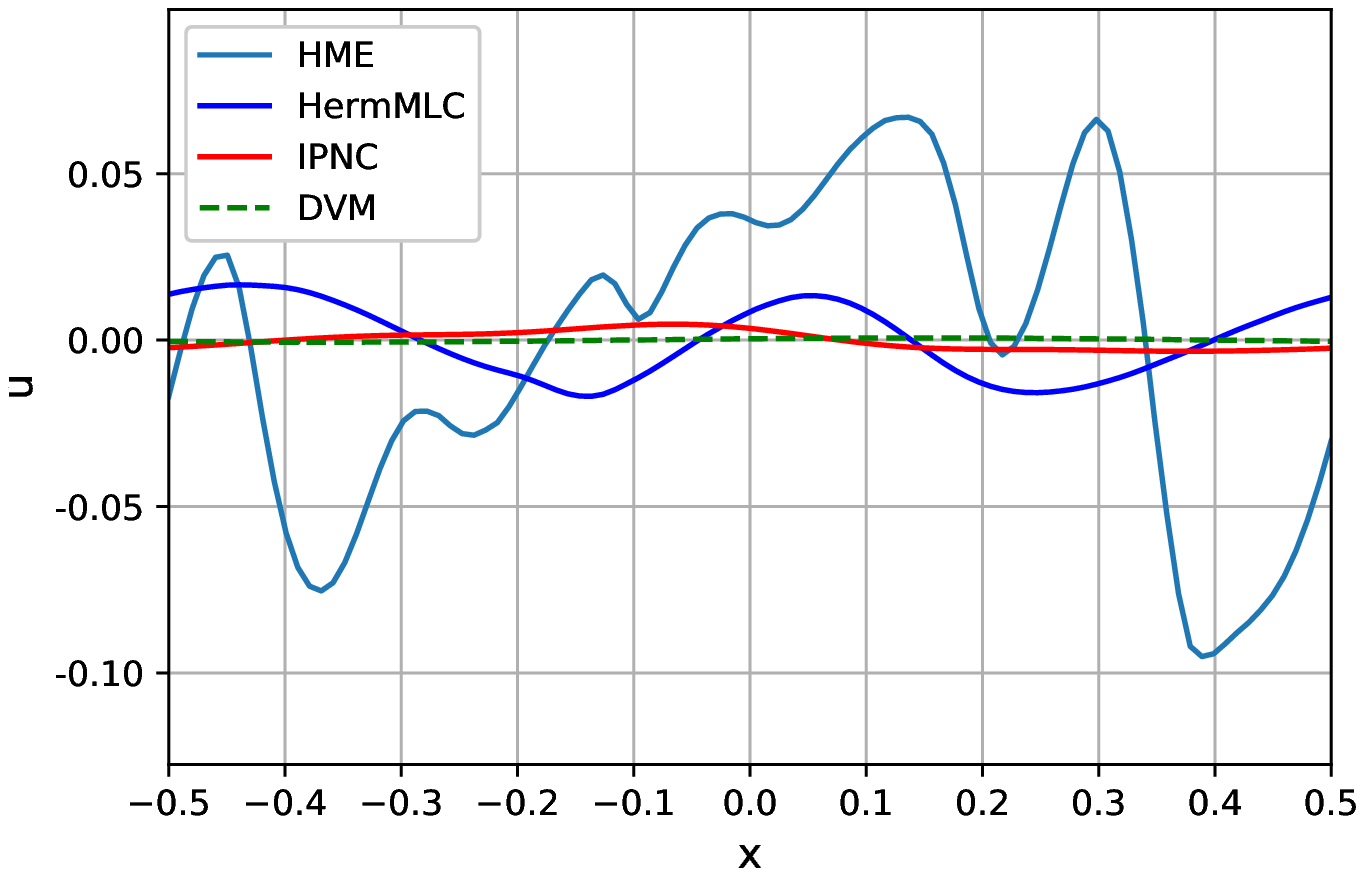}
\end{minipage}
}\quad
\subfloat[$\theta$]{
\begin{minipage}[c]{0.3\textwidth}
\includegraphics[width=\textwidth, height=0.75\textwidth]{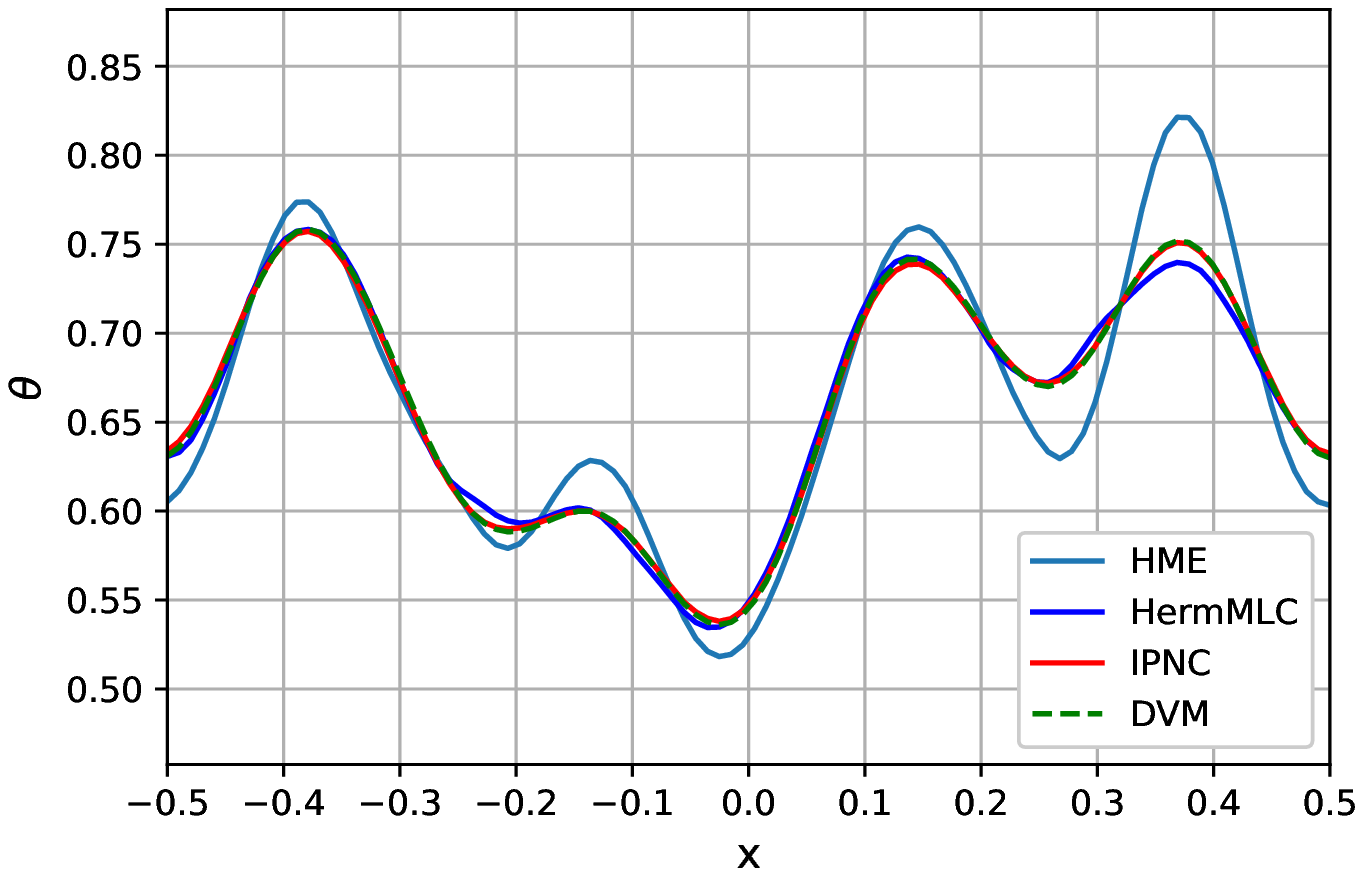}
\includegraphics[width=\textwidth, height=0.75\textwidth]{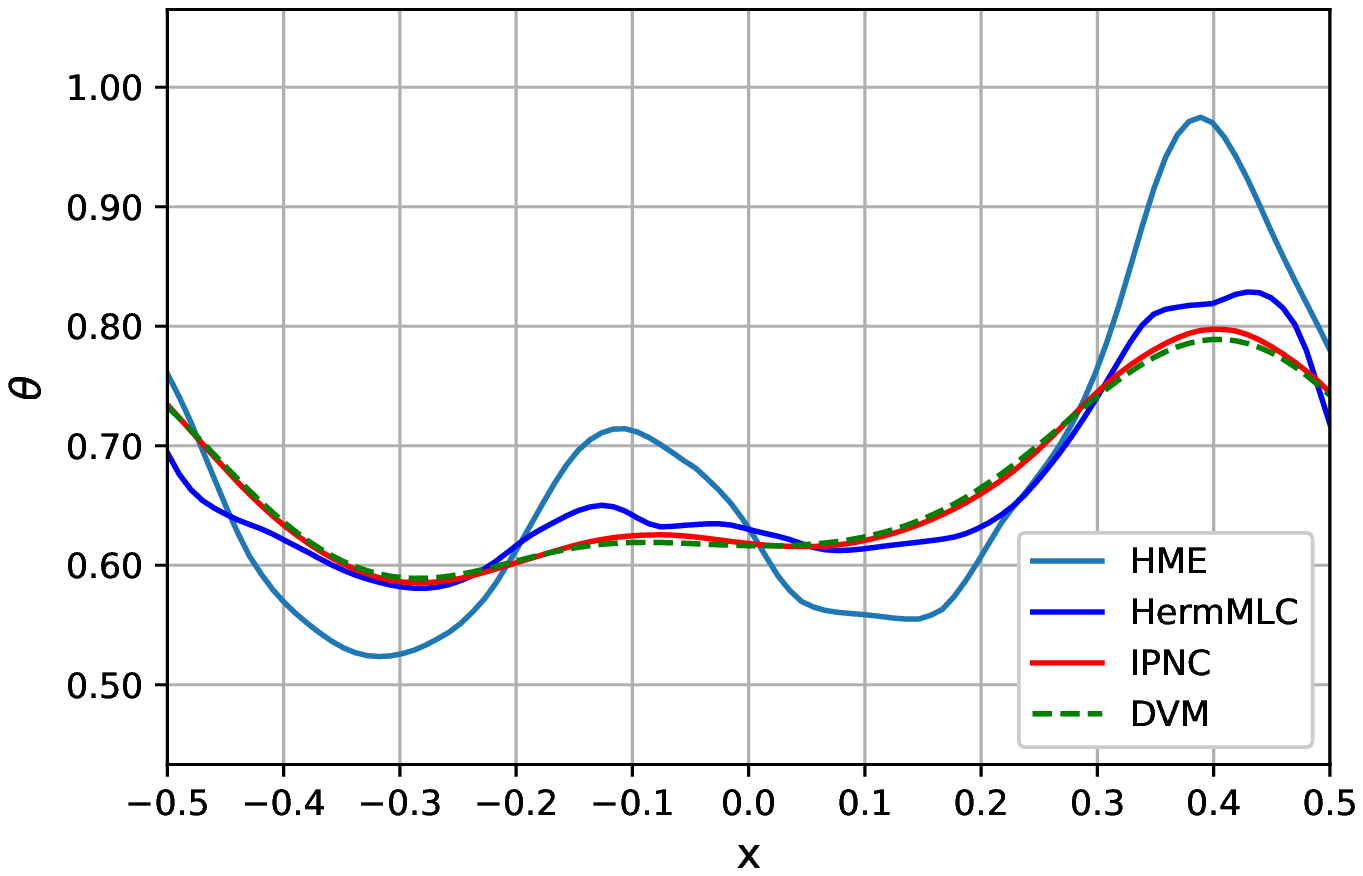}
\includegraphics[width=\textwidth, height=0.75\textwidth]{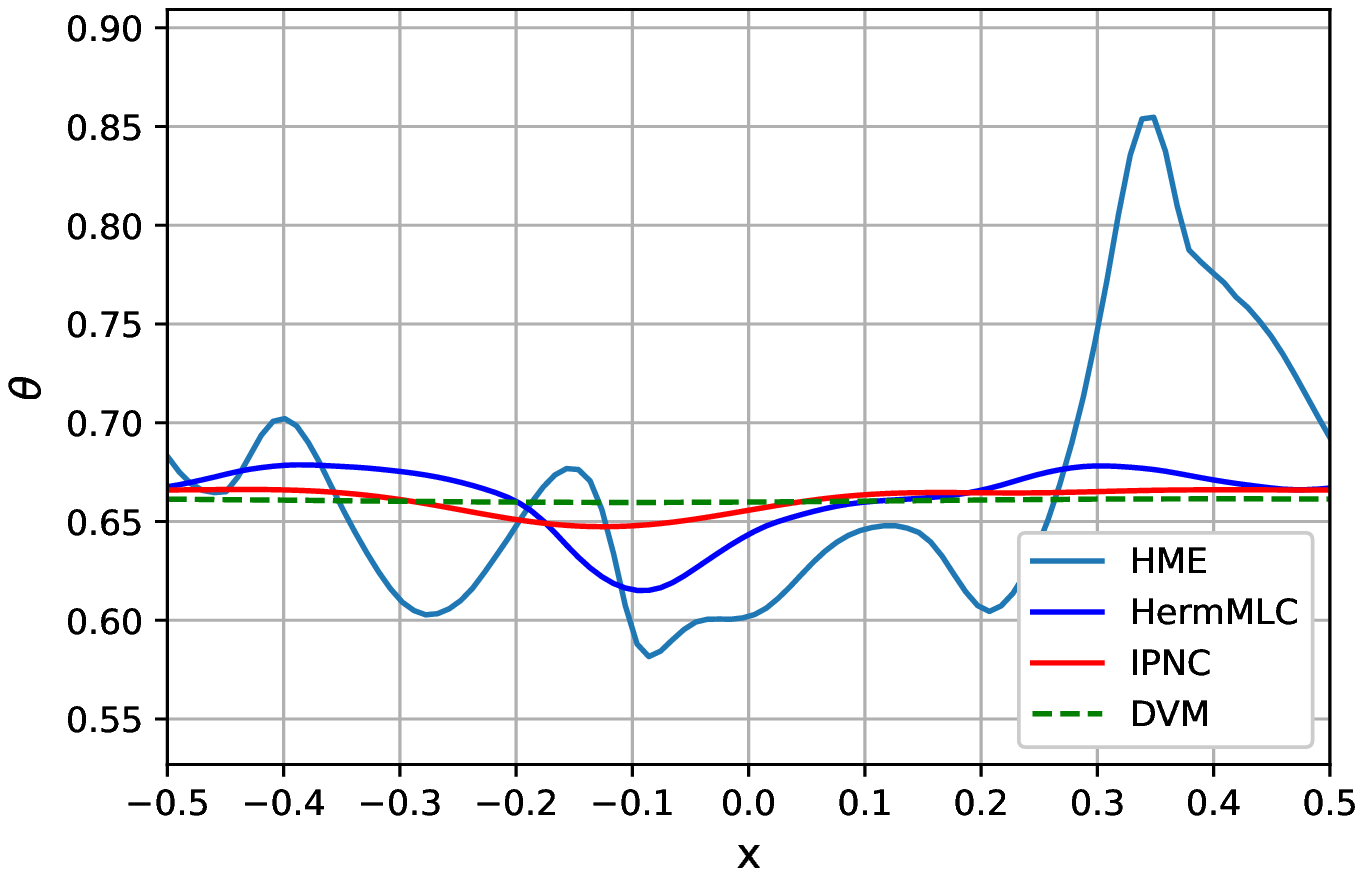}
\end{minipage}
}
\caption {(Sec. \ref{sec:num_smooth_dis}) Density $\rho$, macroscopic velocity $u$, and temperature $\theta$. The three rows are $t = 0.1$, $0.2$ and $1$, respectively. Here, the blue line is got by HME, the black line is got by HermMLC, the red line is got by IPNC, and the green dashed line is the reference solution got by DVM.}
\label{fig:ex1_solution}
\end{figure}

\begin{table}[ht]
\centering
\renewcommand\arraystretch{1.5}
\footnotesize
\begin{tabular}{l|lllll}
$ t = 0.1$ & Kn      & $0.01$            & $0.1$             & $1.0$             & $10$              \\ \hline
           & Euler   & $ 10.66 $         & $ 23.18 $         & $ 25.84 $         & $ 26.14 $         \\
           & HME     & $ 0.59 $          & $ 5.49 $          & $ 8.72 $          & $ 9.15 $          \\
           & HermMLC & $ 0.81 $          & $ 1.48 $          & $ 2.35 $          & $ 2.53 $          \\
           & IPNC    & $ \mathbf{0.47} $ & $ \mathbf{0.42} $ & $ \mathbf{0.57} $ & $ \mathbf{0.61} $ \\ \hline
$t = 0.2$  & Kn      & $0.01$            & $0.1$             & $1.0$             & $10$              \\ \hline
           & Euler   & $13.76$           & $33.66$           & $38.22$           & $38.73$           \\
           & HME     & $0.71$            & $16.37$           & $34.79$           & $37.85$           \\
           & HermMLC & $1.01$            & $4.69$            & $9.97$            & $11.06$           \\
           & IPNC    & $\mathbf{0.51}$   & $\mathbf{0.63}$   & $\mathbf{1.25}$   & $\mathbf{1.50}$   \\ \hline
$t = 1$    & Kn      & $0.01$            & $0.1$             & $1.0$             & $10$              \\ \hline
           & Euler   & $20.24$           & $36.51$           & $40.23$           & $40.52$           \\
           & HME     & $0.34$            & $4.36$            & $32.61$           & $-$               \\
           & HermMLC & $0.75$            & $1.69$            & $14.52$           & $21.20$           \\
           & IPNC    & $\mathbf{0.28}$   & $\mathbf{0.59}$   & $\mathbf{3.76}$   & $\mathbf{6.58}$  
\end{tabular}
\caption{(Sec. \ref{sec:num_smooth_dis})  Average of the relative error \eqref{eq:macro_error} obtained by the different numerical methods with the $100$ samples of the smooth initial condition \eqref{eq:ex1_ini} at time $t = 0.1, 0.2$ and $1$ with different Knudsen number.}
\label{tab:waveH}
\end{table}

\paragraph{Generalization of Knudsen number} To test the generalization ability of IPNC on Knudsen number for the continuous problem, we change the region of $\Kn$ to $\Kn \in [0.1, 1.0]$ in the training data set and keep the remaining remains the same.

\begin{figure}[!htb]
\centering
\subfloat[$\rho$]{
\begin{minipage}[c]{0.3\textwidth}
\includegraphics[width=\textwidth, height=0.75\textwidth]{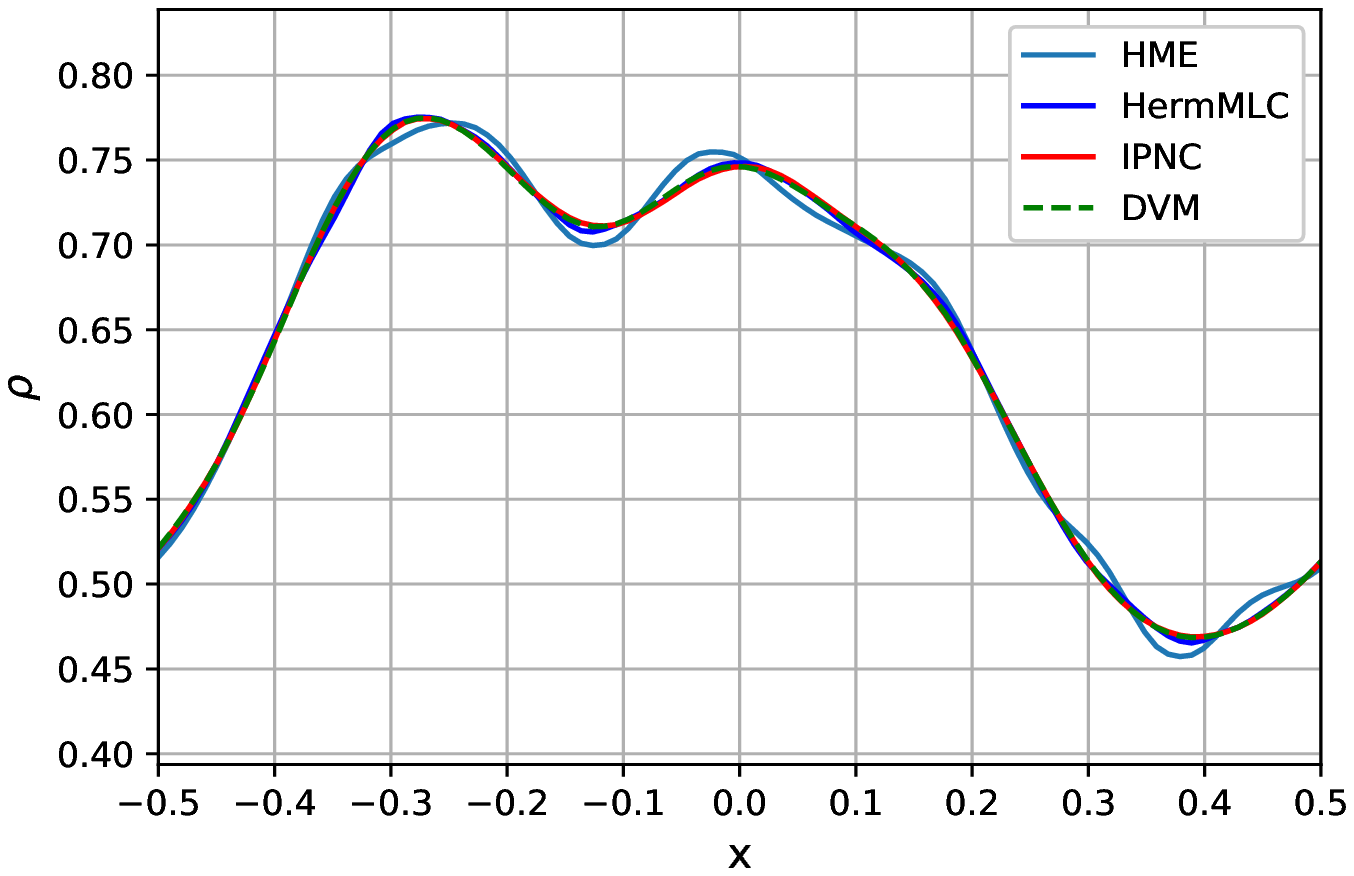}
\end{minipage}
}\quad
\subfloat[$u$]{
\begin{minipage}[c]{0.3\textwidth}
\includegraphics[width=\textwidth, height=0.75\textwidth]{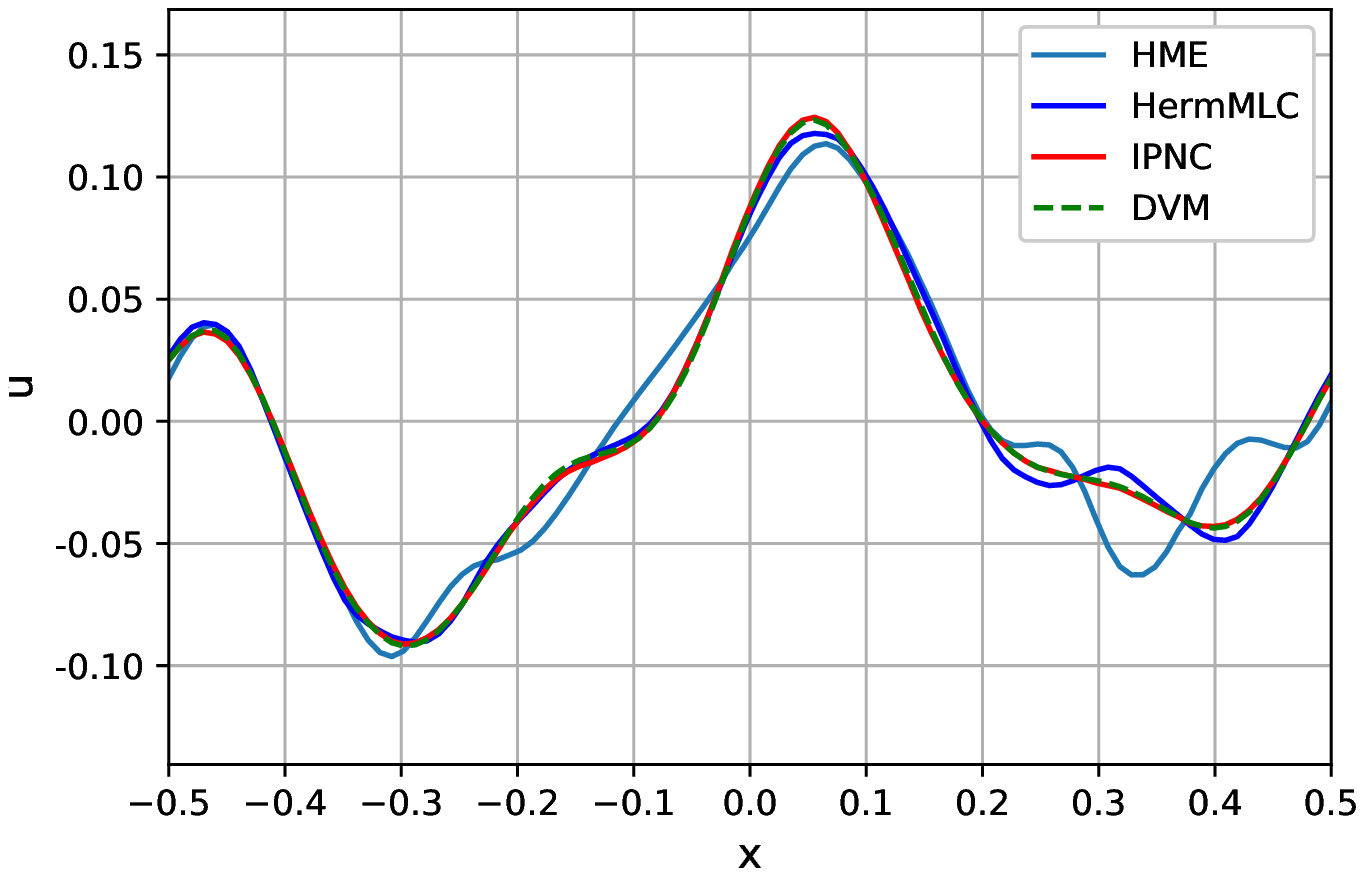}
\end{minipage}
}\quad
\subfloat[$\theta$]{
\begin{minipage}[c]{0.3\textwidth}
\includegraphics[width=\textwidth, height=0.75\textwidth]{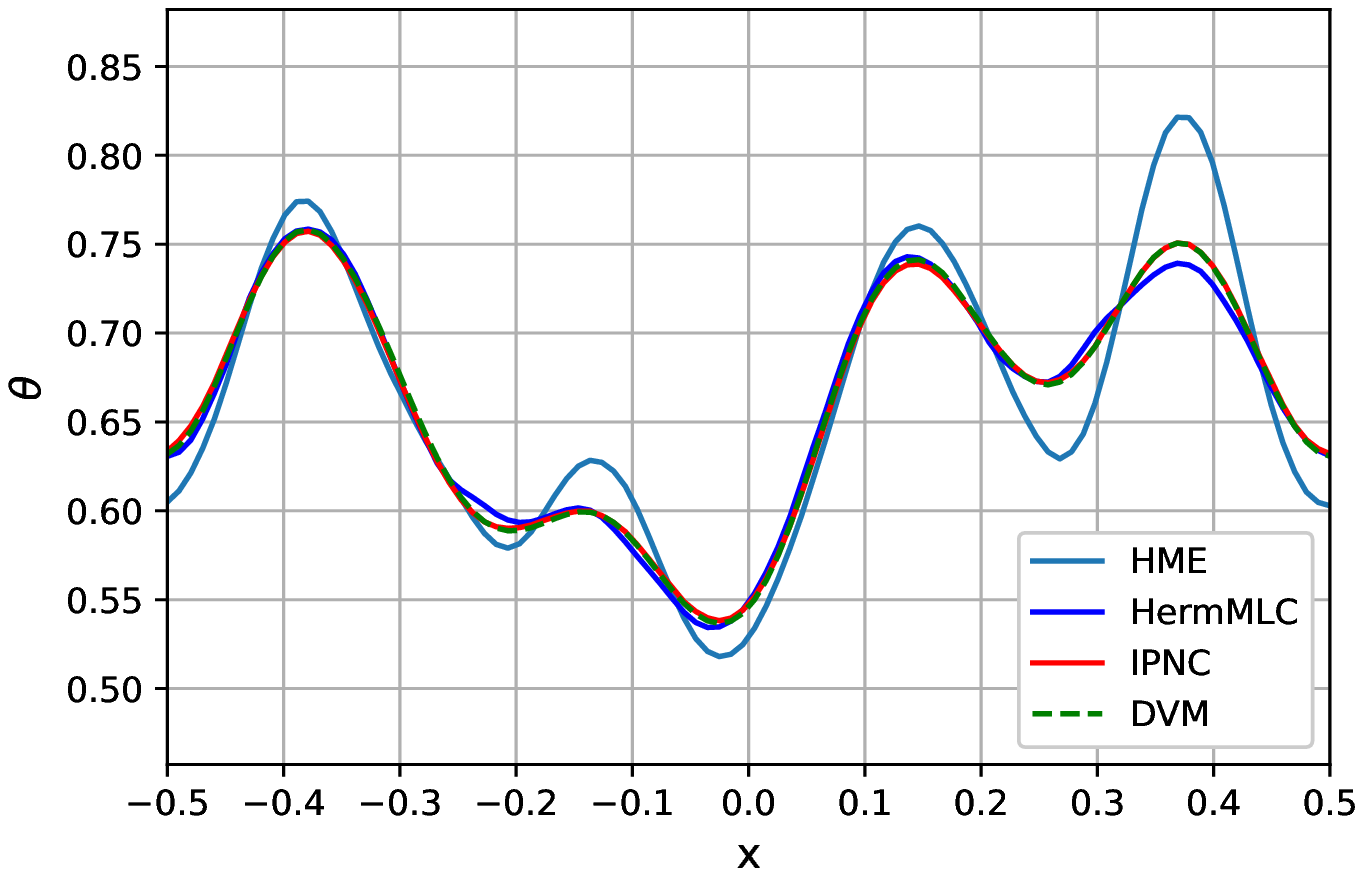}
\end{minipage}
}
\caption { (Sec. \ref{sec:num_smooth_dis}: Generalization of $\Kn$) Density $\rho$, macroscopic velocity $u$, and temperature $\theta$ at $ t= 0.1$ with $\Kn = 10$.  Here, the blue line is got by HME, the black line is got by HermMLC, the red line is got by IPNC, and the green dashed line is the reference solution got by DVM.}
\label{fig:ex1_ger_Kn_sol}
\end{figure}

The density $\rho$, macroscopic velocity $u$, and temperature $\theta$ at $ t = 0.1$ for the same sample as in Figure \ref{fig:ex1_solution} with $\Kn = 10$ is plotted in Figure \ref{fig:ex1_ger_Kn_sol}. We can find that even $\Kn$ is increasing to $\Kn = 10$, which is much larger than those in the training data, the numerical solution matches well with the reference solution by DVM. To show the generalization of IPNC on Knudsen number quantitatively, the similar average of the relative errors for the different Knudsen number at time $t = 0.1, 0.2$ and $1.0$ with the changed training region of $\Kn$ is shown in Table \ref{tab:waveH_ger_Kn}. It shows that even the training data set of $\Kn$ is reduced from $[0.001, 10]$ to $[0.1, 1]$, the average error does not increase, indicating the generalization ability of the Knudsen number for IPNC. 

\begin{table}[ht]
\centering
\renewcommand\arraystretch{1.5}
\footnotesize
\begin{tabular}{l|lllll}
$ t = 0.1 $ & Kn         & $0.01$  & $0.1$  & $1.0$   & $10$   \\ \hline
            & Kn         & $0.47$  & $0.42$ & $0.57$  & $0.61$ \\
            & $\Kn_{\rm lim}$ & $0.49 $ & $0.42$ & $0.53 $ & $0.70$ \\ \hline
$ t = 0.2$  & Kn         & $0.01$  & $0.1$  & $1.0$   & $10$   \\ \hline
            & Kn         & $0.51$  & $0.63$ & $1.25$  & $1.50$ \\
            & $\Kn_{\rm lim}$ & $0.55$  & $0.51$ & $1.49 $ & $1.49$ \\ \hline
$ t = 1$    & Kn         & $0.01$  & $0.1$  & $1.0$   & $10$   \\ \hline
            & Kn         & $0.28$  & $0.59$ & $3.76$  & $6.58$ \\
            & $\Kn_{\rm lim} $ & $0.23$  & $0.49$ & $4.64$  & $6.07$
\end{tabular}
\caption{
(Sec. \ref{sec:num_smooth_dis}: Generalization of $\Kn$)  Average of the relative error \eqref{eq:macro_error} got by IPNC with the $100$ samples of the smooth initial condition \eqref{eq:ex1_ini} at time $t = 0.1, 0.2$ and $1$ with different training Knudsen number. Here, $\Kn$ is  corresponding to the training set $\Kn \in [0.001, 10]$ and $\Kn_{\rm lim}$ is corresponding to the training set $\Kn \in [0.1, 1]$.
}
\label{tab:waveH_ger_Kn}
\end{table}

\begin{figure}[!htb]
\centering
\subfloat[$\rho$]{
\begin{minipage}[c]{0.3\textwidth}
\includegraphics[width=\textwidth, height=0.75\textwidth]{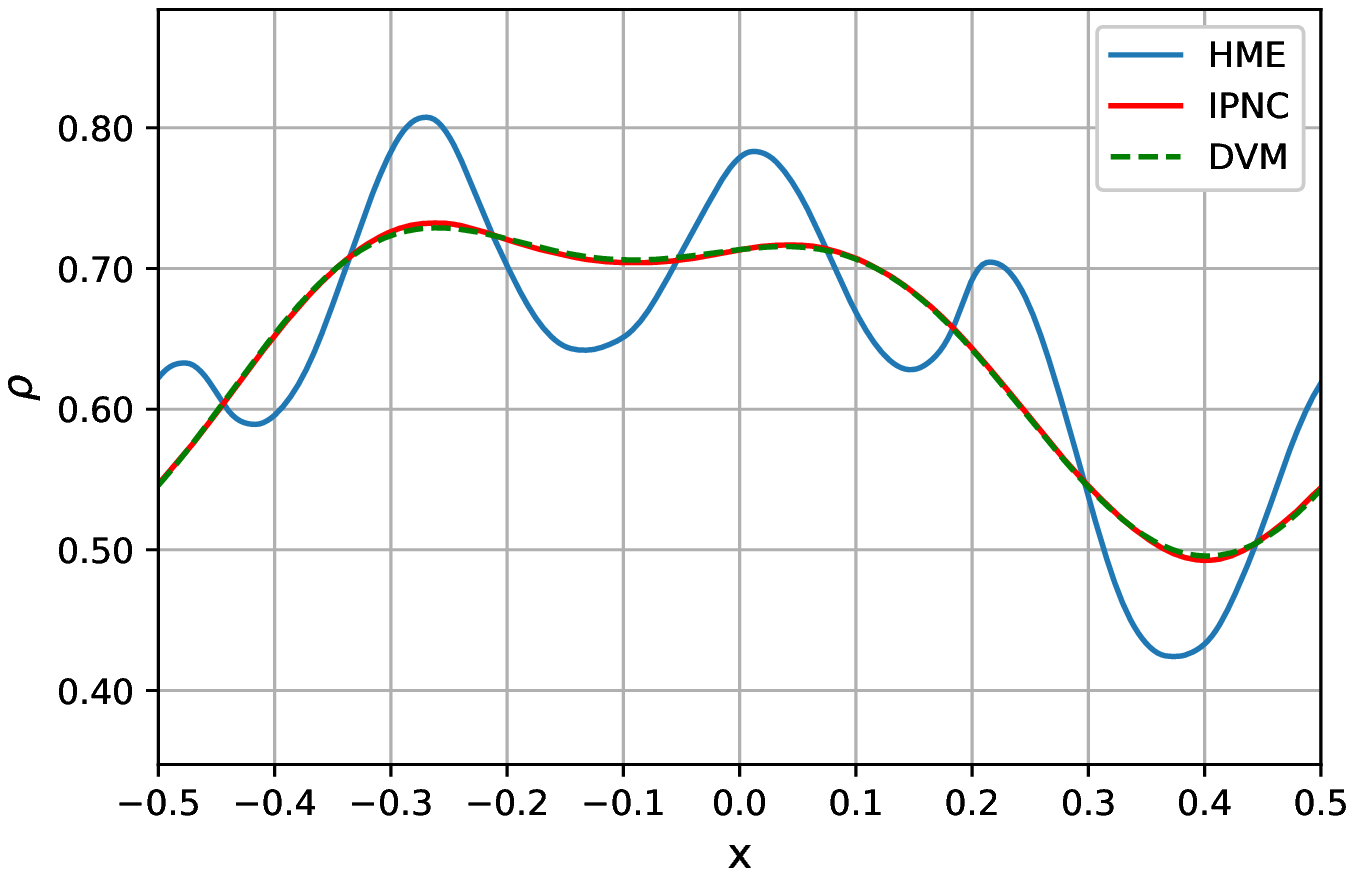}
\end{minipage}
}\quad
\subfloat[$u$]{
\begin{minipage}[c]{0.3\textwidth}
\includegraphics[width=\textwidth, height=0.75\textwidth]{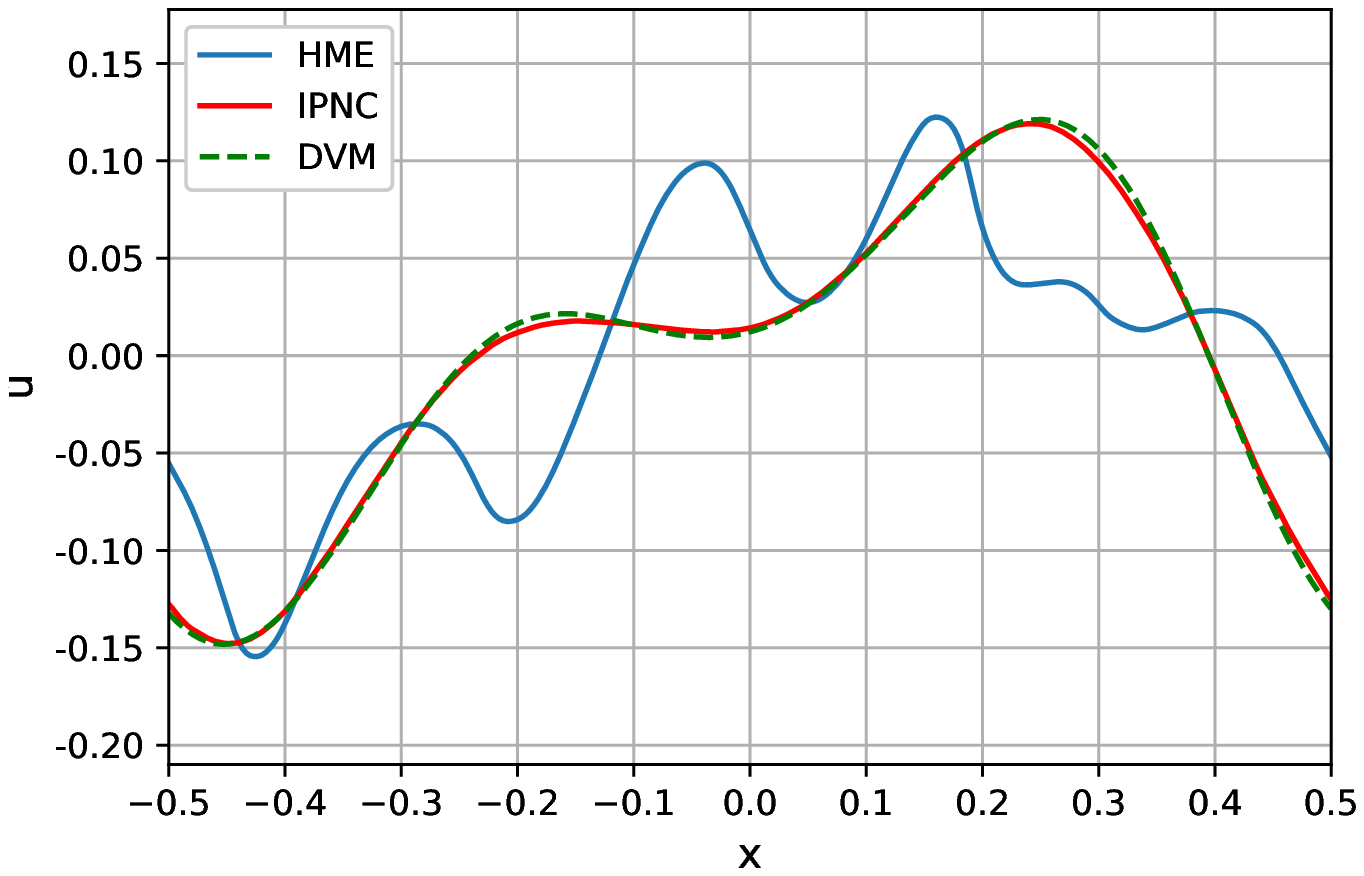}
\end{minipage}
}\quad
\subfloat[$\theta$]{
\begin{minipage}[c]{0.3\textwidth}
\includegraphics[width=\textwidth, height=0.75\textwidth]{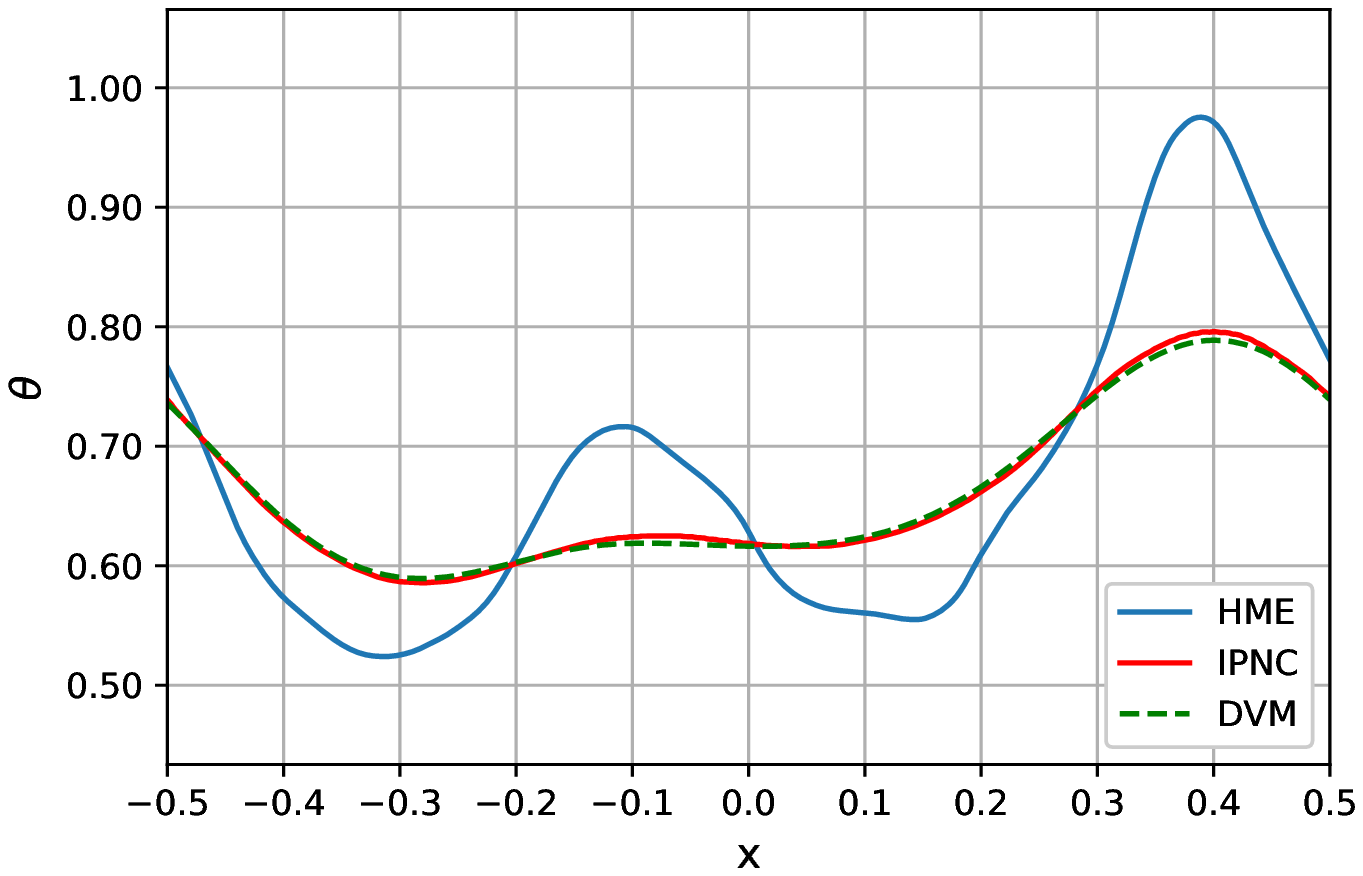}
\end{minipage}
}
\caption {(Sec. \ref{sec:num_smooth_dis}: Generalization of mesh size) Density $\rho$, macroscopic velocity $u$, and temperature $\theta$ at $ t= 0.2$ with mesh size $N_x = 300$.  Here, the blue line is got by HME, the black line is obtained by HermMLC, the red line is obtained by IPNC, and the green dashed line is the reference solution got by DVM.}
\label{fig:ex1_ger_dx_sol}
\end{figure}

\paragraph{Generalization of the mesh size} To show the generalization ability of IPNC on the mesh size for the continuous problem, the tests with different mesh sizes are studied. In the training data set, the mesh size is set at $N_x = 100$ with other settings remaining the same. 
% For problems with different $\Delta x$, we can use resample methods to interpolate to the same grid scale as we used to train the network for prediction and then interpolate back to the scale of the problem. That is, when we want to predict the $N_x=200$ problem, we resample the closure results to a resolution of $N_x=100$ by downsampling the data before the closure, and then upsample the closure results to a resolution of $N_x=200$ for continued calculation. Such a technique makes it easy to scale IPNC to a denser grid without a great loss of accuracy.
In the testing data, the mesh size $N_x = 100, 200, 300$ and $400$ are tested, respectively. To generate the closure network from $N_x = 100$ to large mesh size, the interpolation method is utilized in the closure step. The numerical solution with $N_x = 300$ at $t = 0.2$ is shown in Figure \ref{fig:ex1_ger_dx_sol}. It shows that with increasing mesh size, the numerical solution by IPNC matches well with the DVM solution. To show the generalization quantitatively, the average relative error \eqref{eq:macro_error} got by IPNC with the $100$ samples at time $t = 0.2$ with different mesh size is calculated in Table \ref{tab:ex1_ger_dx_error}, where compared to that obtained by HME, this error is much smaller. Moreover, with the increasing mesh size, this error does not increase greatly, which also illustrates the generalization ability of IPNC on mesh size.

\begin{table}[ht]
\centering
\renewcommand\arraystretch{1.5}
\footnotesize
% \begin{tabular}{l|l|l|l||l|l|l|l}
% Mix  & $N_x=100$ & $N_x=200$ & $N_x=300 $ & Wave & $N_x=100$ & $N_x=200$ & $N_x=300 $ \\ \hline
% HME  & $13.12$   & $14.40$   & $13.52$    & HME  & $16.33$   & $21.18$   & $23.08$    \\
% IPNC & $1.18$    & $1.16$    & $1.14$     & IPNC & $1.09$    & $1.04$    & $1.38$    
% \end{tabular}
% \caption{Network trained on Mix with nx=100, prediction error at T=0.2s for different nx.}
\begin{tabular}{l|llll}
     & $N_x = 100$ & $N_x = 200$ & $N_x = 300$ & $N_x = 400$ \\ \hline
HME  & $16.33$     & $21.18$     & $23.08$     & $26.45$     \\
IPNC & $1.09$      & $1.04$      & $1.38$      & $1.20$     
\end{tabular}
\caption{(Sec. \ref{sec:num_smooth_dis}: Generalization of  mesh size)  Average of the relative error \eqref{eq:macro_error} obtained by IPNC with the $100$ samples of the smooth initial condition \eqref{eq:ex1_ini} at time $t = 0.2$ with different mesh size. }
\label{tab:ex1_ger_dx_error}
\end{table}

\subsection{With discontinuous initial conditions}
\label{sec:num_discon_dis}
In this section, the problems with the discontinuous initial conditions, which we will call mix problem for short, are studied, where the same problems were also tested in \cite{han2019uniformly}. The discontinuous initial condition is generated by mixing the smooth initial values \eqref{eq:ex1_ini} with the initial value of a Riemann problem. To be more precise, the initial condition takes the form as follows, 
\begin{equation}
\label{eq:ex2_ini}
f_{\operatorname{mix}}=\alpha f_{\text {smooth}}+(1-\alpha) \mM^{\boldsymbol{U}_{\rm shock}}(x, v),
\end{equation}
where ${\boldsymbol{U}_{\rm shock}}$ is randomly chosen from $\boldsymbol{U}_{\rm shock}^1$ and $\boldsymbol{U}_{\rm shock}^2$ with equal probability. Here, $\boldsymbol{U}_{\rm shock}^i, i = 1, 2$, satisfies
\begin{equation}
\begin{aligned}
\label{eq:ex2_rieman}
 \boldsymbol{U}_{\rm shock}^1 = \left\{
\begin{array}{ll}
 (\rho_l, u_l, \theta_l),    &  x\in [-0.5, x_1]~{\rm or}~[x_2, 0.5], \\
    (\rho_r, u_r, \theta_r),  &  x \in [x_1, x_2], 
\end{array}
\right.
\\
 {\rm or}\qquad 
\boldsymbol{U}_{\rm shock}^2 = \left\{
\begin{array}{ll}
 (\rho_r, u_r, \theta_r),   &  x\in [-0.5, x_1]~{\rm or}~[x_2, 0.5], \\
    (\rho_l, u_l, \theta_l),  &  x \in [x_1, x_2], 
\end{array}
\right.
\end{aligned}
\end{equation}
where $\rho_l$ and $\theta_l$ are sampled from the uniform distribution on $[1, 2]$, while $\rho_r$ and $\theta_r$ are sampled from the uniform distribution on $[0.55, 0.9]$, with $u_l = u_r = 0$. Here, $x_1$ and $x_2$ are two random variables sampled from the uniform distributions on $[-0.3, -0.1]$ and $[0.1, 0.3]$ respectively. In this experiment, the numerical setting including the methods to generate the training and testing data, the end-to-end training approach \ref{sec:end_end}, the numerical scheme in IPNC etc. is the same as in Sec. \ref{sec:num_smooth_dis}.

\begin{figure}[!htb]
\centering
\subfloat[$\rho$]{
\begin{minipage}[c]{0.3\textwidth}
\includegraphics[width=\textwidth, height=0.75\textwidth]{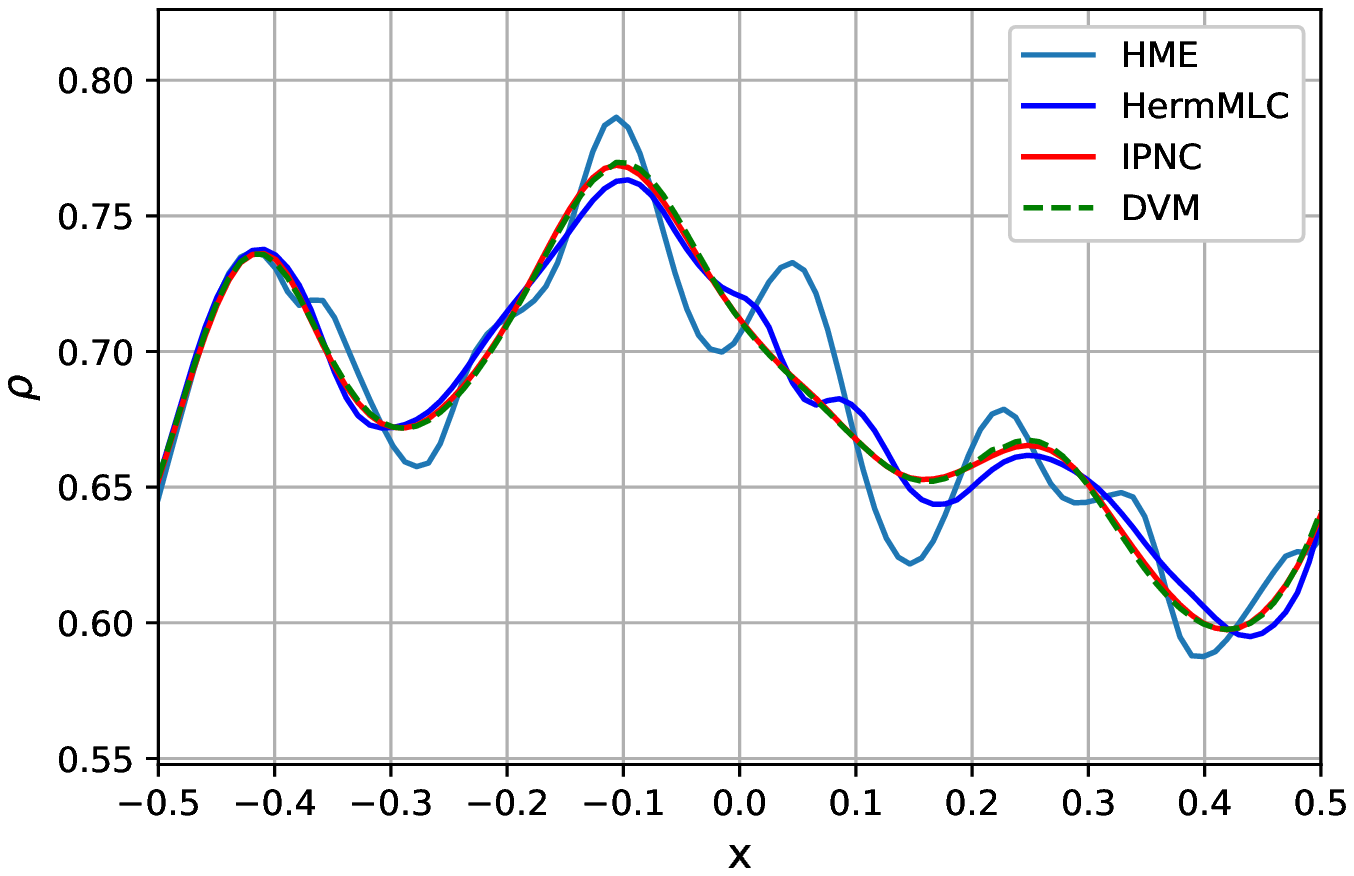}
\includegraphics[width=\textwidth, height=0.75\textwidth]{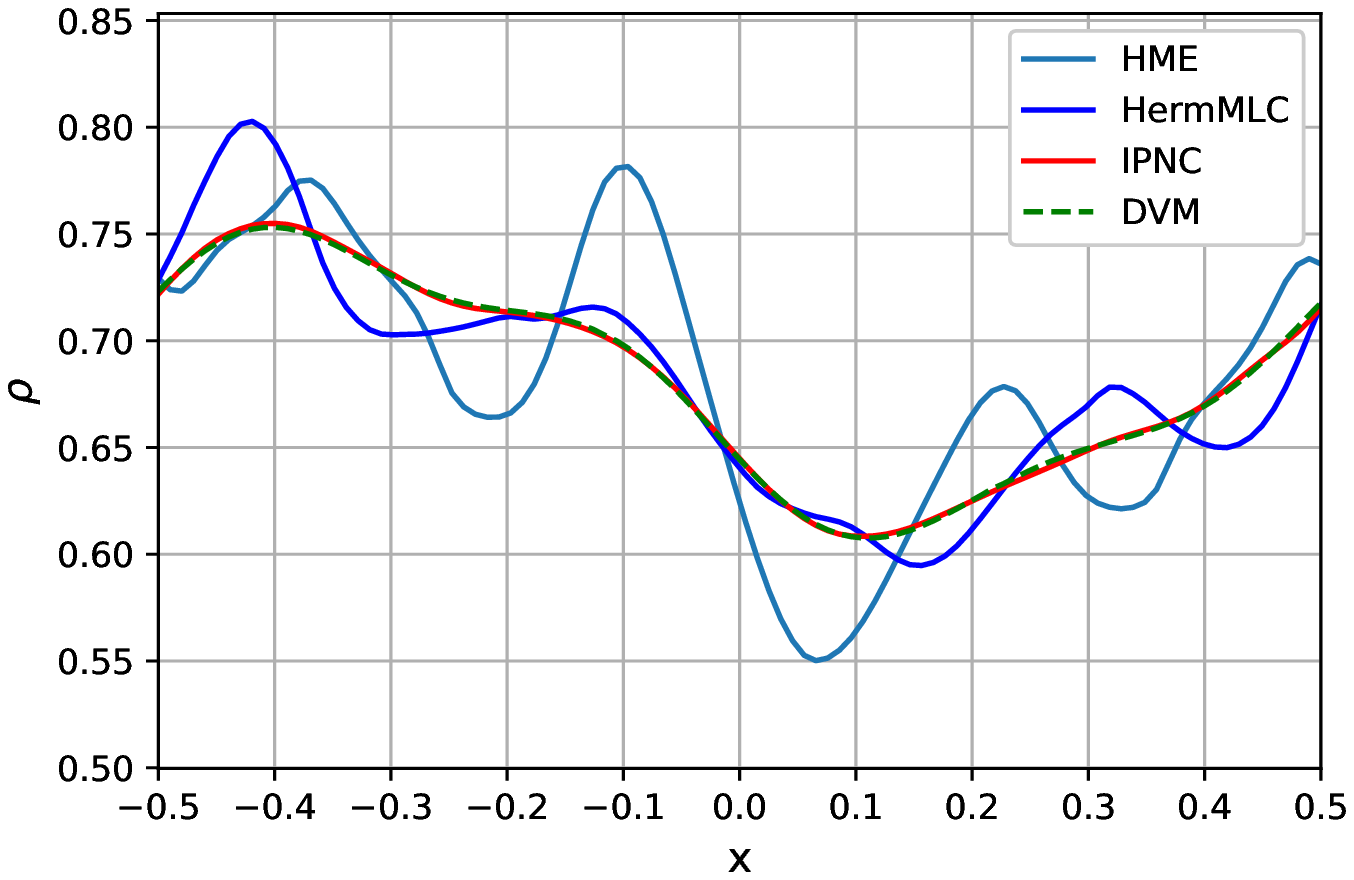}
\includegraphics[width=\textwidth, height=0.75\textwidth]{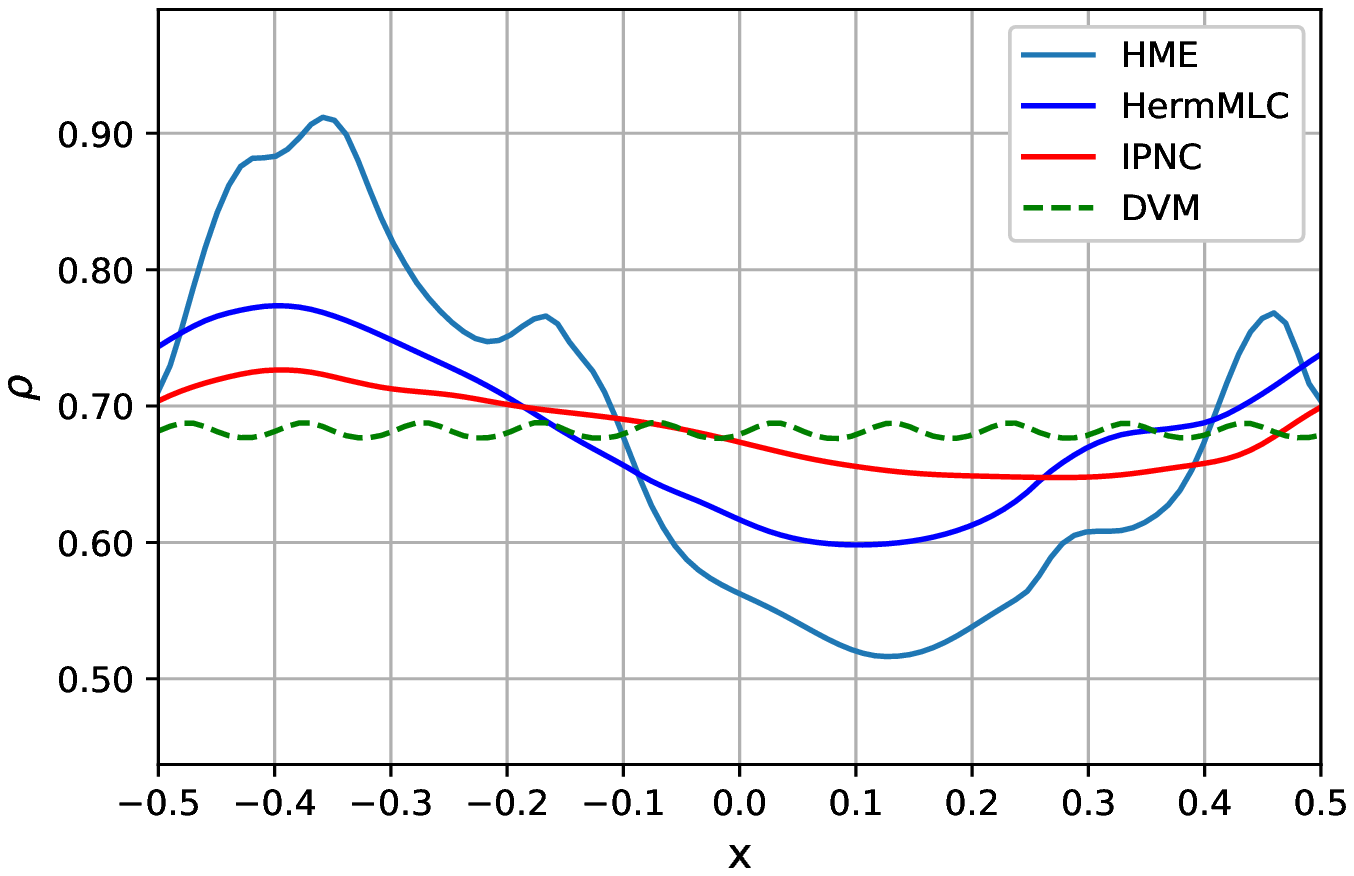}
\end{minipage}
}\quad
\subfloat[$u$]{
\begin{minipage}[c]{0.3\textwidth}
\includegraphics[width=\textwidth, height=0.75\textwidth]{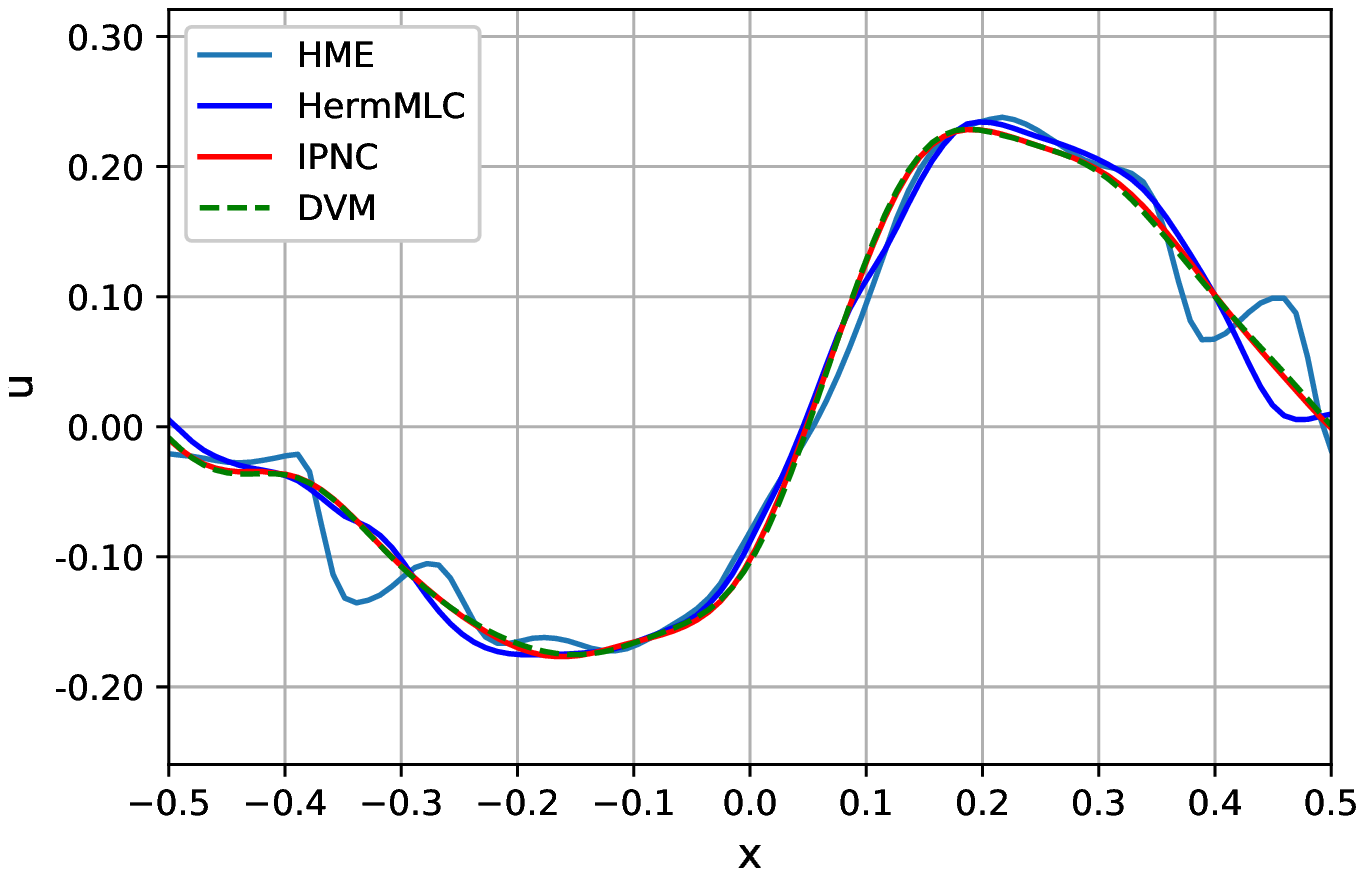}
\includegraphics[width=\textwidth, height=0.75\textwidth]{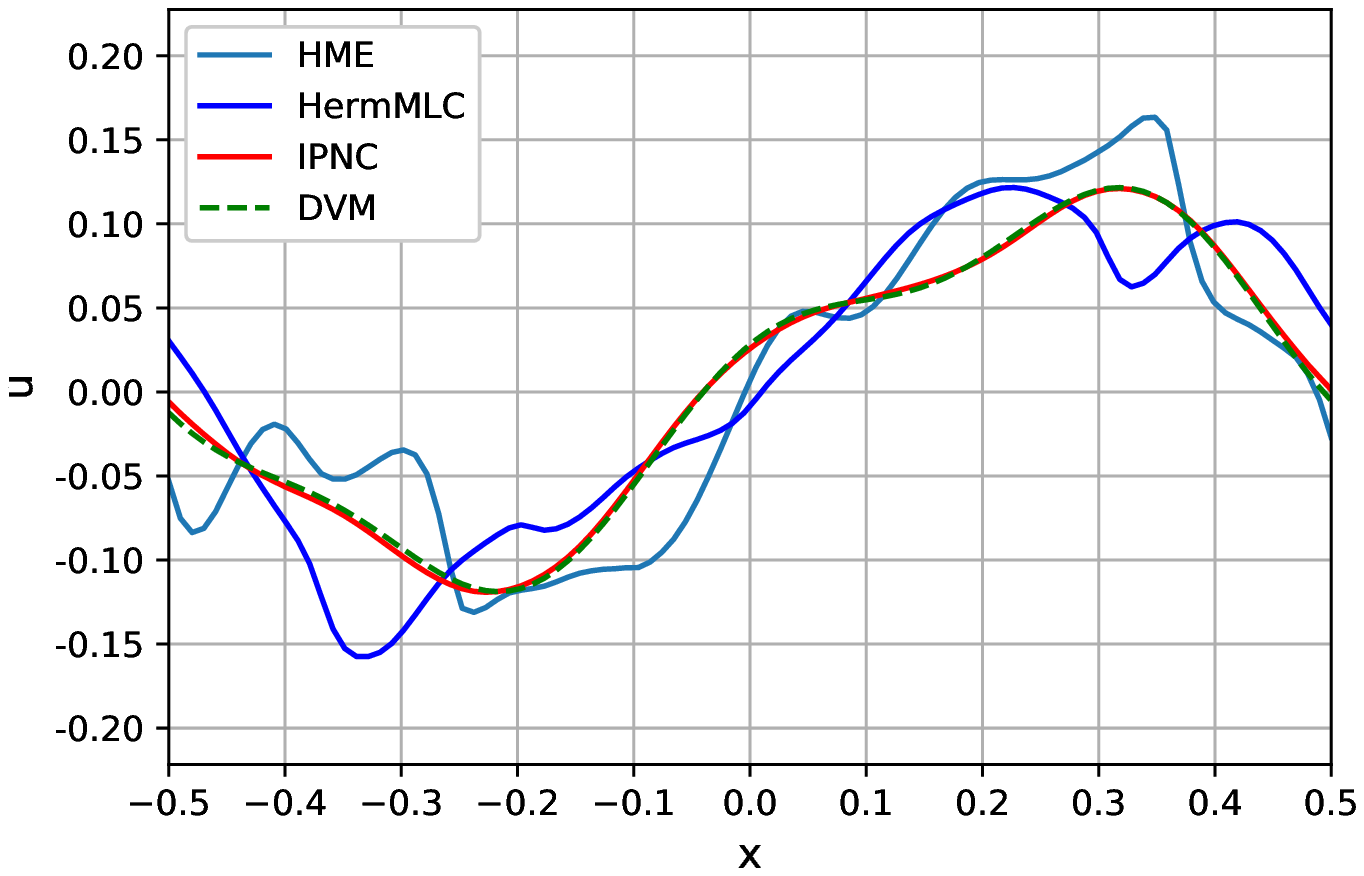}
\includegraphics[width=\textwidth, height=0.75\textwidth]{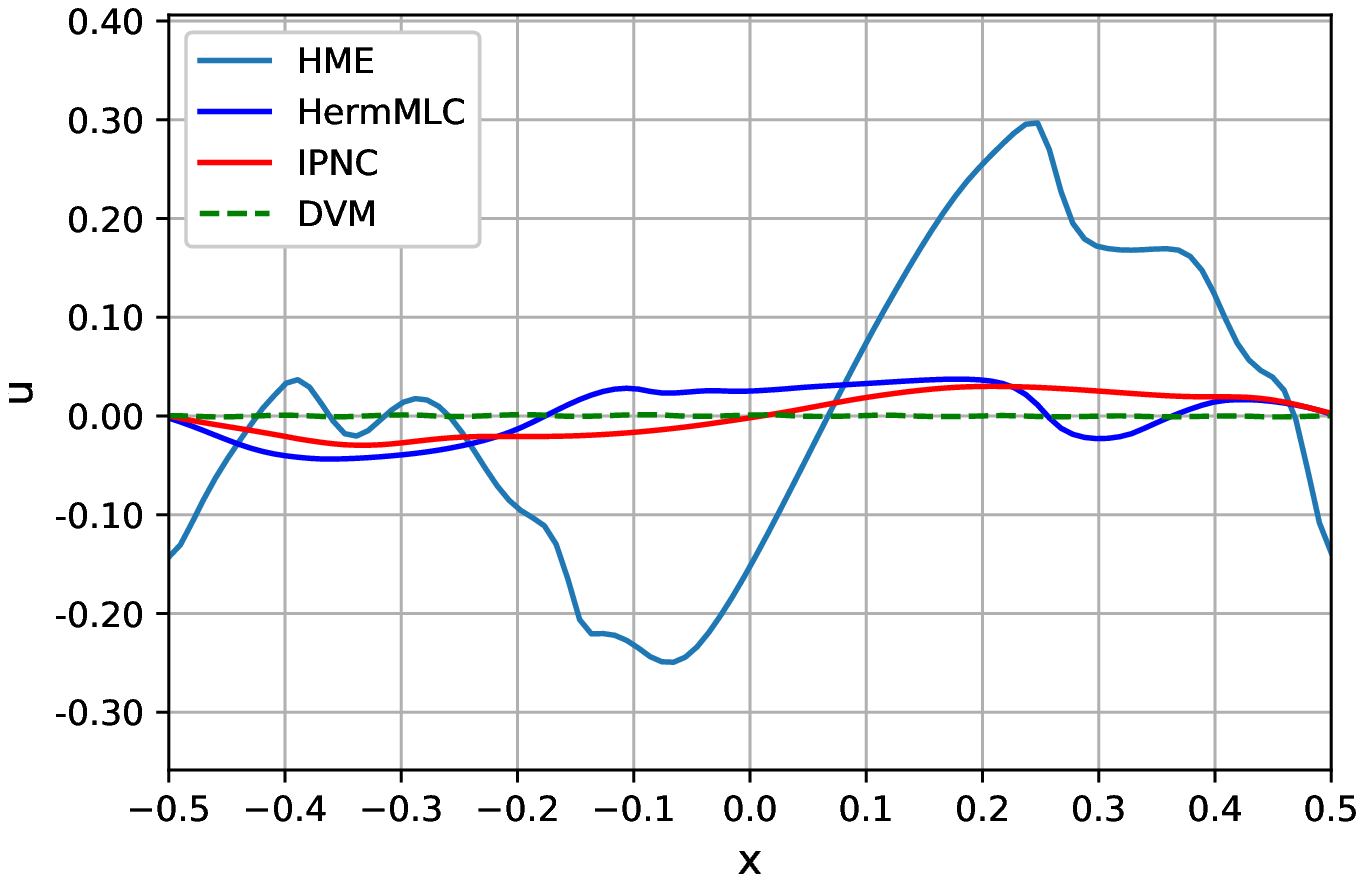}
\end{minipage}
}\quad
\subfloat[$\theta$]{
\begin{minipage}[c]{0.3\textwidth}
\includegraphics[width=\textwidth, height=0.75\textwidth]{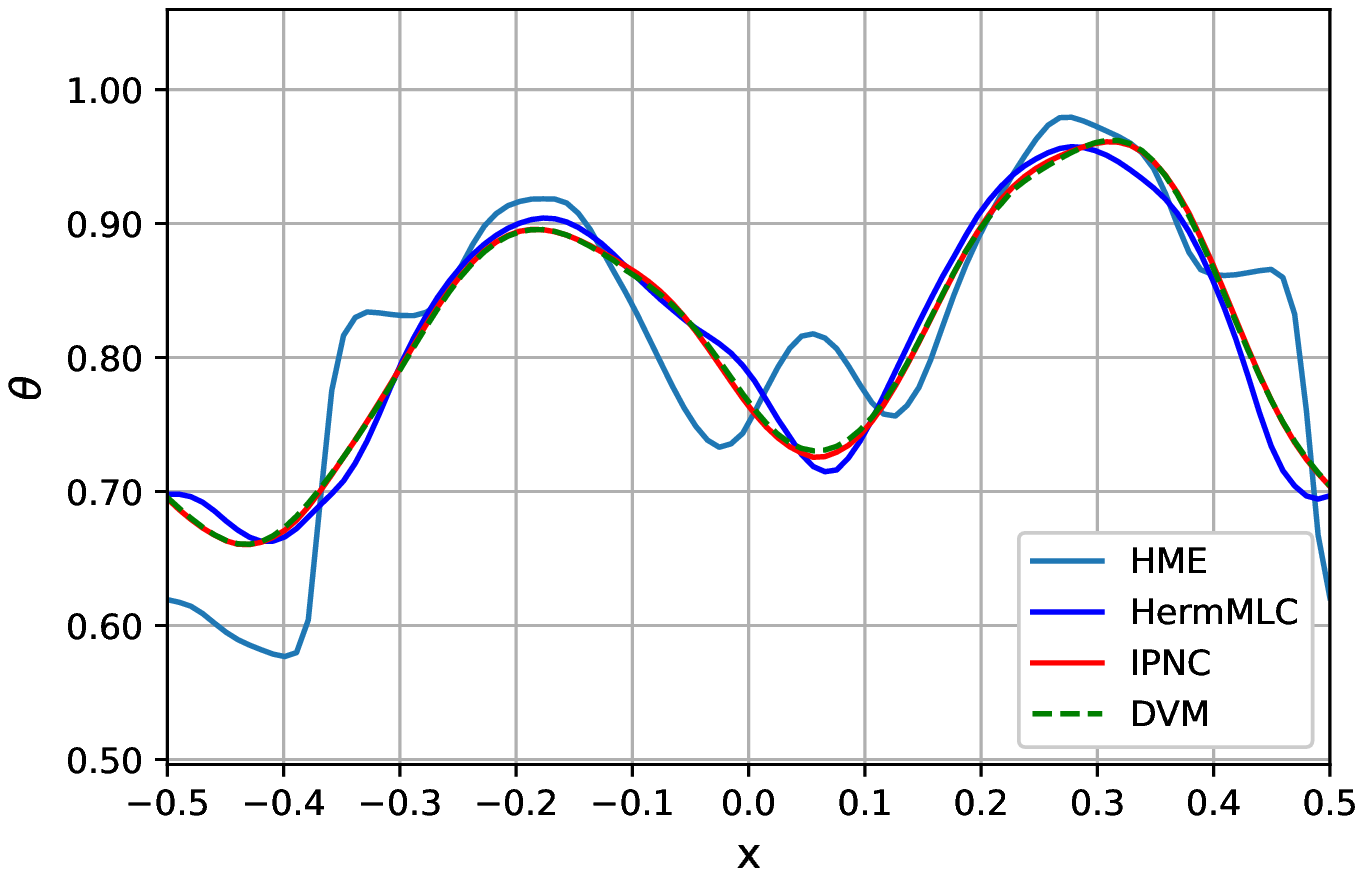}
\includegraphics[width=\textwidth, height=0.75\textwidth]{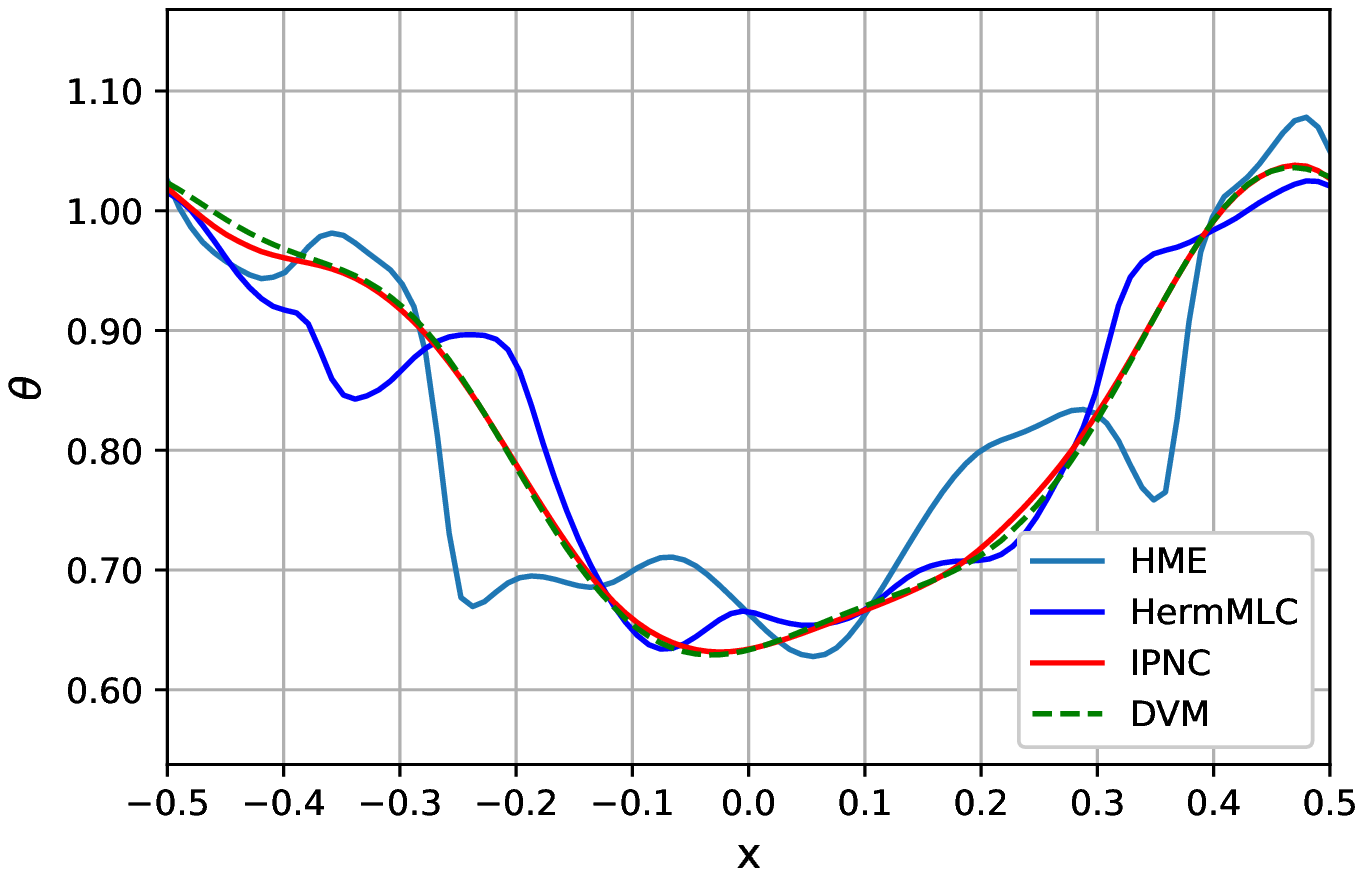}
\includegraphics[width=\textwidth, height=0.75\textwidth]{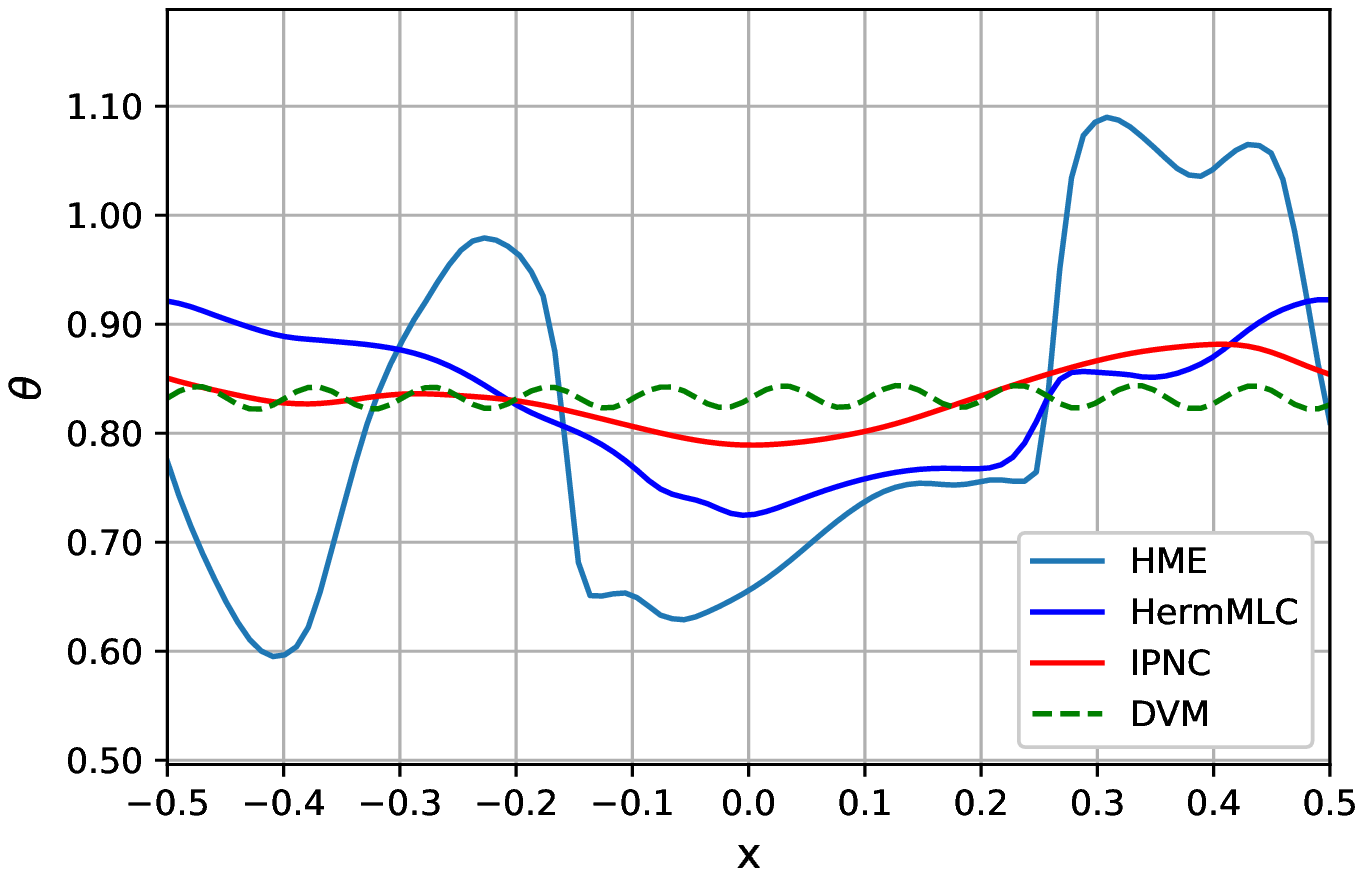}
\end{minipage}
}
\caption{(Sec. \ref{sec:num_discon_dis})  Density $\rho$, macroscopic velocity $u$, and temperature $\theta$ at time $t = 0.1, 0.2$ and $1$ for the discontinuous initial condition problem. The three rows are at $t = 0.1, 0.2$ and $1$, respectively. Here the blue line is obtained by HME, the black line is obtained by HermMLC, the red line is obtained by IPNC, and the green dashed line is the reference solution obtained by DVM.}
\label{fig:ex2_solution}
\end{figure}

\begin{figure}[!htb]
\centering
\subfloat[error quarterlies]{
\includegraphics[width=0.35\textwidth]{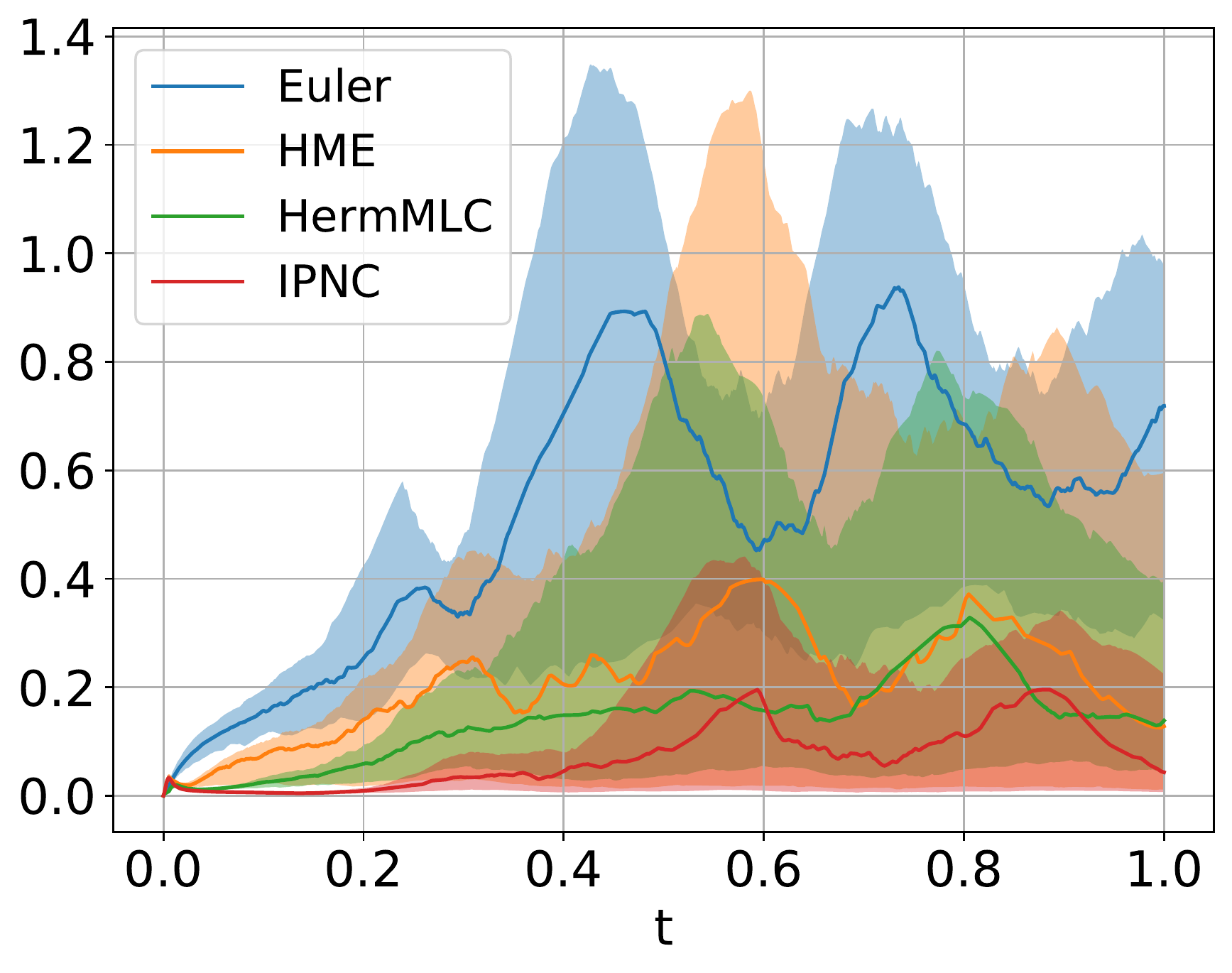}
\label{fig:ex2_error_1_1000}} \qquad
\subfloat[different $\Kn$ at $t=0.1$]{
\includegraphics[width=0.35\textwidth]{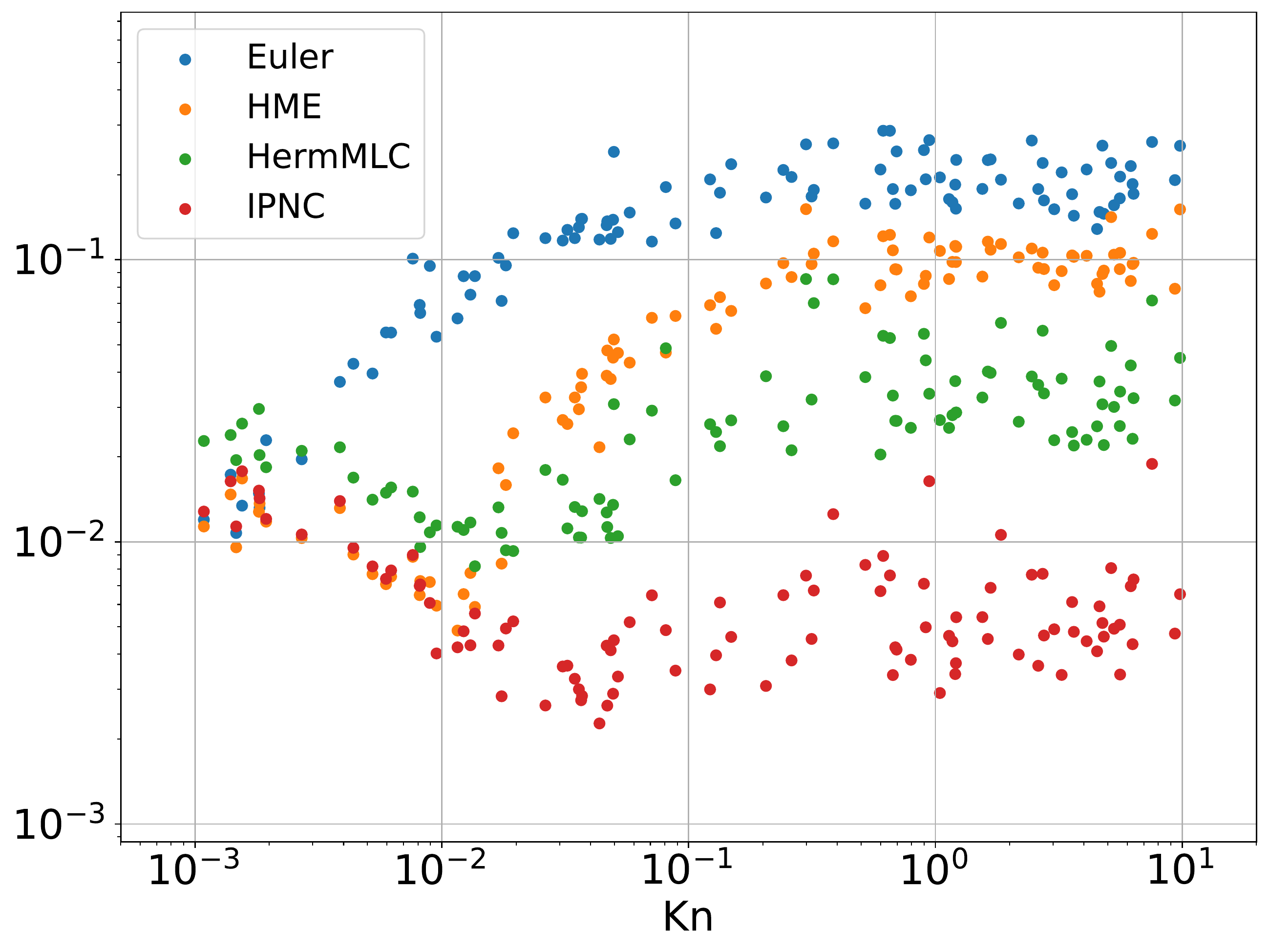}
\label{fig:ex2_error_2_100}}
\\
\subfloat[different $\Kn$ at $t=0.2$]{
\includegraphics[width=0.35\textwidth]{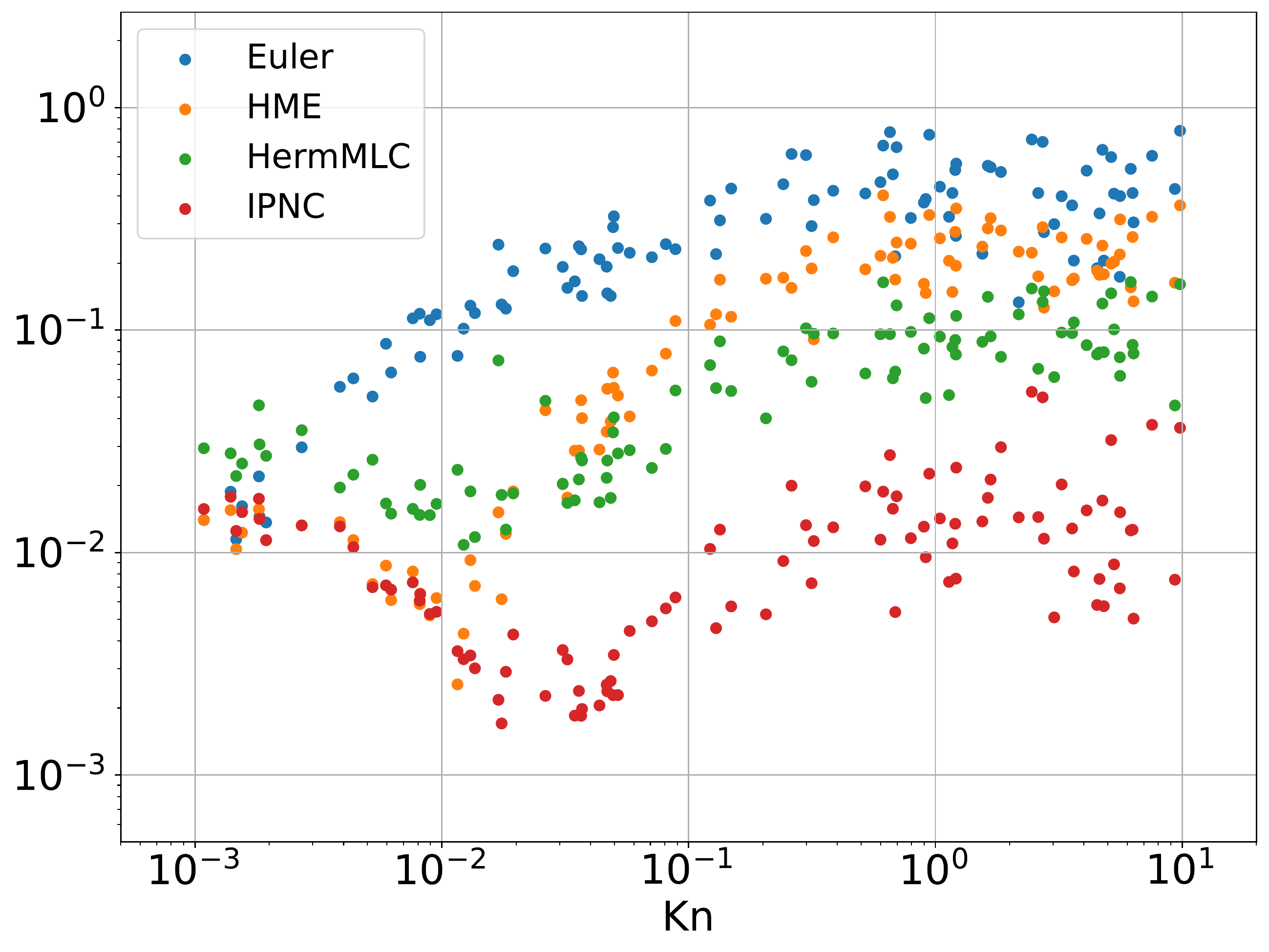}
\label{fig:ex2_error_2_200}}\qquad 
\subfloat[different $\Kn$ at $t=1$]{
\includegraphics[width=0.35\textwidth]{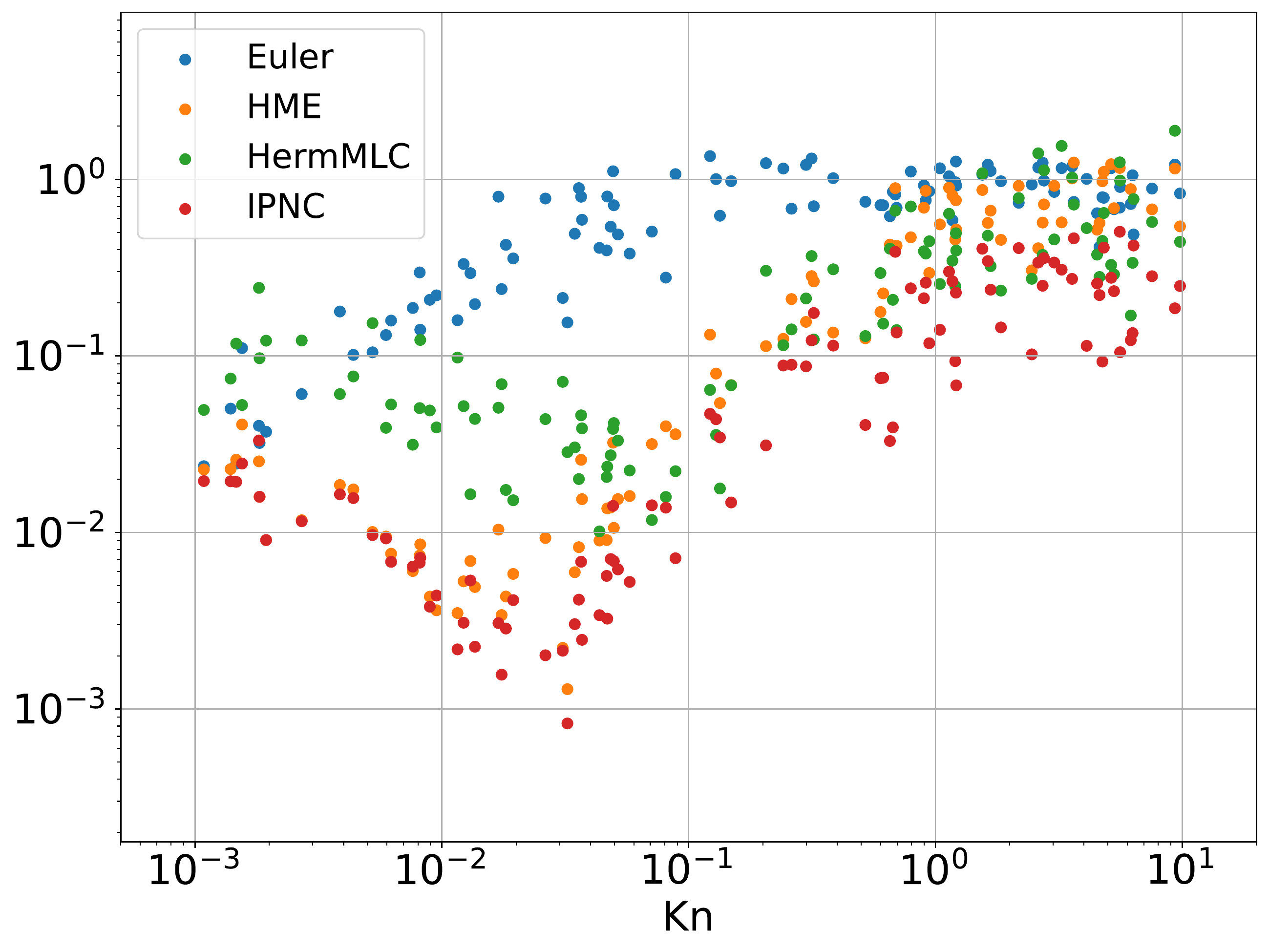}
\label{fig:ex2_error_2_1000}}
\caption{(Sec. \ref{sec:num_discon_dis}) (a) The time evolution of the distribution of the relative error with the $100$ initial samples \eqref{eq:ex2_ini} obtained by the different methods, where the $x$-axis is time $t$ and the $y$-axis is the relative error \eqref{eq:macro_error}. Here, the solid line indicates the median, and the translucent region is the area between the 1/4 and 3/4 quarterlies. (b-d) The relative error \eqref{eq:macro_error} at time $t = 0.1, 0.2$  and $1$, respectively, of the $100$ initial samples with different $\Kn$ and different initial conditions \eqref{eq:ex2_ini}. The $x$-axis is the Knudsen number and the $y$-axis is the relative error \eqref{eq:macro_error}. }
\label{fig:mix_solution_error_200}
\end{figure}

Figure \ref{fig:ex2_solution} shows the behavior of density $\rho$, macroscopic velocity $u$ and temperature $\theta$ at time $t = 0.1, 0.2$ and $1$ for one sample from the testing data, where the numerical solutions by HermMLC, HME, and the reference solution by DVM are all plotted. The initial condition corresponding to this specific case is given in Appendix \ref{app:supp}. We can also see that IPNC matches DVM better than HermMLC and HME. The similar distribution of the relative error \eqref{eq:macro_error} for different initial conditions with different methods is shown in Figure \ref{fig:ex2_error_1_1000}. We can see that both the mean and variance of the relative errors of IPNC are much smaller than those of all the other methods. The stability of IPNC with respect to the Knudsen number $\Kn$ at $t = 0.1, 0.2$ and $1$ is shown from Figure \ref{fig:ex2_error_2_100} to \ref{fig:ex2_error_2_1000}, respectively, where it is apparent that the behavior of IPNC is more stable compared to other methods. The average values of the relative errors for the same five different $\Kn$ as in Sec. \ref{sec:num_smooth_dis} are presented in Table \ref{tab:mixH}. We can see that IPNC is the most accurate method for this discontinuous problem.

\begin{table}[ht]
\centering
\renewcommand\arraystretch{1.5}
\footnotesize
\begin{tabular}{l|lllll}
$ t = 0.1$ & Kn      & $0.01$          & $0.1$           & $1.0$            & $10$             \\ \hline
           & Euler   & $7.89$          & $16.89$         & $19.62$          & $19.97 $         \\
           & HME     & $0.72 $         & $6.66 $         & $9.83$           & $10.25 $         \\
           & HermMLC & $0.90$          & $2.60$          & $3.88 $          & $4.05$           \\
           & IPNC    & $\mathbf{0.43}$ & $\mathbf{0.49}$ & $\mathbf{0.66}$  & $\mathbf{0.72}$  \\ \hline
$t = 0.2$  & Kn      & $0.01$          & $0.1$           & $1.0$            & $10$             \\ \hline
           & Euler   & $10.81$         & $31.79$         & $42.78$          & $44.35$          \\
           & HME     & $0.58$          & $10.24$         & $21.26$          & $23.13$          \\
           & HermMLC & $1.47$          & $4.96$          & $9.55$           & $10.39$          \\
           & IPNC    & $\mathbf{0.35}$ & $\mathbf{0.74}$ & $\mathbf{1.56} $ & $\mathbf{1.62}$  \\ \hline
$t = 1$    & Kn      & $0.01$          & $0.1$           & $1.0$            & $10$             \\ \hline
           & Euler   & $26.54$         & $80.00$         & $85.56$          & $85.68$          \\
           & HME     & $0.64$          & $7.83$          & $-$              & $-$              \\
           & HermMLC & $5.44$          & $6.45$          & $47.95$          & $63.76$          \\
           & IPNC    & $\mathbf{0.42}$ & $\mathbf{3.34}$ & $\mathbf{18.64}$ & $\mathbf{23.04}$
\end{tabular}
\caption{(Sec. \ref{sec:num_discon_dis})  Average of the relative error \eqref{eq:macro_error} with $100$ initial samples obtained by different numerical methods at time $t = 0.1, 0.2$ and $1$ with the discontinuous initial condition \eqref{eq:ex2_ini}.}
\label{tab:mixH}
\end{table}

\begin{figure}[!htb]
\centering
\subfloat[$\rho$]{
\begin{minipage}[c]{0.3\textwidth}
\includegraphics[width=\textwidth, height=0.75\textwidth]{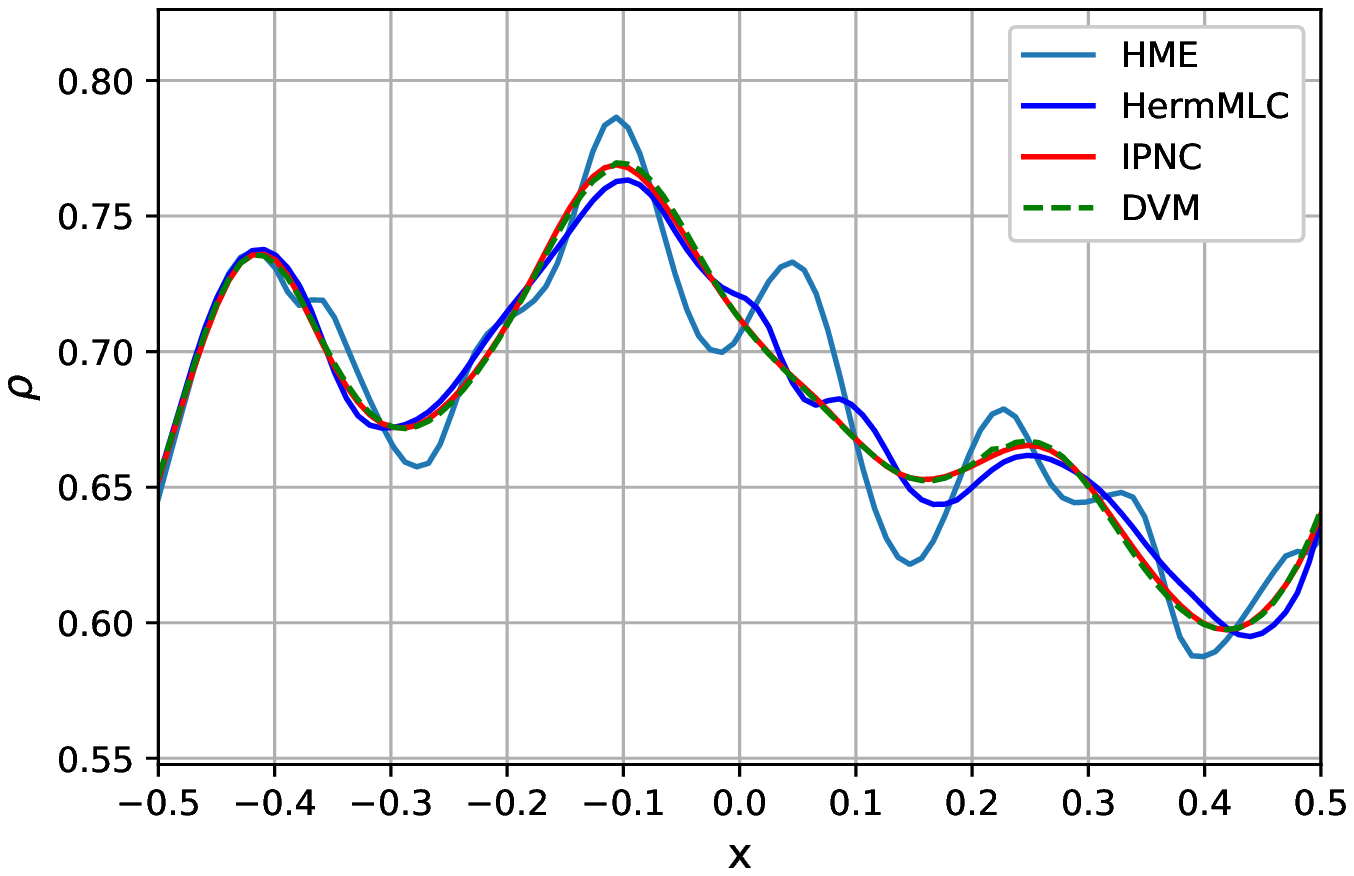}
\end{minipage}
}\quad
\subfloat[$u$]{
\begin{minipage}[c]{0.3\textwidth}
\includegraphics[width=\textwidth, height=0.75\textwidth]{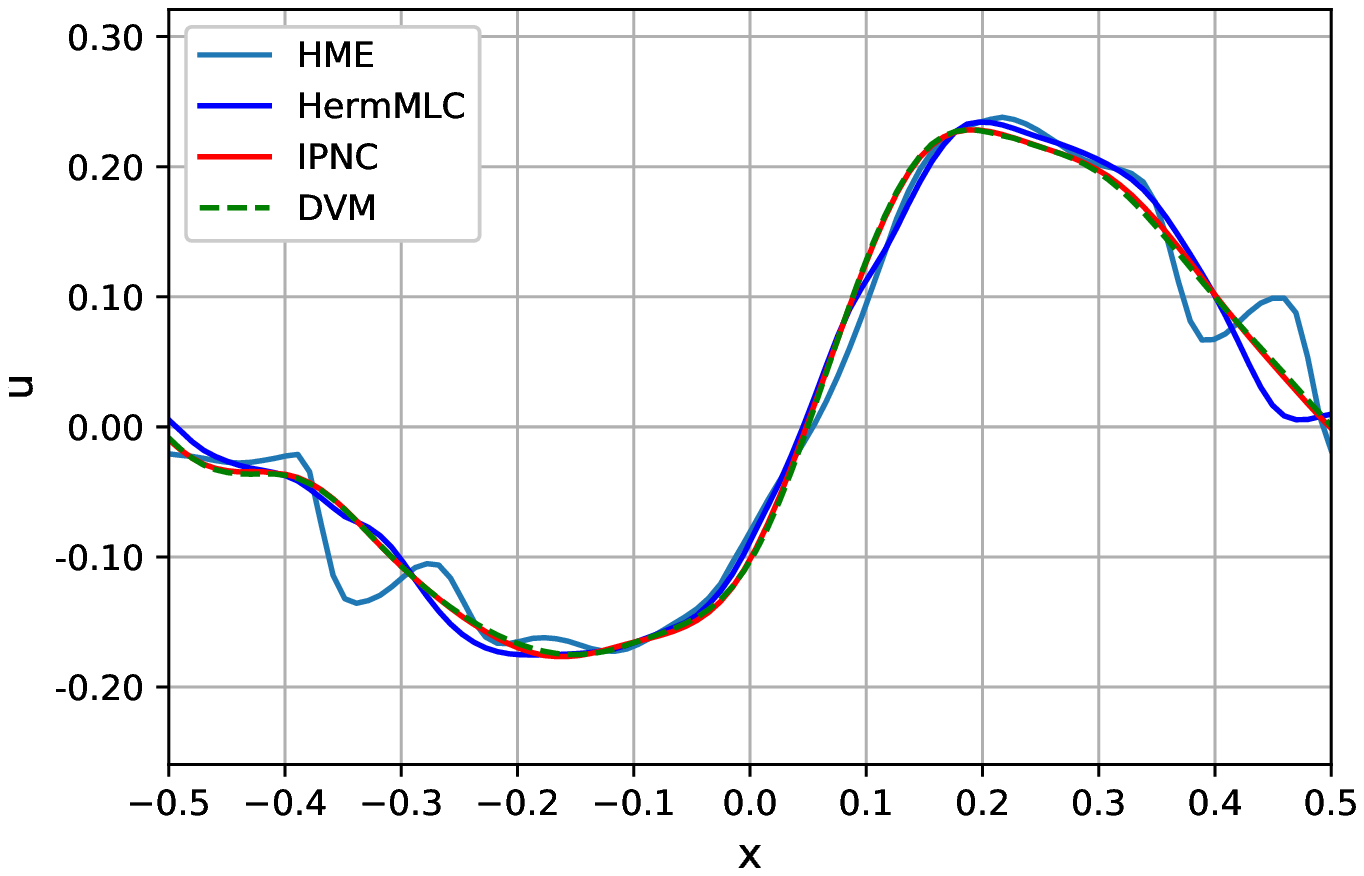}
\end{minipage}
}\quad
\subfloat[$\theta$]{
\begin{minipage}[c]{0.3\textwidth}
\includegraphics[width=\textwidth, height=0.75\textwidth]{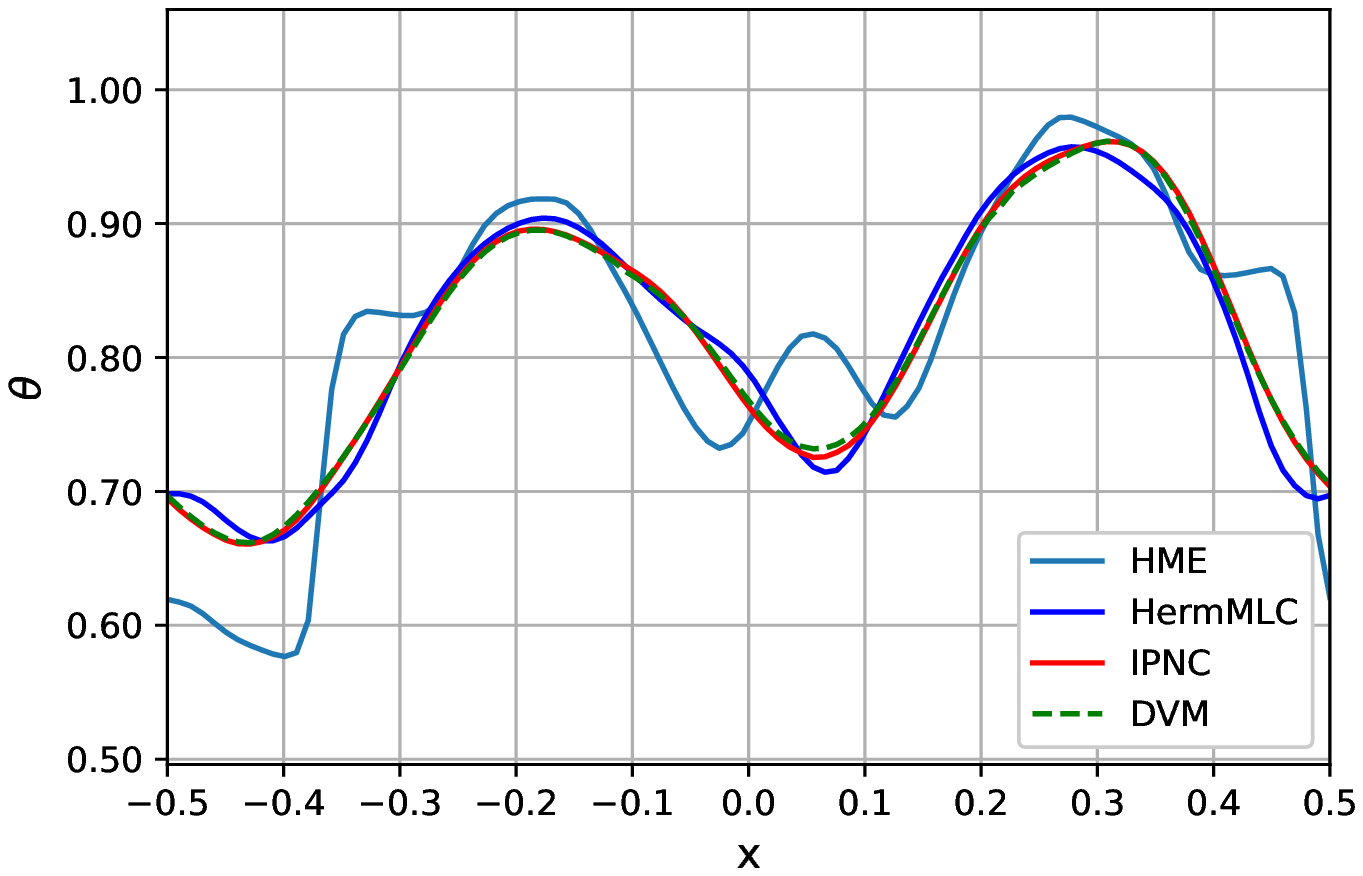}
\end{minipage}
}
\caption { (Sec. \ref{sec:num_discon_dis}: Generalization of $\Kn$) Density $\rho$, macroscopic velocity $u$, and temperature $\theta$ at $ t= 0.1$ with $\Kn = 10$.  Here, the blue line is got by HME, the black line is got by HermMLC, the red line is got by IPNC, and the green dashed line is the reference solution got by DVM.}
\label{fig:ex2_ger_Kn_sol}
\end{figure}

\paragraph{Generalization of Knudsen number} To test the generalization ability of IPNC on Knudsen number for the discontinuous problem, the similar tests as in Sec. \ref{sec:num_smooth_dis} are studied, where the region of $\Kn$ in the training set is changed to $\Kn \in [0.1, 1.0]$ and keep the remains the same. The density $\rho$, macroscopic velocity $u$, and temperature $\theta$ at $ t = 0.1$ for the same initial data as in Figure \ref{fig:ex2_solution} with $\Kn = 10$ is plotted in Figure \ref{fig:ex2_ger_Kn_sol}. A similar average of the relative errors for the different Knudsen numbers at time $t = 0.1, 0.2$ and $1$ with the changed training region of $\Kn$ is shown in Table \ref{tab:mixH_ger_Kn}. These results show that for the discontinuous problem, the numerical solution is still quite close to the reference solution by DVM, and the average relative error is not increasing with the change of the training set for the Knudsen number. The generalization ability of the IPNC on the Knudsen number for the discontinuous is validated.

\begin{table}[ht]
\centering
\renewcommand\arraystretch{1.5}
\footnotesize
\begin{tabular}{l|lllll}
$ t = 0.1 $ & Kn         & $0.01$  & $0.1$   & $1.0$   & $10$    \\ \hline
            & Kn         & $0.43$  & $0.49$  & $0.66$  & $0.72$  \\
            & $Kn_{\rm lim}$ & $0.71 $ & $0.48$  & $0.67$  & $1.24$  \\ \hline
$ t = 0.2$  & Kn         & $0.01$  & $0.1$   & $1.0$   & $10$    \\ \hline
            & Kn         & $0.35$  & $0.74$  & $1.56$  & $1.62$  \\
            & $\Kn_{\rm lim}$ & $0.97$  & $0.76 $ & $1.55 $ & $2.24 $ \\ \hline
$ t = 1$    & Kn         & $0.01$  & $0.1$   & $1.0$   & $10$    \\ \hline
            & Kn         & $0.42$  & $3.34$  & $18.64$ & $23.04$ \\
            & $\Kn_{\rm lim}$ & $1.98$  & $2.90$  & $18.37$ & $16.98$
\end{tabular}
\caption{
(Sec. \ref{sec:num_discon_dis}: Generalization of $\Kn$)  Average of the relative error \eqref{eq:macro_error} got by IPNC with the $100$ samples of the discontinuous initial condition \eqref{eq:ex2_ini} at time $t = 0.1, 0.2$ and $1$ with different training Knudsen number. Here, $\Kn$ is  corresponding to the training set $\Kn \in [0.001, 10]$ and $\Kn_{\rm lim}$ is corresponding to the training set $\Kn \in [0.1, 1]$.
}
\label{tab:mixH_ger_Kn}
\end{table}

\begin{figure}[!htb]
\centering
\subfloat[$\rho$]{
\begin{minipage}[c]{0.3\textwidth}
\includegraphics[width=\textwidth, height=0.75\textwidth]{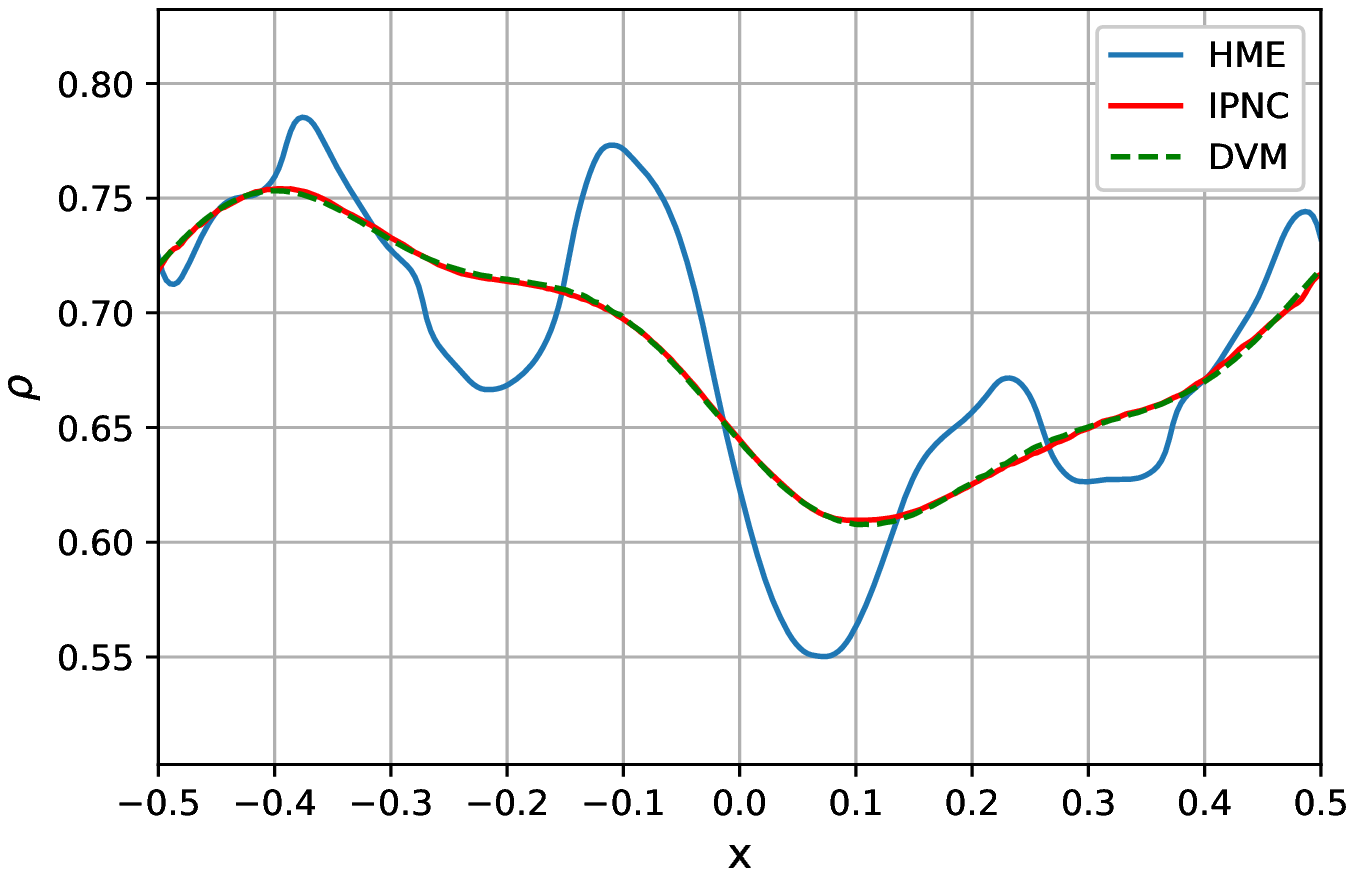}
\end{minipage}
}\quad
\subfloat[$u$]{
\begin{minipage}[c]{0.3\textwidth}
\includegraphics[width=\textwidth, height=0.75\textwidth]{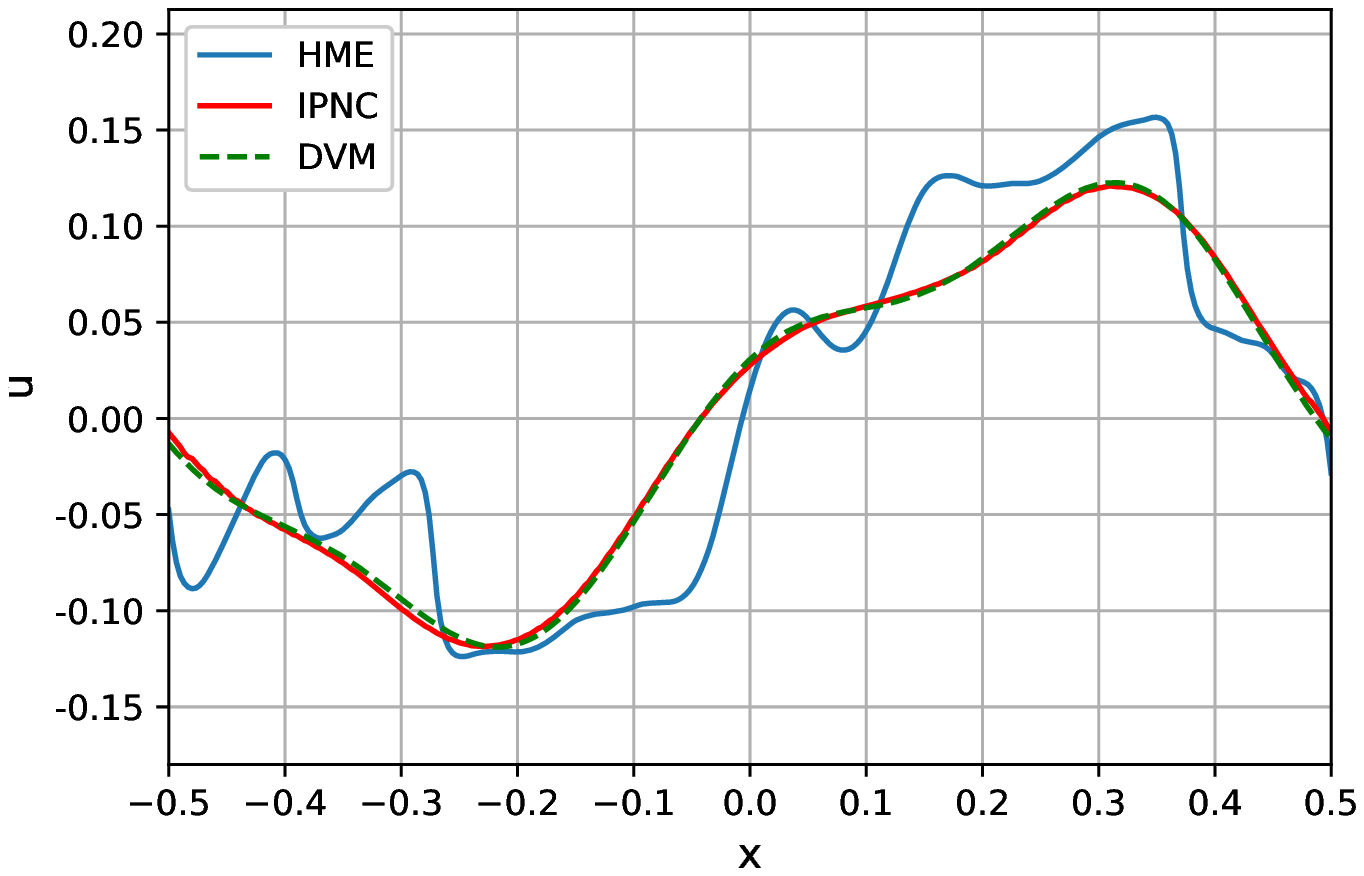}
\end{minipage}
}\quad
\subfloat[$\theta$]{
\begin{minipage}[c]{0.3\textwidth}
\includegraphics[width=\textwidth, height=0.75\textwidth]{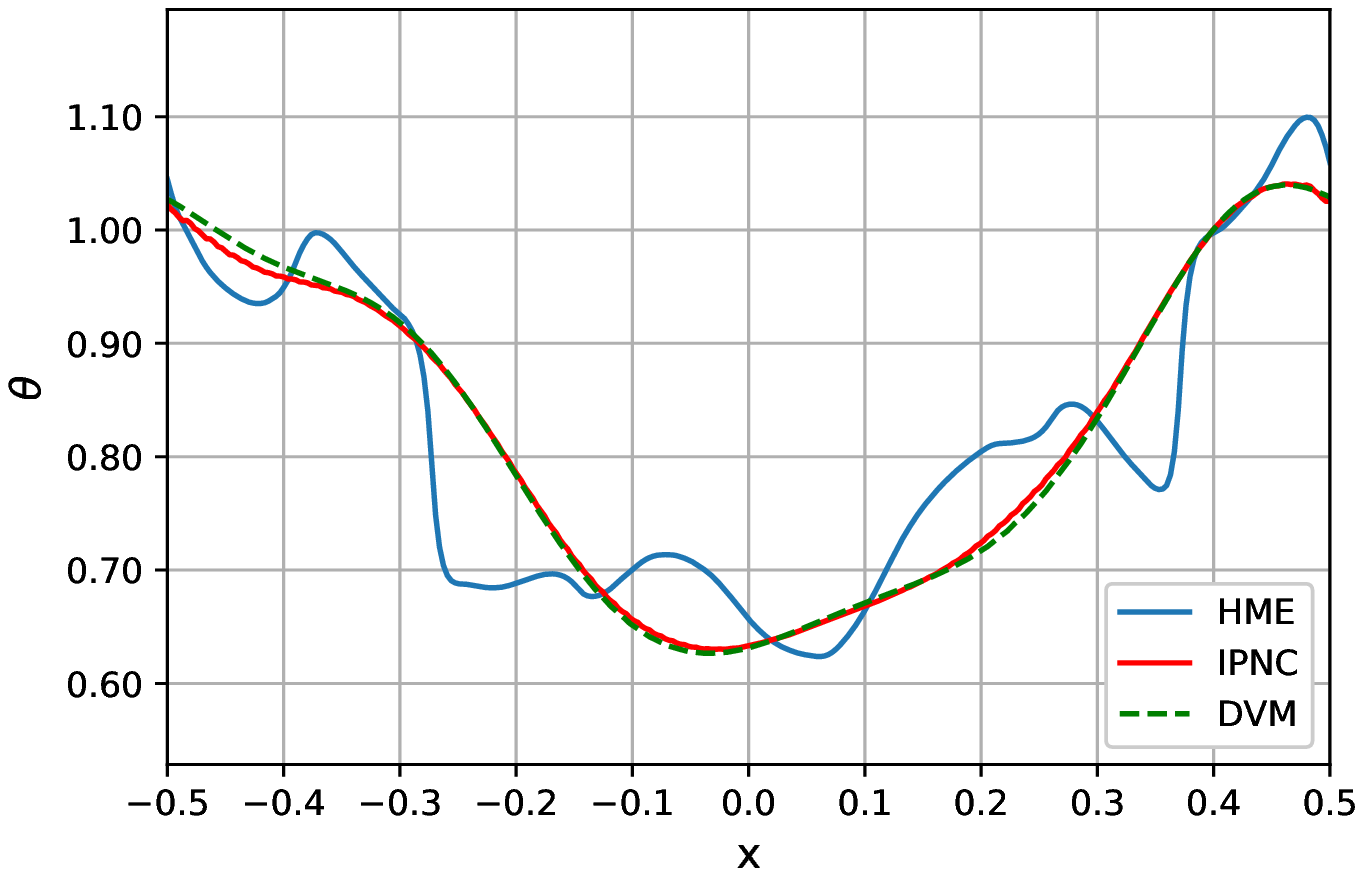}
\end{minipage}
}
\caption { (Sec. \ref{sec:num_discon_dis}: Generalization of mesh size) Density $\rho$, macroscopic velocity $u$, and temperature $\theta$ at $ t= 0.2$ with mesh size $N_x = 300$.  Here, the blue line is obtained by HME, the black line is obtained by HermMLC, the red line is got by IPNC, and the green dashed line is the reference solution obtained by DVM.}
\label{fig:ex2_ger_dx_sol}
\end{figure}

\paragraph{Generalization of the mesh size} To show the generalization ability of IPNC on mesh size for the discontinuous problem, the tests with different mesh sizes are studied. Similarly, the mesh size $N_x = 100$ is utilized in the training data set, while $N_x = 100, 200, 300$ and $400$ is adopted for the testing data set. The numerical solution with $N_x = 300$ at $ t = 0.2 $ for the discontinuous problem is shown in Figure \ref{fig:ex2_ger_dx_sol} and the average relative error \eqref{eq:macro_error} obtained by IPNC with the $100$ samples at time $t = 0.2$ with different mesh size is calculated in Table \ref{tab:ex2_ger_dx_error}. The numerical results show behavior similar to those of the continuous problem, which validates the generalization ability of IPNC on mesh size for the discontinuous problem. However, when the mesh size is increasing to $N_x = 400$,
there are two failures for the $100$ samples, and the error calculated from the rest samples is much larger compared to other mesh sizes. This means that there exists a limitation for IPNC on the mesh generalization by adopting only the interpolation method. 

\begin{table}[ht]
\centering
\renewcommand\arraystretch{1.5}
\footnotesize
% \begin{tabular}{l|l|l|l||l|l|l|l}
% Mix  & $N_x=100$ & $N_x=200$ & $N_x=300 $ & Wave & $N_x=100$ & $N_x=200$ & $N_x=300 $ \\ \hline
% HME  & $13.12$   & $14.40$   & $13.52$    & HME  & $16.33$   & $21.18$   & $23.08$    \\
% IPNC & $1.18$    & $1.16$    & $1.14$     & IPNC & $1.09$    & $1.04$    & $1.38$    
% \end{tabular}
% \caption{Network trained on Mix with nx=100, prediction error at T=0.2s for different nx.}
\begin{tabular}{l|llll}
     & $N_x = 100$ & $N_x = 200$ & $N_x = 300$ & $N_x = 400$ \\ \hline
HME  & $13.12$     & $14.40$     & $13.52$     & $15.24$     \\
IPNC & $1.18$      & $1.16$      & $1.14$      & $2.29$ ($2$ failures)        
\end{tabular}
\caption{(Sec. \ref{sec:num_discon_dis}: Generalization of  mesh size)  Average of the relative error \eqref{eq:macro_error} obtained by IPNC with the $100$ samples of the discontinuous initial condition \eqref{eq:ex2_ini} at time $t = 0.2$ with different mesh size. }
\label{tab:ex2_ger_dx_error}
\end{table}

\subsection{Sod's shock tube problem}
\label{sec:num_sod}
In this section, we test the generalization of the IPCN network trained on the mix problem to the Sod's shock tube problem without retraining. Here, the classical Sod's shock tube problem is studied, where the initial condition is chosen as 
\begin{equation}
(\rho, u, \theta) = \left\{\begin{array}{ll}
\left(1, 0, 1\right), & x \geqslant 0, \\
\left(0.125, 0, 0.8\right), & x < 0.
\end{array}\right. 
\end{equation}
In this experiment, the reference solution is obtained by DVM using a $2$nd order LF scheme with the number of grid $N_x = 400$ on $[-0.5, 0.5]$. The computational region in the microscopic velocity space is $[-10, 10]$ with $400$ grid points. The same numerical scheme as in Sec. \ref{sec:num_discon_dis} is adopted but with the grid size increased to $N_x = 400$ so as to improve the precision of the network and the numerical solution. The set $\mathcal{I}$ is chosen as $\mathcal{I}=\{0\leqslant i \leqslant 396,i\in \mathbb{Z}\}$ with $t^i=0.00025i$ and $N_s=100, B = 4$.
Here, we want to emphasize that the network in Sec. \ref{sec:num_discon_dis} is retrained with the increased number of the spatial grid $N_x = 400$. 

\begin{figure}[!htb]
\centering
\includegraphics[width=0.3\textwidth]{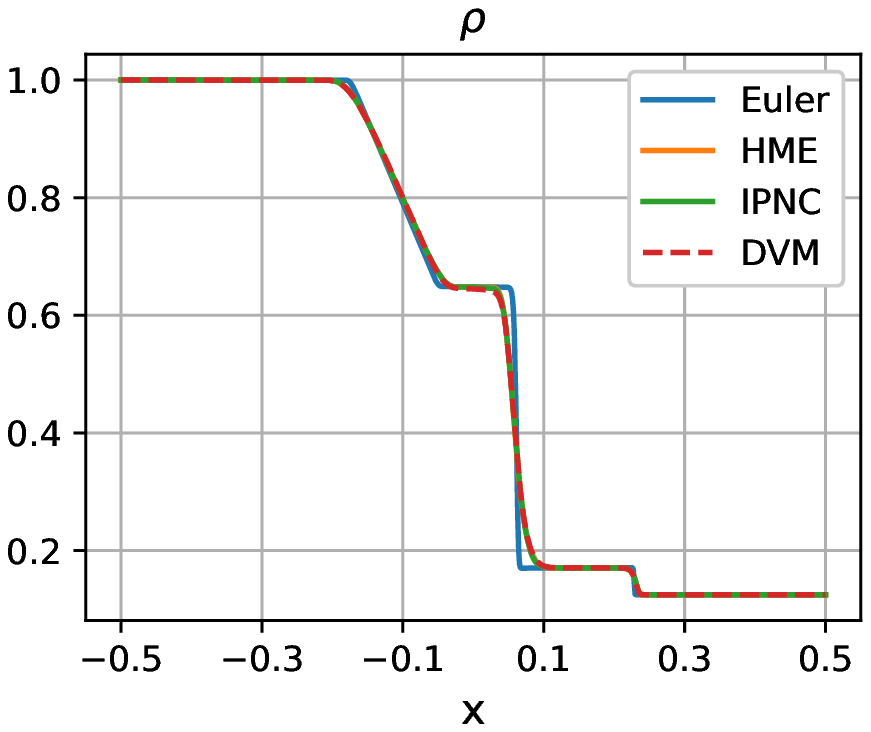}\qquad
\includegraphics[width=0.3\textwidth]{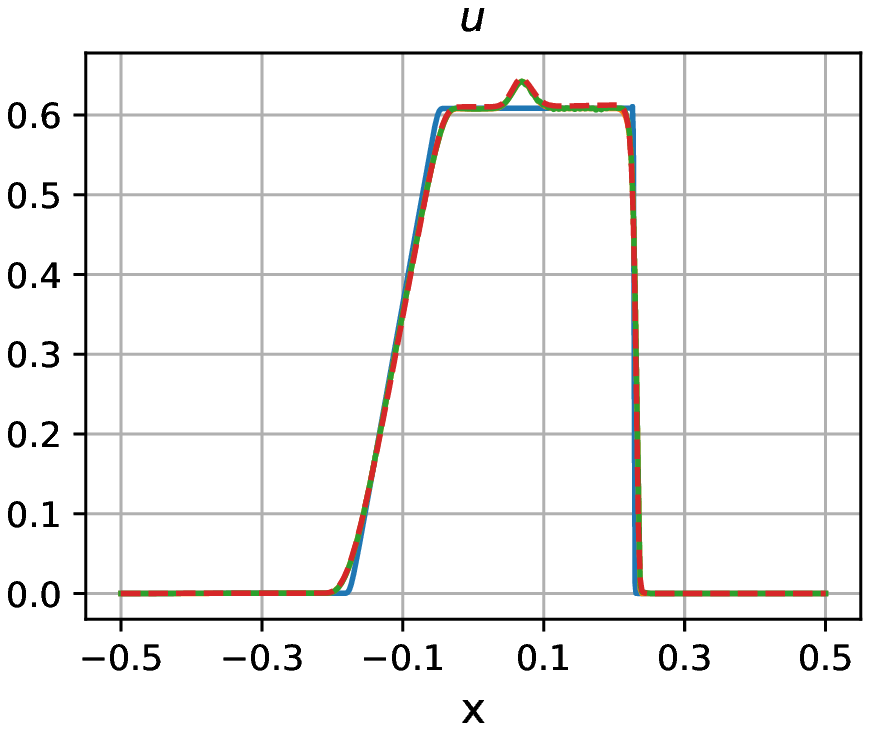}\qquad
\includegraphics[width=0.3\textwidth]{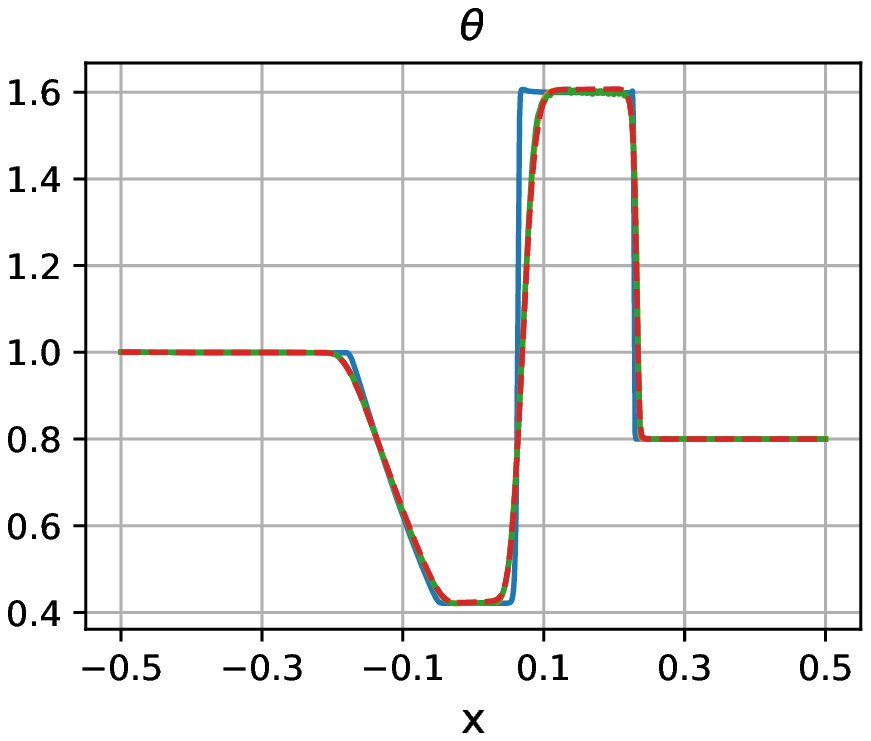}\\
\includegraphics[width=0.3\textwidth]{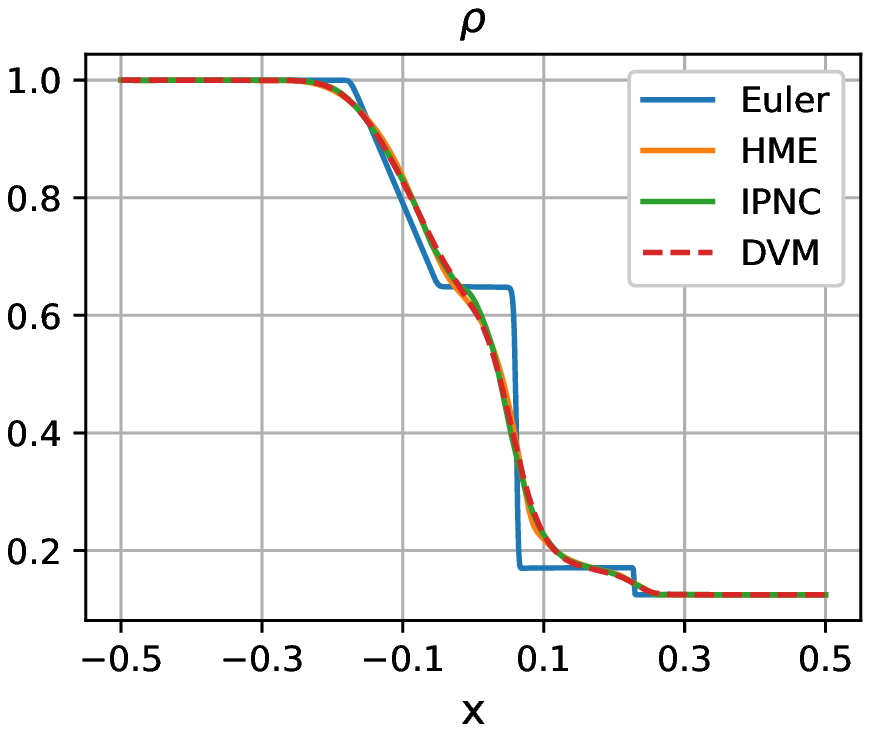} \qquad
\includegraphics[width=0.29\textwidth]{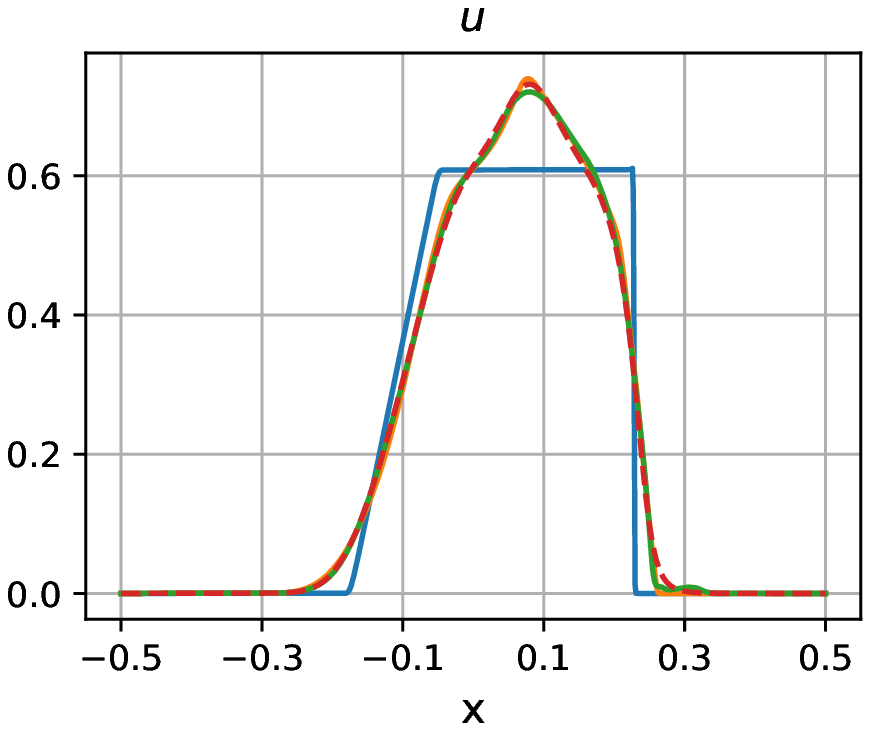} \qquad
\includegraphics[width=0.3\textwidth]{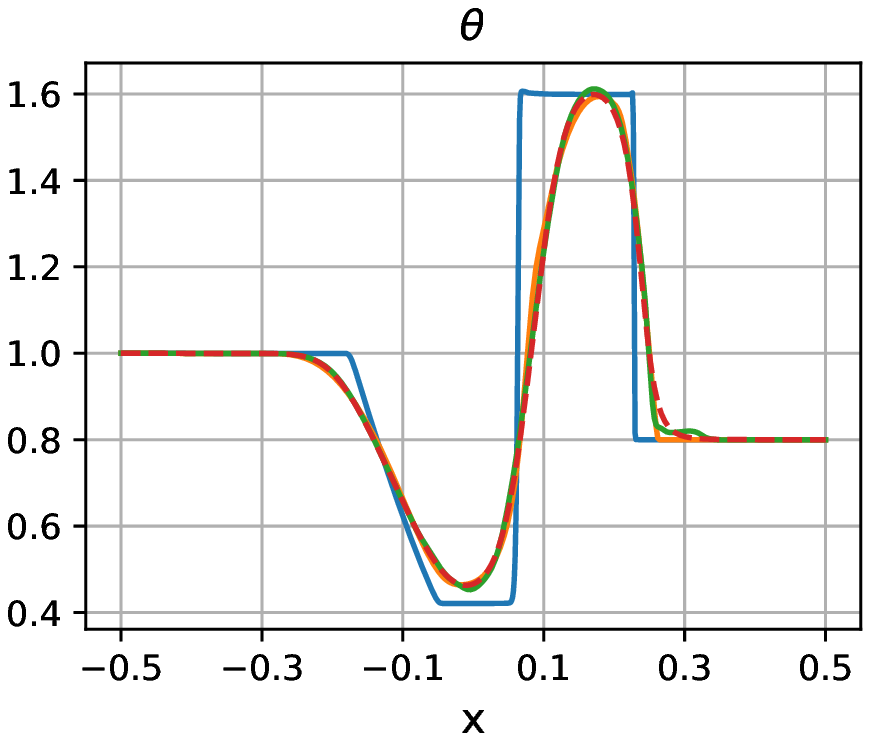} \\
\includegraphics[width=0.3\textwidth]{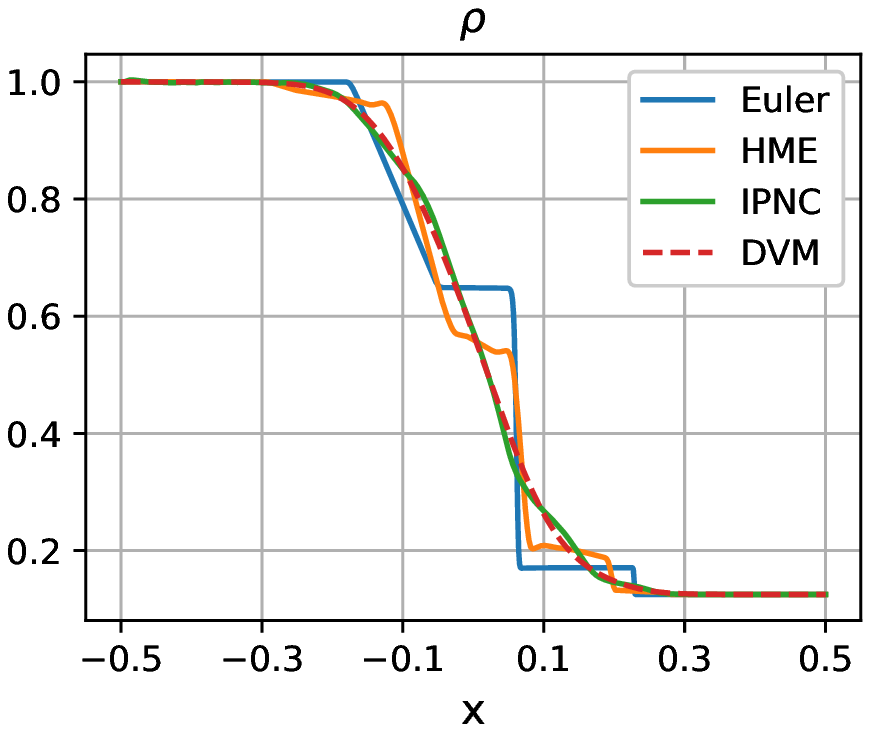}\qquad
\includegraphics[width=0.3\textwidth]{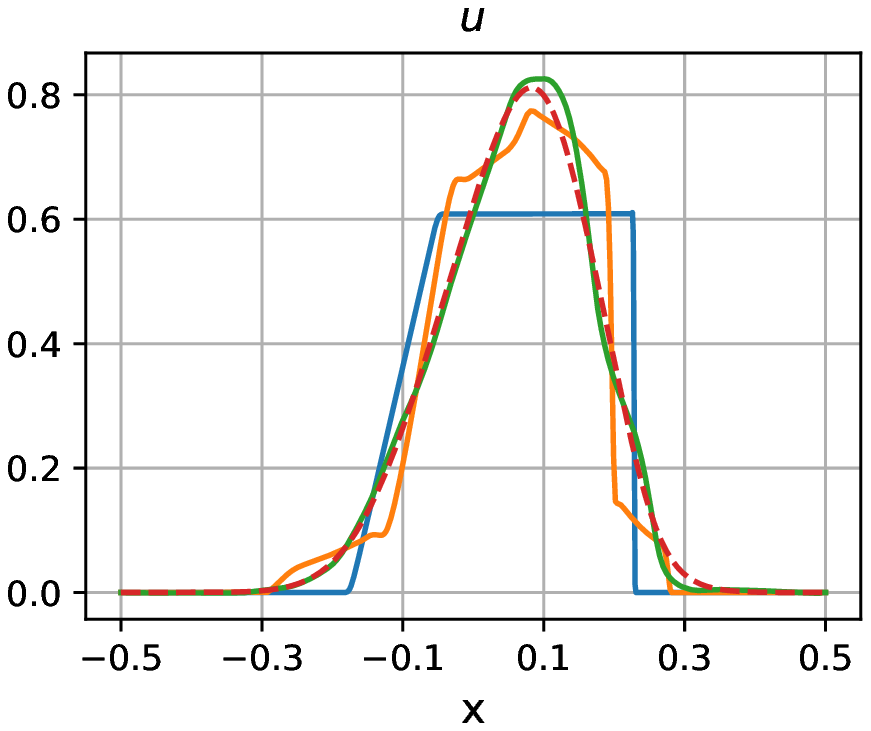}\qquad
\includegraphics[width=0.3\textwidth]{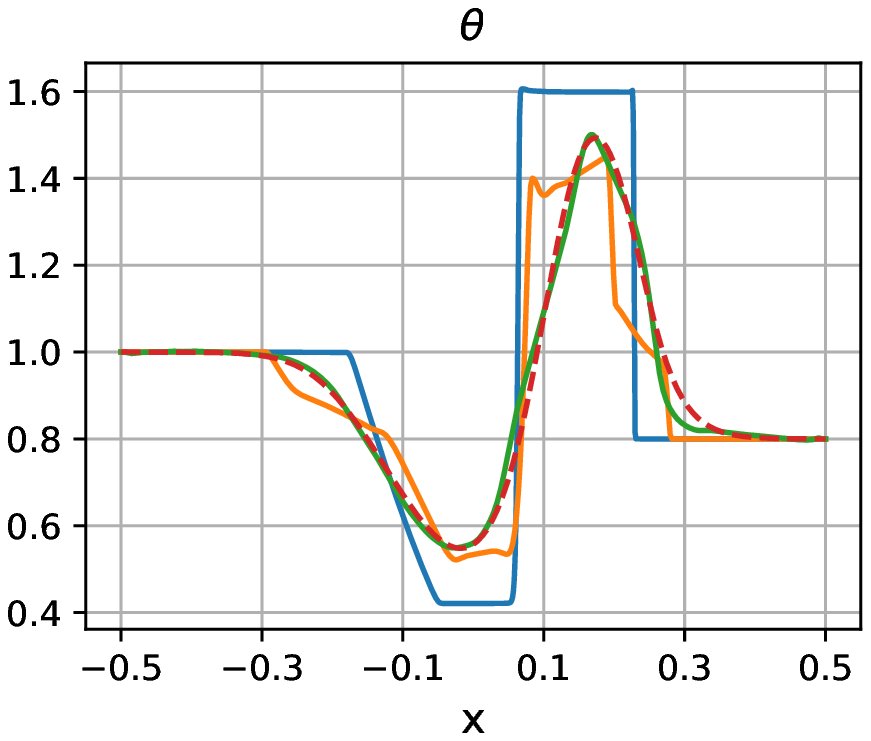}\\
\includegraphics[width=0.3\textwidth]{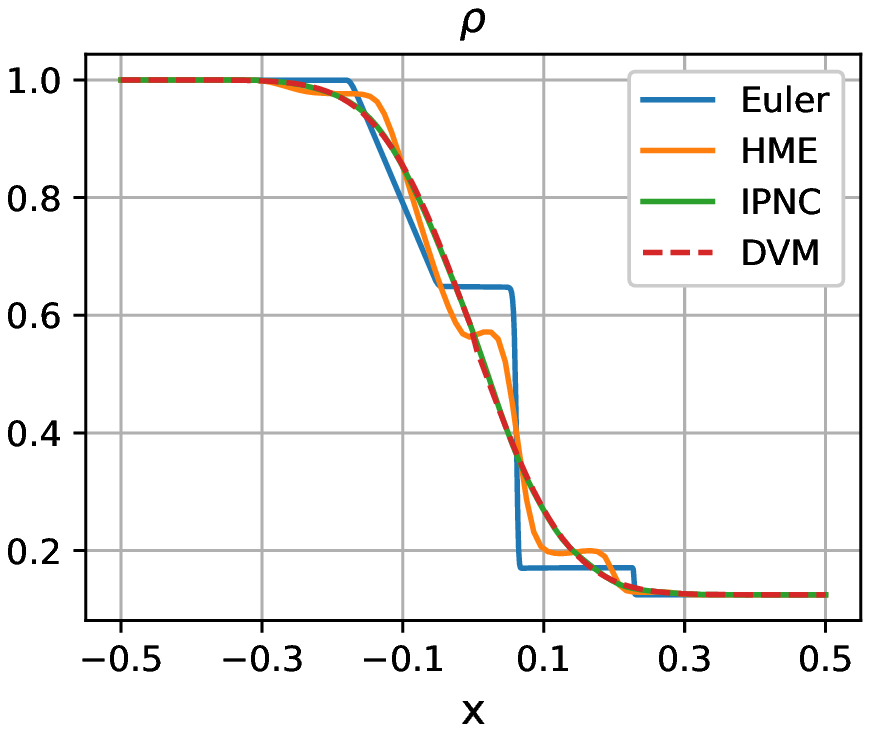}\qquad
\includegraphics[width=0.3\textwidth]{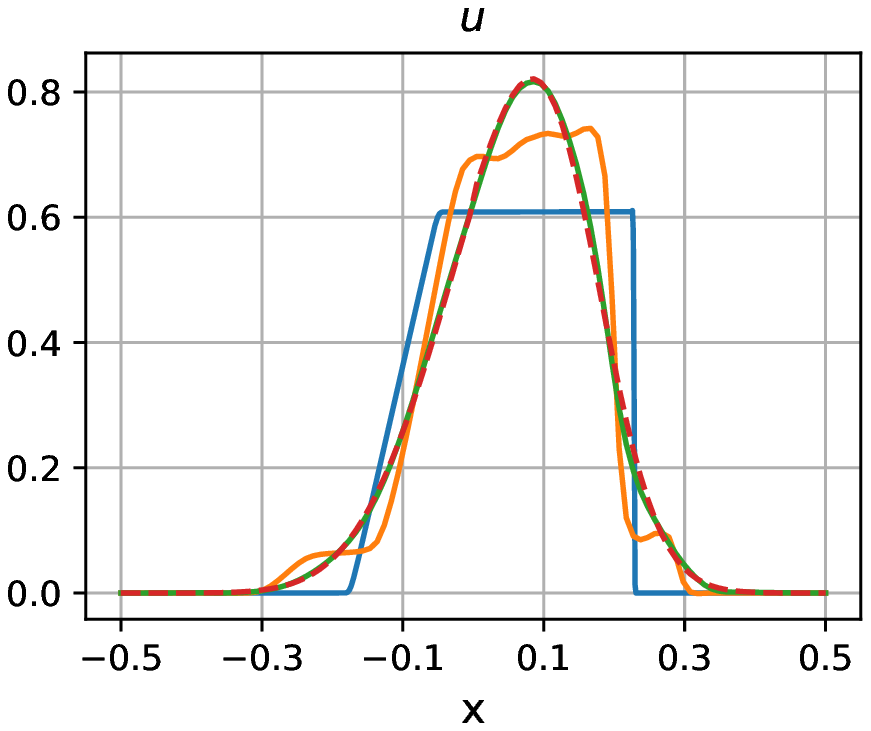}\qquad
\includegraphics[width=0.3\textwidth]{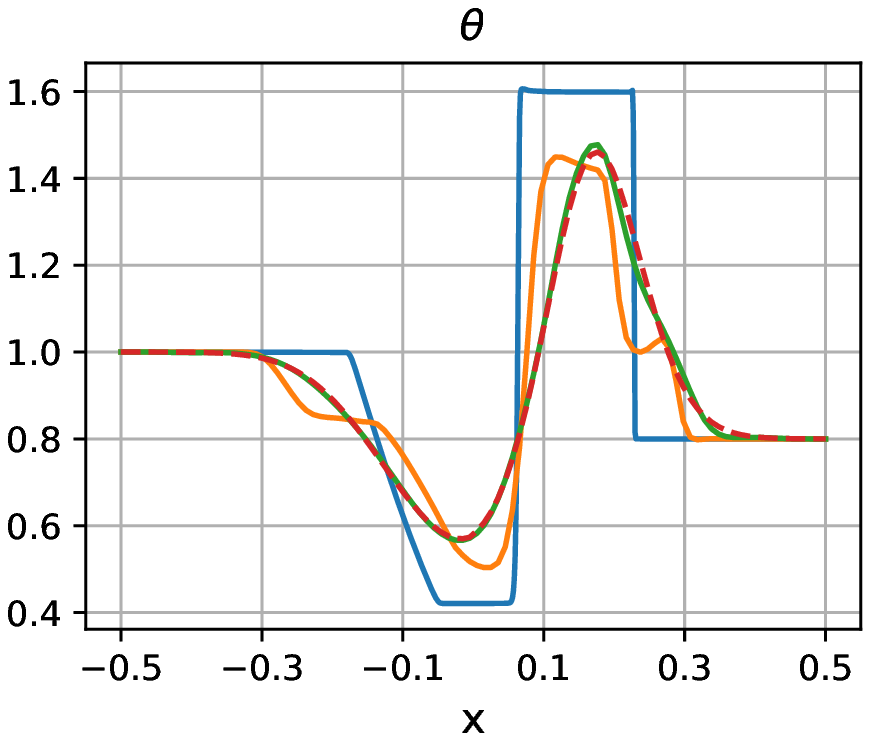}\\
\includegraphics[width=0.3\textwidth]{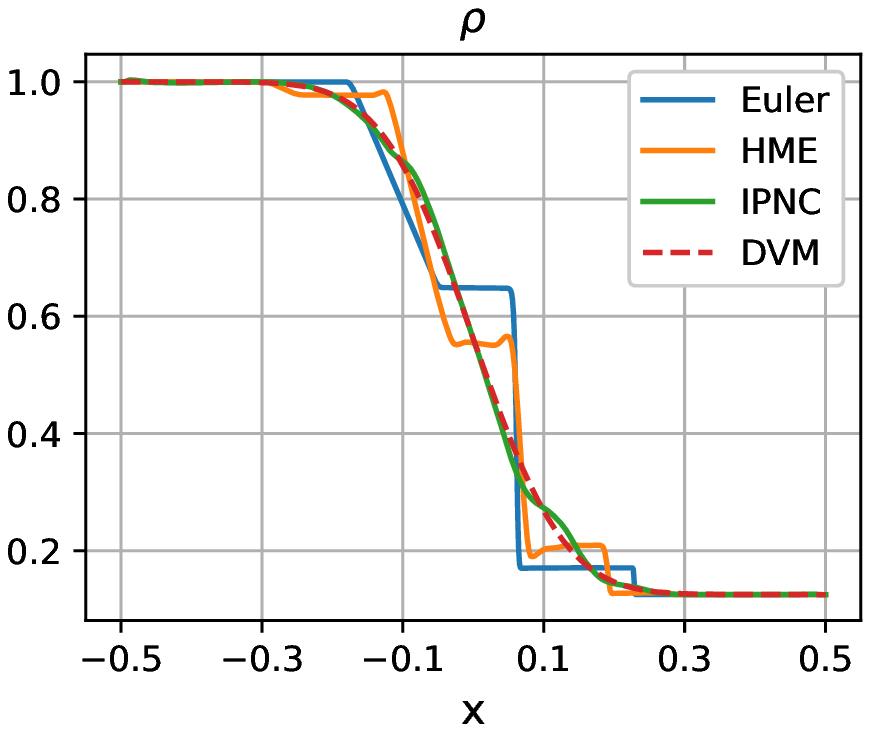}\qquad
\includegraphics[width=0.3\textwidth]{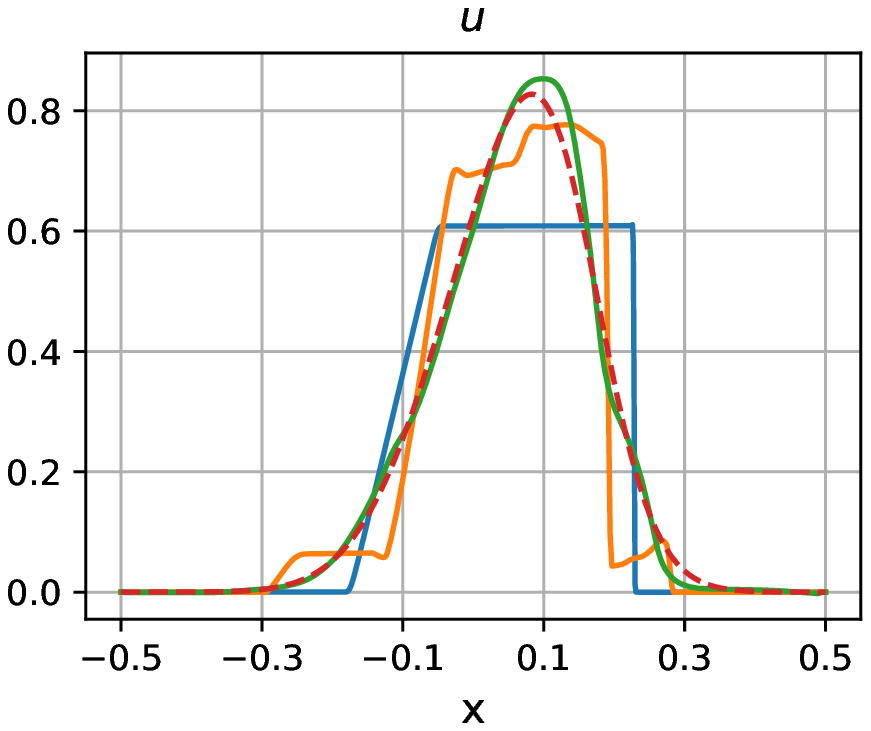}\qquad
\includegraphics[width=0.3\textwidth]{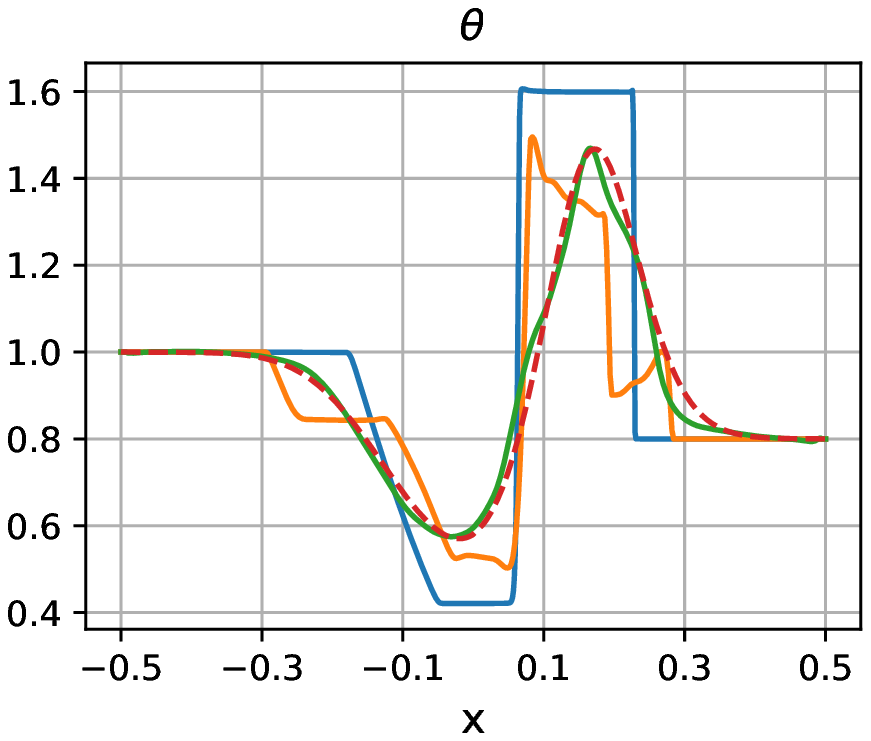}
\caption{(Sec. \ref{sec:num_sod}) The density $\rho$, macroscopic velocity $u$ and temperature $\theta$ of the Sod shock tube problem with different Knudsen numbers at $t = 0.1$. The Knudsen number for each row is $\Kn = 0.001, 0.01, 0.1, 1$ and $10$, respectively. }
\label{fig:sod}
\end{figure}

The numerical solutions by the Euler model, HME, IPNC, and the reference solution by DVM with different Knudsen numbers are studied. For the Euler model, the software Clawpack \cite{mandli2016clawpack} is applied. For HME and IPNC, we use the same Lax-Friedrichs scheme with a $2$nd order linear reconstruction scheme. To test the generalization of IPNC, the network learned from the mix problem in Sec. \ref{sec:num_discon_dis} is directly applied without retraining. The macroscopic variables such as the density $\rho$ and macroscopic velocity $u$ and the temperature $\theta$ at time $t = 0.1$ are plotted. The numerical results with $\Kn = 0.001, 0.01, 0.1, 1.0, 10$ are shown in Figure \ref{fig:sod}. When $\Kn = 0.001$, we can see that the numerical solutions by IPNC and HME are almost the same as the reference solution. Furthermore, there is only a small gap between these numerical solutions and those obtained from the Euler model. With an increase of $\Kn$, this gap is increasing as well. When $\Kn = 0.01$, the numerical solutions by IPNC, HME, and DVM are still close to each other, but the discrepancy between those and the numerical results by the Euler model is larger. When $\Kn$ is increasing to $\Kn = 0.1$, IPNC behaves distinctively better than HME, which indicates that the moment number $M = 5$ of HME is not enough to describe the behavior of the system for $\Kn = 0.1$. Even when $\Kn$ is increasing to $\Kn = 10$, we can still see that the numerical solution by IPNC and the reference solution are almost identical. This shows that even for the rarefied gases, IPNC can still generate satisfactory results with the neural network trained from the mix problem in Sec. \ref{sec:num_discon_dis}. This demonstrates a strong generalization ability of IPNC.

\subsection{Shock structure}
\label{sec:num_shock}
In this section, we test the generalization of the IPNC with respect to the Mach number of the shock structure problem. We shall demonstrate that the IPNC trained on initial values with the Mach numbers within a certain interval can well generalize to the initial values with the Mach number beyond the interval.

For the classical shock structure problem, the initial condition is chosen as 
\begin{equation}
\label{eq:shock_ini}
(\rho, u, \theta)=\left\{\begin{array}{ll}
\left(\rho_{l}, u_{l}, \theta_{l}\right), & x \leqslant 0, \\
\left(\rho_{r}, u_{r}, \theta_{r}\right), & x>0,
\end{array}\right.
\end{equation}
where
\begin{equation}
\begin{gathered}
\rho_{l}=1, \qquad 
u_{l}=\sqrt{3} \Ma, \qquad 
\theta_{l}=1, \qquad \\
\rho_{r}=\frac{2 \Ma^{2}}{\Ma^{2}+1}, \qquad u_{r}=\frac{\sqrt{3} \Ma}{\rho_r}, 
\qquad \theta_{r}=\frac{3 \Ma^{2}-1}{2 \rho_{r}}.
\end{gathered}
\end{equation}
For the shock structure problem, the fluid is initially in local equilibrium, and when the solution evolves for a long enough time, a stable shock structure will be formed.

\begin{figure}[!htb]
\centering
\includegraphics[width=1\textwidth]{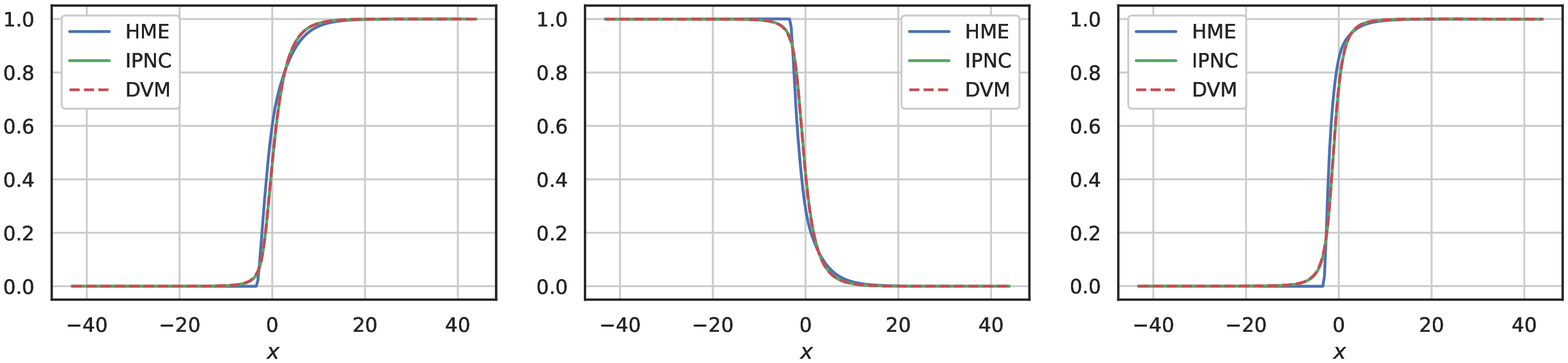}\\
\includegraphics[width=1\textwidth]{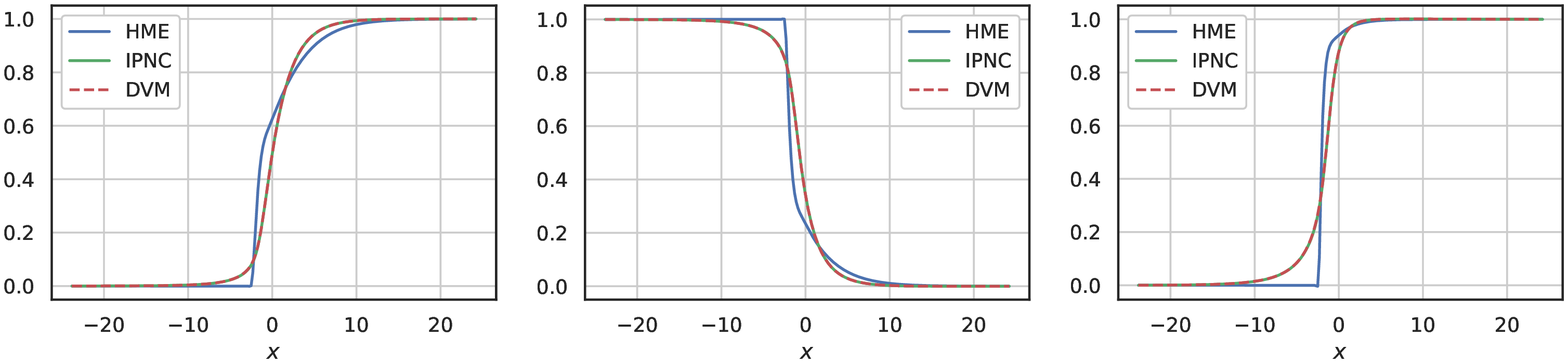}\\
\includegraphics[width=1\textwidth]{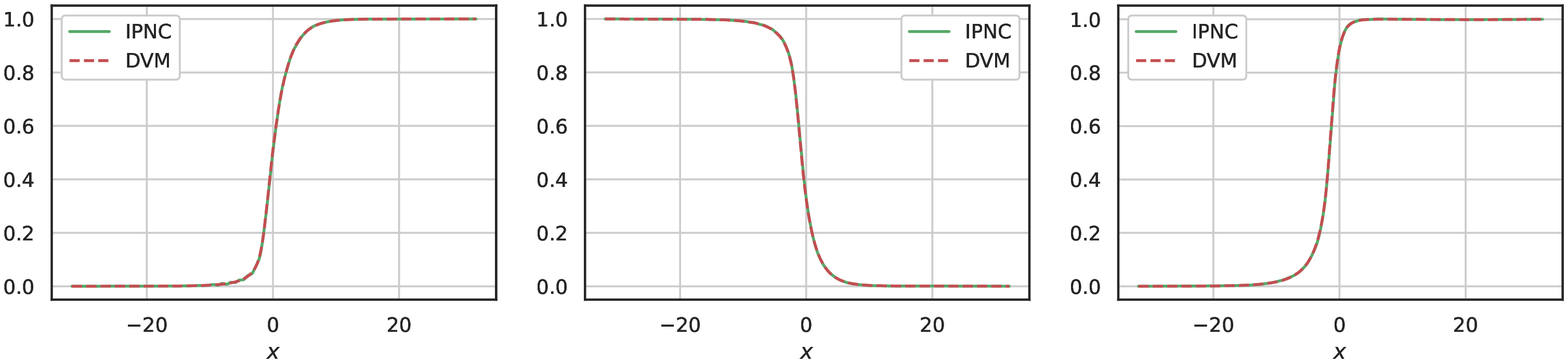}\\
\includegraphics[width=1\textwidth]{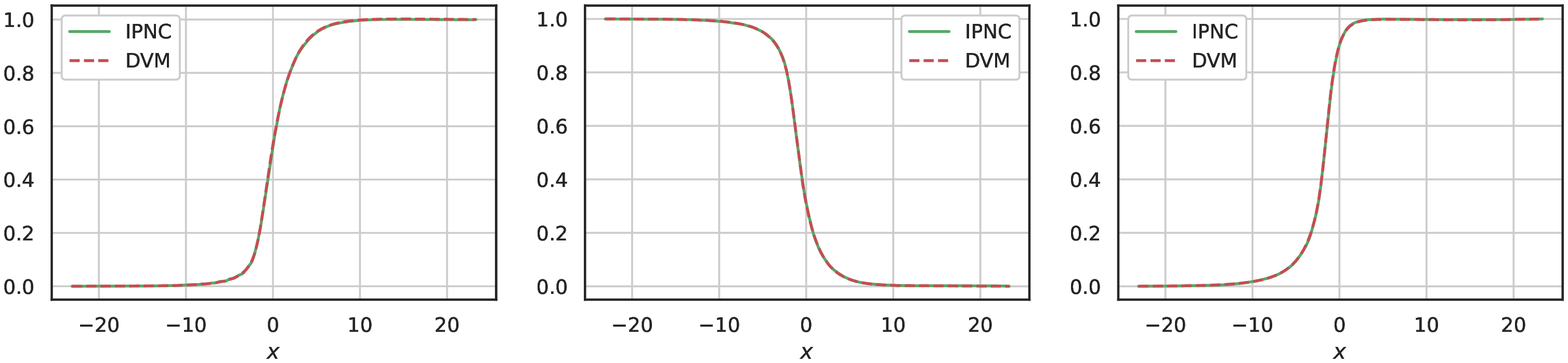}
\caption{(Sec. \ref{sec:num_shock}) The steady state of density $\rho$, macroscopic velocity $u$ and temperature $\theta$ of the shock structure problem with different Mach number and $\Kn = 0.1$. The Mach number for each row is $\Ma = 2, 7, 11$ and $21$, respectively.}
\label{fig:shock}
\end{figure}

In the numerical test, the scaling invariance 
\begin{equation}
\label{eq:shock_ini2}
(\rho, u, \theta)=\left\{\begin{array}{ll}
\left(\hat \rho_{l}, \hat u_{l}, \hat \theta_{l}\right), & x \leqslant 0, \\
\left(\hat \rho_{r}, \hat u_{r}, \hat \theta_{r}\right), & x>0,
\end{array}\right.
\end{equation}
where
\begin{equation}
\begin{gathered}
\hat \rho_{l}=\frac{\rho_l}{\rho_r}, \qquad 
\hat u_{l}=\frac{u_{l}}{\sqrt{\theta_r}}, \qquad
\hat \theta_{l}=\frac{\theta_{l}}{\theta_{r}}, \qquad \qquad 
\hat \rho_{r}=1, \qquad 
\hat u_{r}=\frac{u_r}{\sqrt{\theta_r}}, \qquad \hat \theta_{r}=1,
\end{gathered}
\end{equation}
is utilized to normalize the initial values in \eqref{eq:shock_ini}, which is an equivalent problem but will make the initial values of the density and temperatures laid in the range $[0, 1]$. 

In this experiment, the direct training method in Sec. \ref{sec:dir_learn} is used, and no invariance is constrained on the neural network. For the training data set, it is generated with DVM using the initial value \eqref{eq:shock_ini2}. Instead of U-Net, we use a shared MLP as the backbone of the closure network. This MLP consists of $12$ hidden layers with skip connections, where each layer contains $128$ hidden neurons. LeakyReLU is used as the activation function.

For the training data obtained by DVM, the computational region is $[-10,10]$, with the number of grid $N_x = 400$. The computational region in the microscopic velocity space is $[-20, 20]$ with the number of grid $N_v = 400$. The training set is evolved to $t = 20$ with $100$ time steps taken uniformly within time $[0, 0.1]$, $100$ time steps within $[0.1,1]$, and $200$ time steps within $[1, 20]$. That is, $\mathcal{I}=\mathcal{I}_1\cap\mathcal{I}_2\cap\mathcal{I}_3$, $\mathcal{I}_{1}=\{0\leqslant i \leqslant 100, i\in \mathbb{Z}, t^i \sim{U(0, 0.1)}\}$, $\mathcal{I}_{2}=\{101\leqslant i \leqslant 200,i\in \mathbb{Z}, t^i \sim{U(0.1, 1.0)}\}$, $\mathcal{I}_{2}=\{201\leqslant i \leqslant 400, i\in \mathbb{Z},t^i \sim{U(1.0, 20.0)}\}$. The initial value of the Mach number is from the set $\Ma\in\{1.05+0.1k; 0 < k \leqslant 50, k\in \mathbb{Z}\}$. For the moment method, the same discretization in the spatial space is utilized and the expansion number is set as $M = 5$. 

The numerical solutions with the Mach number $\Ma = 2, 7, 11$ and $21$ are shown in  Figure \ref{fig:shock}, where the density $\rho$, macroscopic velocity $u$, and temperature $\theta$ at the steady-state are plotted. Each variable in the figure is normalized to the range [0,1]. We can see that when the Mach number is small, e.g.,  $\Ma = 2$, the numerical solution by IPNC is almost the same as the reference solution. Here, the numerical solution by HME is also reasonable. When the Mach number increases to $\Ma = 7$, IPNC matches DVM much better than HME. When the Mach number is as large as $\Ma = 11$ and $21$, we can see that the numerical solutions of IPNC are still very close to those of DVM, while we can no longer simulate this problem with such a large Mach number using HME. 

\begin{remark}
We have tried to use the same U-Net as other examples for the shock structure problems. However, with U-Net, we could not obtain satisfactory results for the moment. We have also tested the network with all invariances included and the network trained with the end-to-end method. However, we did not manage to obtain satisfactory results either. Right now, we could not provide a reasonable explanation for this phenomenon or if it is just a problem with our training. Nonetheless, we report this finding and leave the careful diagnosis of this problem as future work. 
\end{remark}

\subsection{Extension to 1D-3D problems}
\label{sec:1d-3d}
\begin{figure}[!htb]
\centering
\subfloat[$\rho$]{
\begin{minipage}[c]{0.3\textwidth}
\includegraphics[width=\textwidth, height=0.75\textwidth]{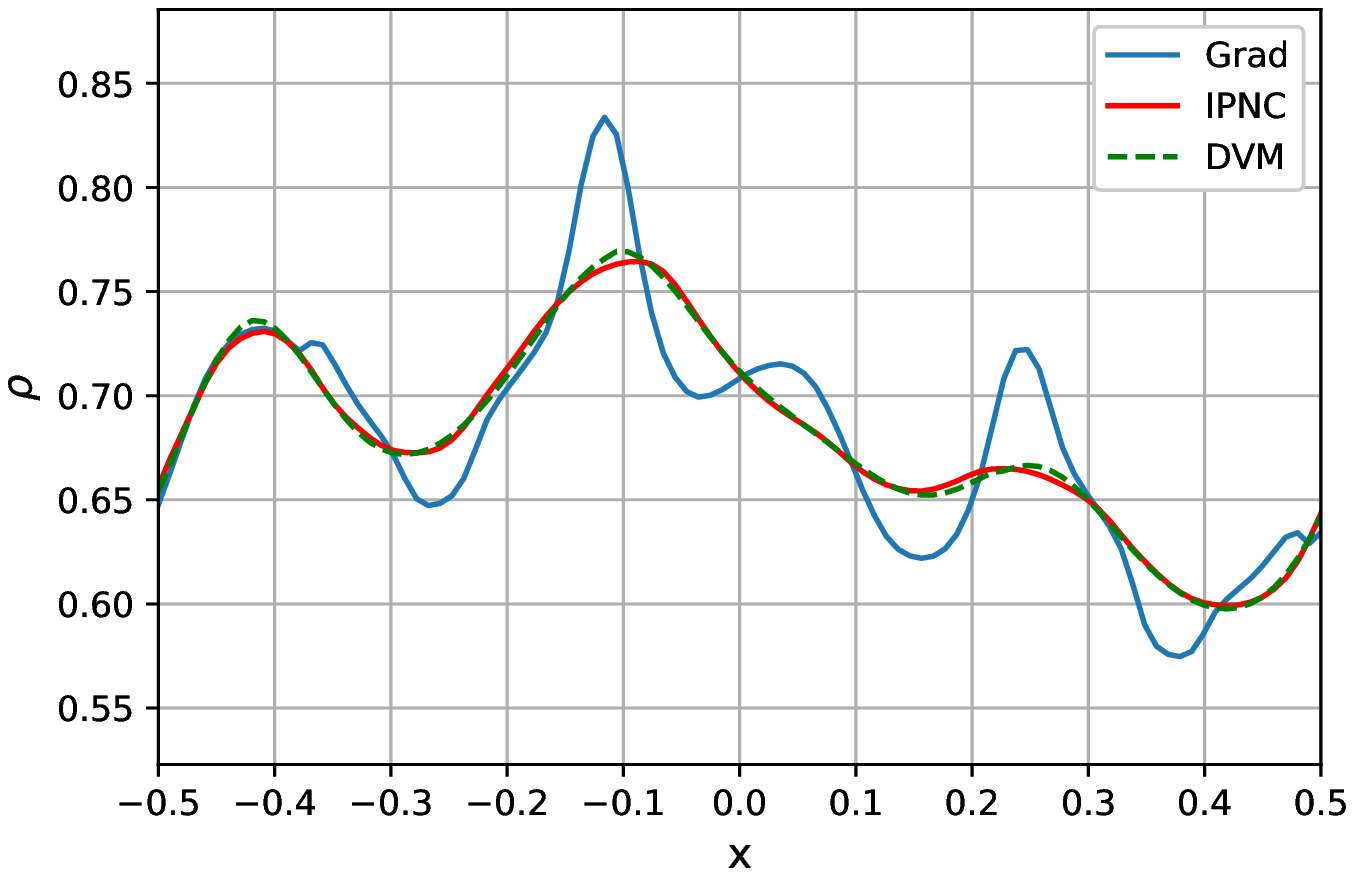}
\includegraphics[width=\textwidth, height=0.75\textwidth]{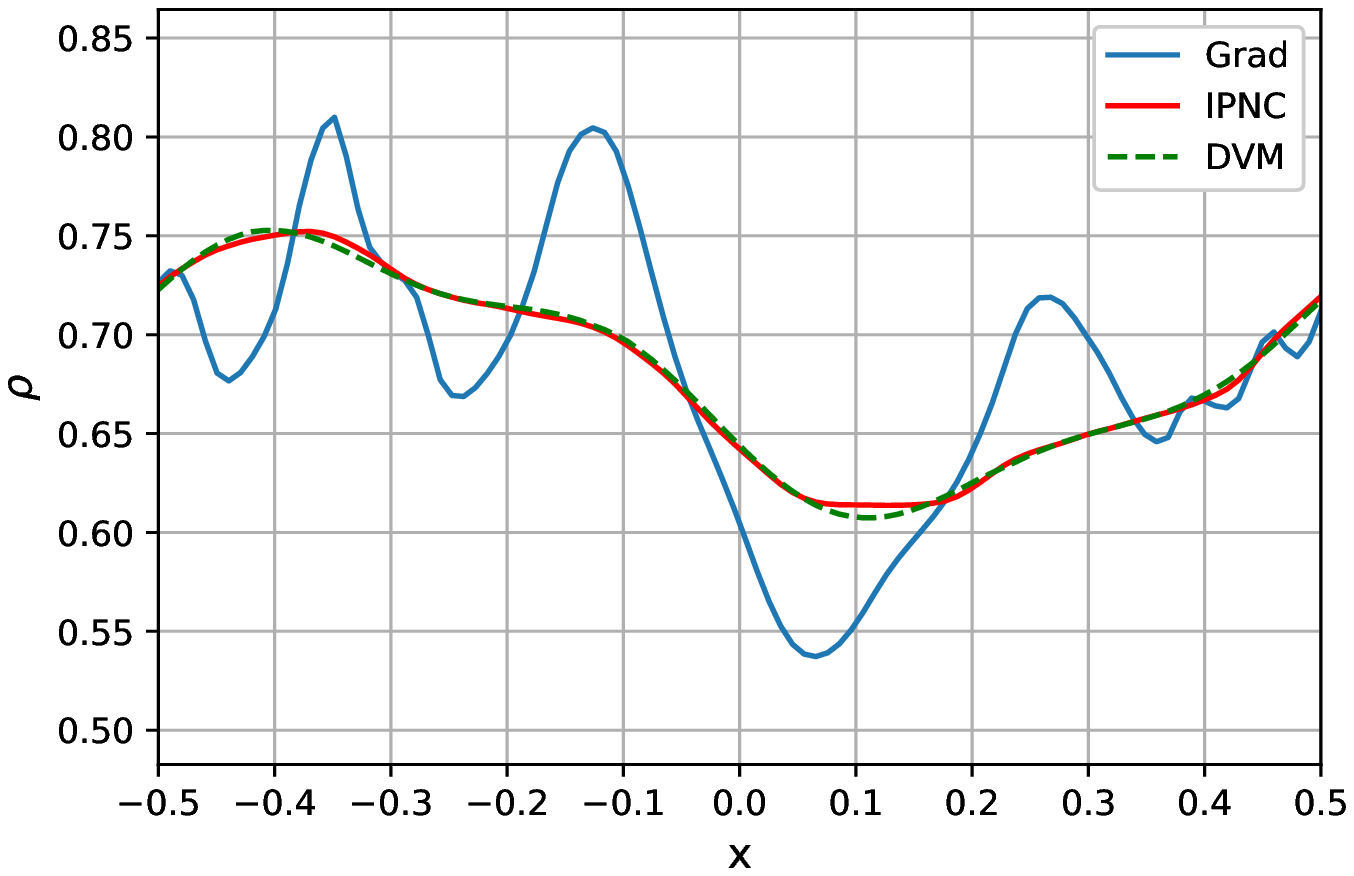}
\end{minipage}
}\quad
\subfloat[$u$]{
\begin{minipage}[c]{0.3\textwidth}
\includegraphics[width=\textwidth, height=0.75\textwidth]{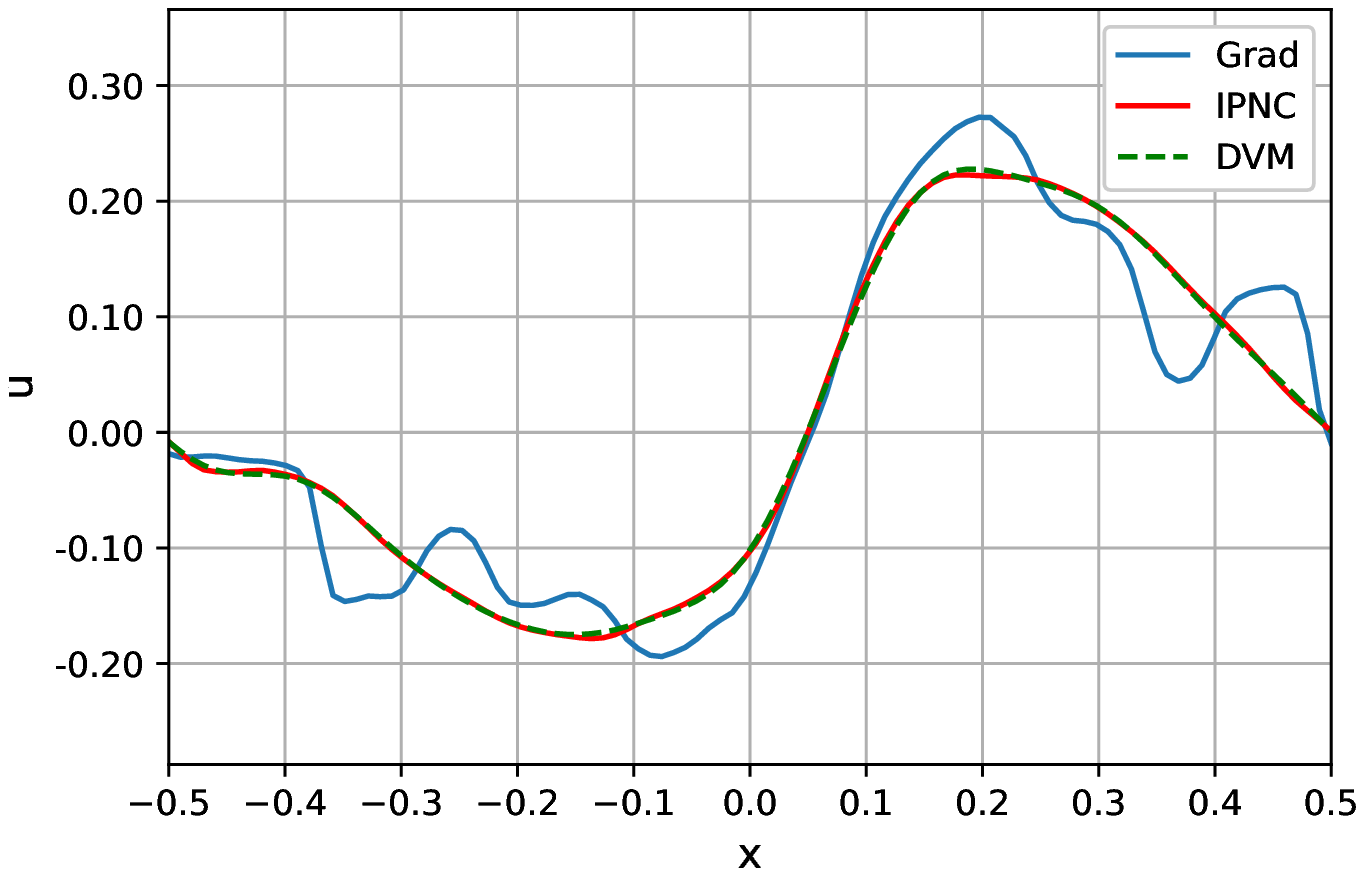}
\includegraphics[width=\textwidth, height=0.75\textwidth]{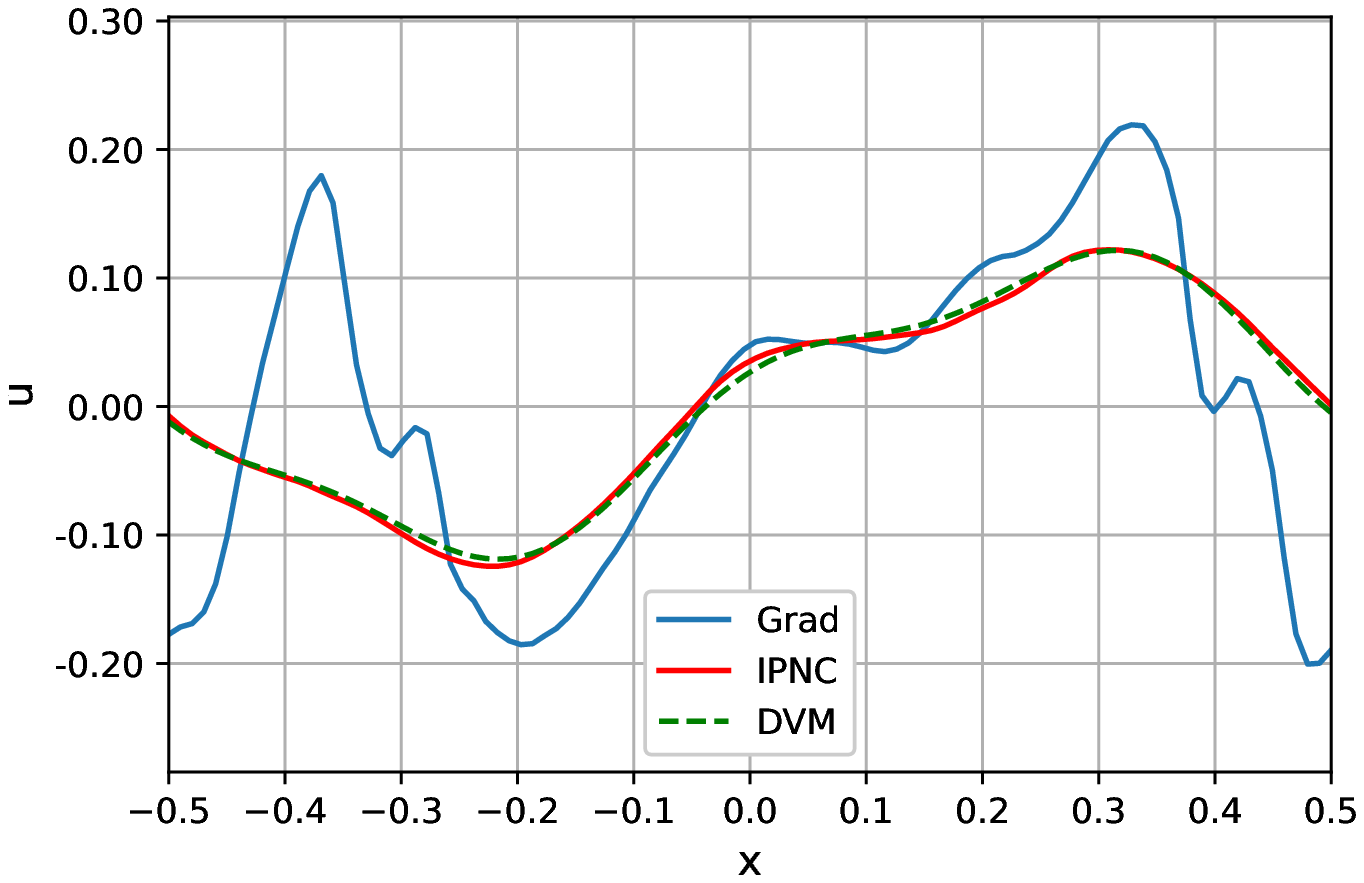}
\end{minipage}
}\quad
\subfloat[$\theta$]{
\begin{minipage}[c]{0.3\textwidth}
\includegraphics[width=\textwidth, height=0.75\textwidth]{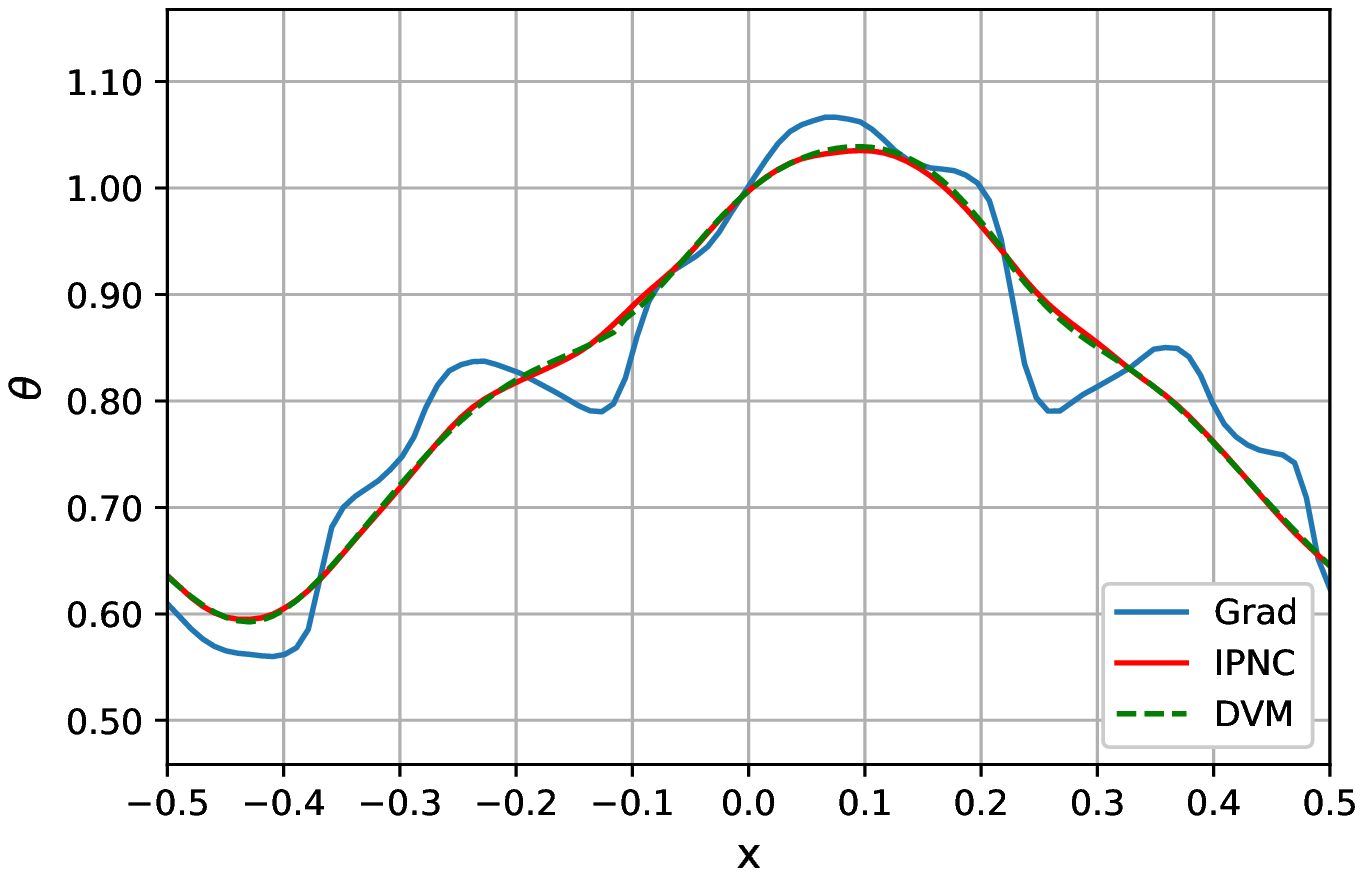}
\includegraphics[width=\textwidth, height=0.75\textwidth]{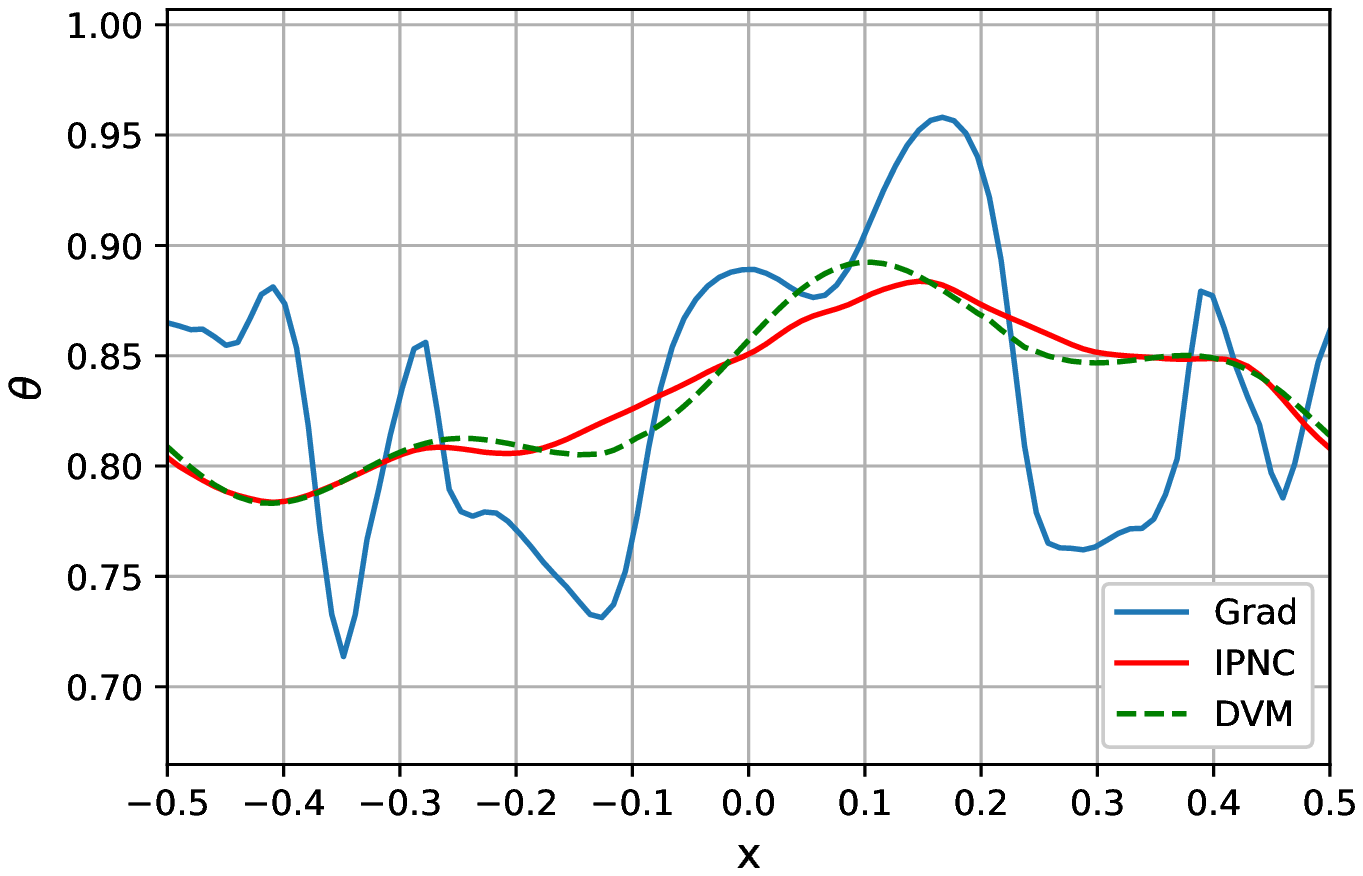}
\end{minipage}
}
\caption {
(Sec. \ref{sec:1d-3d})  Density $\rho$, macroscopic velocity $u$, and temperature $\theta$ at time $t = 0.1$ and $0.2$ for the 1d-3d problem. The two rows are at $t = 0.1$ and $0.2$, respectively. Here the blue line is obtained by Grad, the red line is obtained by IPNC, and the green line is the reference solution obtained by DVM.}
\label{fig:BTE_solution_mix}
\end{figure}

To further validate the IPNC method, we will extend the numerical study to problems in the 1 dimensional spatial and 3-dimensional microscopic velocity space. The same initial condition as 
\eqref{eq:ex2_ini} and \eqref{eq:ex2_rieman} is adopted, but the microscopic velocity space in the distribution function is changed to the case of 3-dimension. The reference solution is still obtained by the DVM method, where the microscopic space chosen is $[-10,10]^3$ and the number of discrete points is $[400,40,40]$. In Grad and IPNC methods, the moments up to the $4$th order are utilized, while the moments at the $5$-th order need to be closed, and compared to the 1D-1D case, the number of moments to close here is $21$. 

\begin{figure}[!htb]
\centering
\subfloat[error quarterlies for Mix]{
\includegraphics[width=0.4\textwidth]{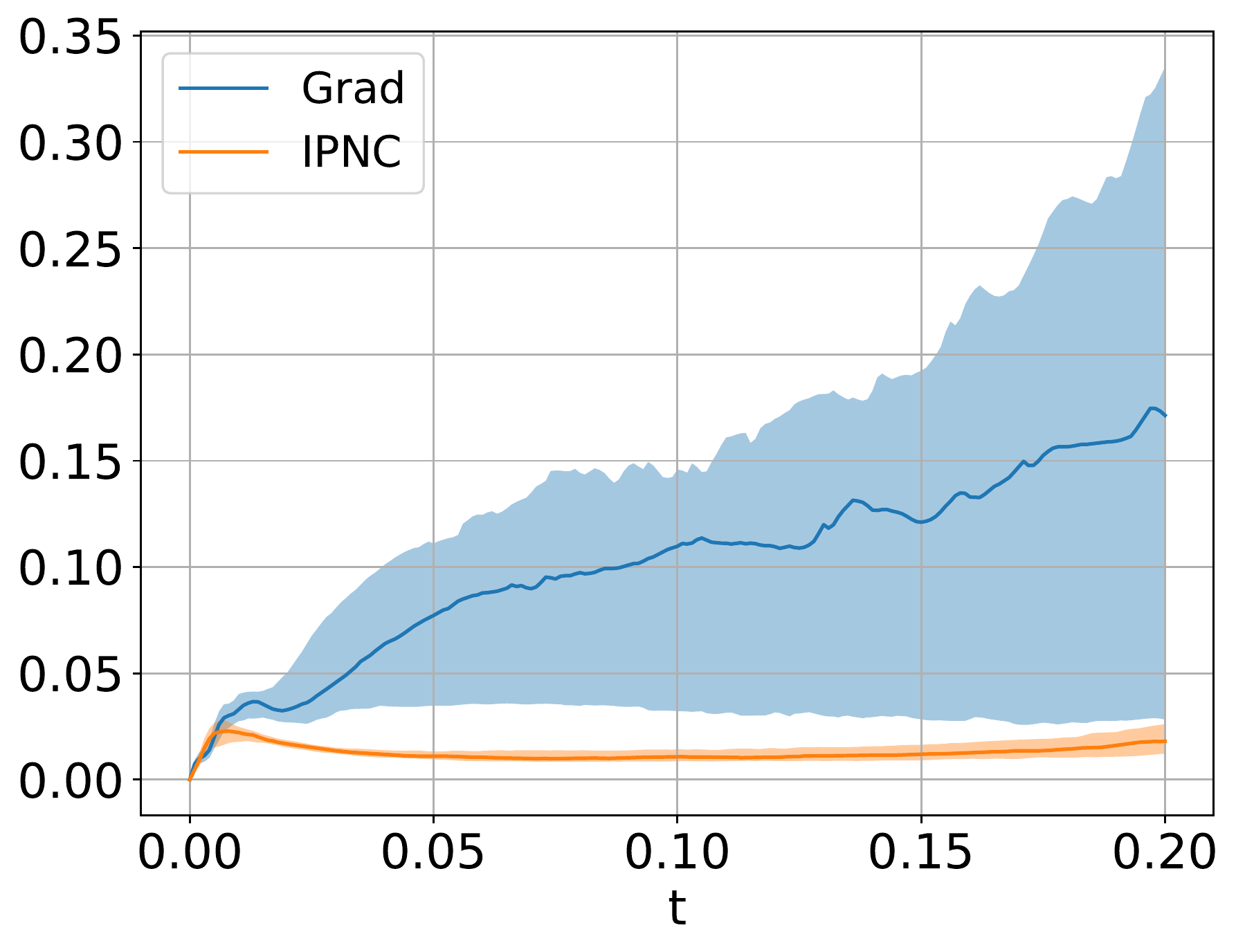}
\label{fig:mixH_solution_error_BTE}} 
\qquad 
\subfloat[different Kn at $t=0.2$ for Mix]{
\includegraphics[width=0.4\textwidth]{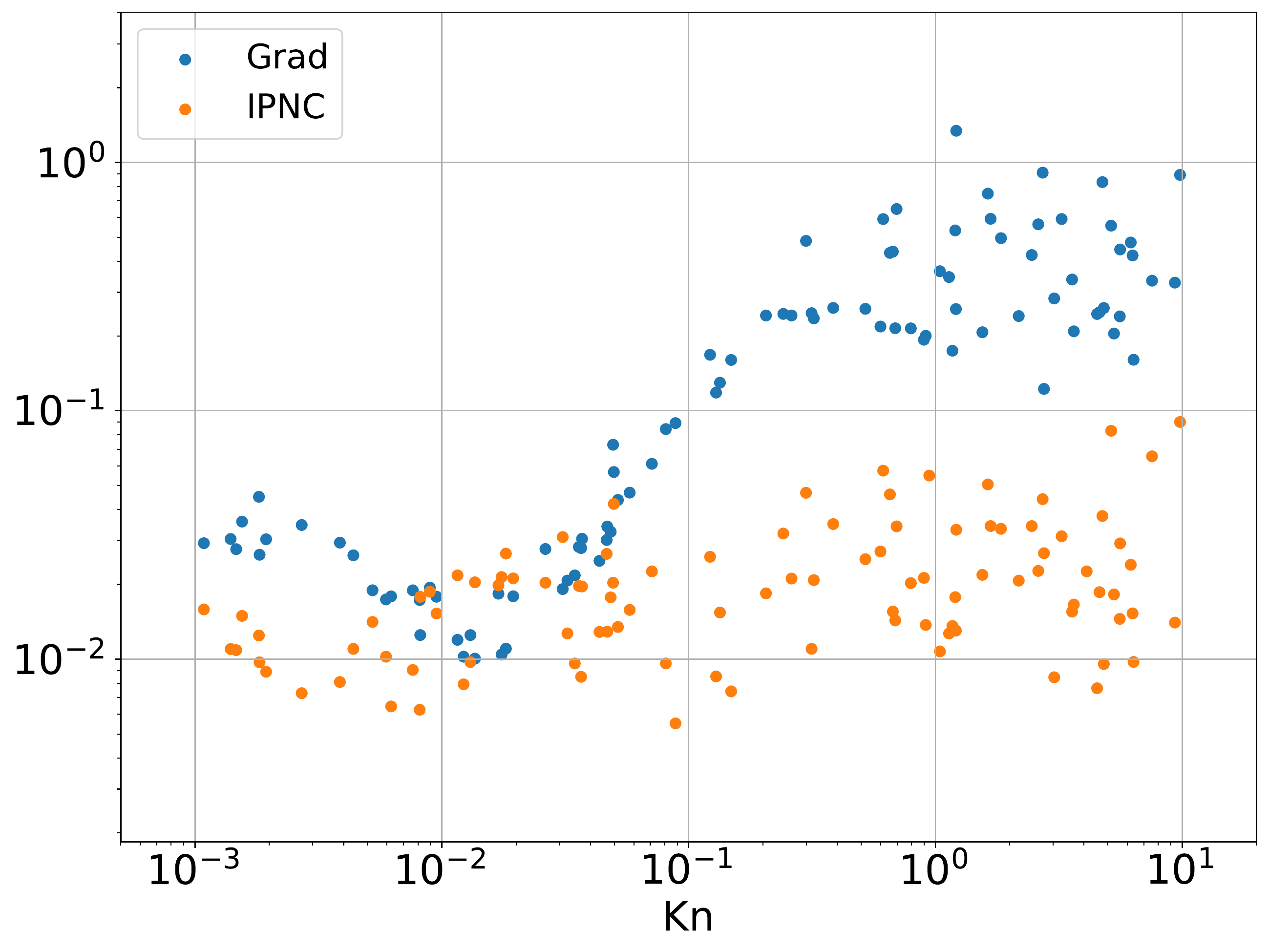}
\label{fig:mixH_scatter_error_BTE}}
\caption{ (Sec. \ref{sec:1d-3d}) (a) The time evolution of the distribution of the relative error with the $100$ initial samples obtained by the different methods, where the $x$-axis is the time $t$ and the $y$-axis is the relative error. (b) The relative error \eqref{eq:macro_error} at time $t = 0.2$ of the $100$ initial samples with different $\Kn$ and different initial conditions. The $x$-axis is the Knudsen number and the $y$-axis is the relative error.}
\label{fig:mixH_error_BTE3D}
\end{figure}

The same training data set and testing data set as in Sec. \ref{sec:num_discon_dis} are adopted. Figure \ref{fig:BTE_solution_mix} shows the behavior of density $\rho$, macroscopic velocity $u$ and temperature $\theta$ at time $t = 0.1$ and $0.2$ for the first sample of the testing data, where the numerical solutions by Grad and the reference solution by DVM are all plotted. It shows that for the 1D-3D case, at $t = 0.1$ the numerical solution by IPNC matches well with the reference solution, while those by Grad have some error from the DVM solution. However, at $t=0.2$, there exists some distance between the numerical solution by IPNC and DVM. For the moment, the generalization ability of IPNC in the 1D-3D case does not perform so well, which may be due to the same network being utilized here as the 1D-1D case, but the freedom to close is increasing greatly. The similar distribution of the relative error \eqref{eq:macro_error} for different initial conditions with different methods is shown in Figure \ref{fig:mixH_solution_error_BTE}, and the stability of IPNC with respect to the Knudsen number $\Kn$ at $t = 0.2$ is shown in Figure \ref{fig:mixH_scatter_error_BTE}. Both figures validate the stability of the IPNC on the Knudsen number for the 1D-3D problem.

For the 1D-3D problem, the efficiency of the DVM is significantly decreased, while the moment method shows a great speed advantage at this point. Moreover, since the same network as the 1D-1D case is utilized, the computational cost for closure also remains the same as in the 1D-1D case. Therefore, the IPNC method may greatly reduce computational costs compared to the DVM method. 
We have counted the computational time with $\Kn=10$ to $t = 0.2$ for DVM, Grad and IPNC, which is shown in Table \ref{tab:BTE3Dtime}. Here, the mesh size in microscopic velocity space is $[120, 40, 40]$ in the DVM method, and $[60, 20, 20]$ for the coarse-DVM method. For the moment method, Grad and IPNC, the $4$th order expansion is utilized. The mesh size in macroscopic space is $N_x=400$ for all methods. All these methods use the GPU for acceleration, and the calculations are done on a single RTX3090. The average of the relative error \eqref{eq:macro_error} is also shown in Table \ref{tab:BTE3Dtime}, where the reference solution is obtained by the DVM method with the discretization $[120, 40, 40]$. Compared to the Grad method, the computational cost of IPNC is not greatly increased and is much lower than both of the DVM methods, while the average of the relative error is much smaller than the other two methods, which validates the high efficiency of the IPNC method.

\begin{table}[!htb]
\centering
\renewcommand\arraystretch{1.5}
\footnotesize
% \begin{tabular}{l||l|l|l}
% Method  & Grad & IPNC & DVM \\ \hline  
% Time (sec.)  &  3   &  7  & 332 \\ 
% \end{tabular}
\begin{tabular}{l|cccc}
Method      & Grad  & IPNC &  coarse-DVM & DVM   \\ \hline
Time (sec.) & 6.77  & 8.38 & 17.5      & 192.4  \\
Error (Wave) & 53.9 & 4.31 & 15.4   & 0     \\
Error (Mix)  & -     & 3.44 & 5.48  & 0    
\end{tabular}
\caption{(Sec. \ref{sec:1d-3d}) The computational time taken to $t = 0.2$ by Grad, IPNC, coarse-DVM, and DVM for one initial condition, and the average of the relative error for these methods with the wave and mix initial conditions, where the numerical solution obtained by DVM is treated as the reference solution.
}
\label{tab:BTE3Dtime}
\end{table}

\subsection{Ablation experiment}
\label{sec:abl}
In this section, we demonstrate the necessity of the invariances we have proposed in Sec. \ref{sec:sym_net}. The performance of the different training approaches as end-to-end learning in Sec. \ref{sec:end_end} and direct learning in Sec. \ref{sec:dir_learn} is also studied. 

\subsubsection{Invariance preservation}
\label{sec:inv_pre}
In Sec. \ref{sec:sym_net},  we have proposed the network with the three invariances. Here, we quantitatively show the effects of these invariances on the numerical results. The experiments are done on the wave problem in Sec. \ref{sec:num_smooth_dis} and the mix problem in Sec. \ref{sec:num_discon_dis}. In the experiments, two different kinds of testing data sets are chosen. In the first kind of data set, samples are generated with the same distribution as the training data set, while in the second data set, the samples are generated outside the distribution of the training data set. The testing on data points that are generated outside of the distribution of the training data is called the out-of-distribution (OoD) data, which is quite important for a practical concern.

For both wave and mix problems, the training and first testing data set is the same as that in Sec. \ref{sec:num_smooth_dis} and \ref{sec:num_discon_dis}, respectively. However, for the second testing data, i.e., the OoD testing data set of the wave problem (denoted as wave-OoD), $a_{\rho}$ and $a_{\theta}$ in \eqref{eq:ex1_ini} is uniformly sampled from $[0.04, 0.06]$, and $b_{\rho}$ and $b_{\theta}$ are uniformly sampled from $[0.10, 0.14]$. For the mix-OoD data set, $\rho_l$ and $\theta_l$ are sampled from the uniform distribution on $[0.2, 0.4]$ in \eqref{eq:ex2_rieman} while $\rho_r$ and $\theta_r$ are sampled from the uniform distribution on $[0.11, 0.18]$. 

We first train the IPNC network with the data from the wave problem and then test the effect of this network on the wave and mix problems with both testing data sets. The same numerical setting as in Sec. \ref{sec:num_smooth_dis} is utilized here. For example, $100$ samples of training and testing data are generated, and the terminal time is $t = 0.1$. The same average of the relative error of the numerical solutions among the $100$ samples is calculated. 

\begin{table}[!htb]
\centering
\renewcommand\arraystretch{1.5}
\footnotesize
\begin{tabular}{ccc||cccc||cccc}
\multirow{2}{*}{GI} & \multirow{2}{*}{SI} & \multirow{2}{*}{RI} & \multicolumn{4}{c||}{Train on Wave}                                            & \multicolumn{4}{c}{Train on Mix}                                              \\ \cline{4-11} 
                    &                     &                     & Wave              & Wave-OoD          & Mix               & Mix-OoD           & Wave              & Wave-OoD          & Mix               & Mix-OoD           \\ \hline \hline
                    &                     &                     & $ 0.63 $          & $ 2.16 $          & $ 4.33 $          & $ 2.44 $          & $ 1.50 $          & -                 & $ 0.62 $          & -                 \\
                    &                     & $\checkmark$        & $0.61$            & $1.89$            & $4.25$            & $2.12$            & $1.38$            & -                 & $ 0.60 $          & -                 \\
$\checkmark$        &                     &                     & $ 0.63 $          & $ 2.9 $           & $ 1.53 $          & $ 3.43 $          & $ 1.62 $          & -                 & $ 0.61 $          & -                 \\
                    & $\checkmark$        &                     & $ 0.61 $          & $ \mathbf{0.25} $ & $ 1.53 $          & $ 0.86 $          & $ 0.96 $          & $ 0.32 $          & $ 0.61 $          & $ \mathbf{0.52} $ \\
                    & $\checkmark$        & $\checkmark$        & $ \mathbf{0.58} $ & $ \mathbf{0.25} $ & $ 1.48 $          & $ \mathbf{0.85} $ & $ \mathbf{0.92} $ & $ \mathbf{0.31} $ & $ \mathbf{0.59} $ & $ \mathbf{0.52} $ \\
$\checkmark$        & $\checkmark$        &                     & $ 0.62 $          & $ 0.26 $          & $ 1.51 $          & $ 0.91 $          & $ 1.11 $          & $ 0.36 $          & $ 0.63 $          & $ 0.53 $          \\
$\checkmark$        & $\checkmark$        & $\checkmark$        & $ 0.59 $          & $\mathbf{ 0.25 }$ & $ \mathbf{1.47} $ & $ 0.90 $          & $ 1.05 $          & $ 0.33 $          & $ \mathbf{0.59} $ & $ \mathbf{0.52} $
\end{tabular}
\caption{(Sec. \ref{sec:abl}) Results of the ablation experiments on invariance. The left 3 columns label which invariants are embedded in the network. Columns 4-7 show the test errors on different datasets for the network obtained by training on the Wave problem. Columns 8-11 show the test errors on different datasets for the network trained on the Mix problem. Each row represents the prediction error of a network embedded with different invariances on four data sets. ``-'' means that the method failed at least in one test sample.}
\label{tab:ablation-invariants}
\end{table}

To demonstrate the effect of embedded invariances, we designed comparative experiments to embed only partial symmetries in the network, respectively. We use GI to denote Galilean invariance, SI to denote scaling invariance, and RI to denote the reflecting invariance.
Table \ref{tab:ablation-invariants} shows the average of the relative errors \eqref{eq:macro_error} of IPNC with different invariances. We can see that for the first testing data set of the wave problem, this error is always relatively small. The invariances of the network do not improve much on this result. However, for wave-OoD, it is obvious that the error is greatly reduced for the neural network with invariances. The error reaches the minimum value with the neural network having all three invariances. A similar phenomenon can be observed on mix and mix-OoD, where having as many invariances in the neural network as possible is beneficial.

On the other hand, we train the IPNC network on the mix problem and redo the tests on the four testing data sets. The last four columns in Table \ref{tab:ablation-invariants} show the performance of networks trained on mix and generalized on wave, wave-OoD, and mix-OoD by enforcing different sets of invariances. A similar conclusion can be made as in the previous tests that having more invariances helps with the trained neural network in terms of generalization.

% In another experiment, we tested small-sample learning, and it can be seen from the experimental results (Table \ref{tab:ablation-invariants-2}) that the closure with the invariant performed significantly better than the closure with the normal structure.

% \begin{table}[]
% \centering
% \begin{tabular}{|l|l|l|}
% \hline
%           & Mix    & Wave   \\ \hline
% None      & 0.6538 & 0.5546 \\ \hline
% Invariant & 0.6988 & 0.5170 \\ \hline
% \end{tabular}
% \caption{}
% \label{tab:ablation-invariants-2}
% \end{table}

\subsubsection{End-to-end and direct learning approach}
\label{sec:end-dir}

In this section, we compare the performance of the end-to-end learning and the direct learning approach.
We retrain the network in Sec.\ref{sec:num_smooth_dis} and \ref{sec:num_discon_dis} using direct learning and end-to-end learning with different blocks. The length of the block in end-to-end learning is chosen as $B = 1, 2, 4, 8$. The numerical setting here is the same as that in Sec. \ref{sec:num_smooth_dis} and \ref{sec:num_discon_dis}. Table \ref{tab:ablation-invariants-2} presents the average of the relative errors for different methods. We can see that end-to-end learning is better than direct learning, and having more blocks for end-to-end learning is beneficial, which is also consistent with the analysis in Sec. \ref{sec:spe_app}.

\begin{table}[ht]
\centering
\renewcommand\arraystretch{1.5}
\footnotesize
\begin{tabular}{l||lllll}
     & Direct & Block 1 & Block 2 & Block 4 & Block 8 \\ \hline \hline
Wave &   $0.750$     &  $ 0.645 $   &  $ 0.639 $   &  $ 0.618 $   &  $ \mathbf{0.616 }$   \\
Mix  &   $0.792$     &  $ 0.662 $   &  $ 0.657 $   &  $ 0.634 $   &  $ \mathbf{0.593 }$  
\end{tabular}
\caption{(Sec. \ref{sec:end-dir}) Average of the relative error \eqref{eq:macro_error} for the wave and mix problem with different learning approaches.}
\label{tab:ablation-invariants-2}
\end{table}

%%% Local Variables: 
%%% mode: latex
%%% TeX-master: "article"
%%% End: 

%% file: article_conclusion.tex
\section{Conclusion}
\label{sec:conclusion}
In this paper, we proposed an invariance preserving neural closure (IPNC) method for the Boltzmann-BGK equation under the framework of the Grad-type moment method. The network was particularly designed so that the physical system preserves the Galilean, reflecting and scaling invariance. We tested the performance of the IPNC method on the problems with smooth and discontinuous initial conditions. Numerical results showed that compared to the benchmark methods HermMLC and HME, the averaged relative error of IPNC was much smaller, and the IPNC was also more stable to the Knudsen number. The generalization of IPNC was tested on the Sod's shock tube and shock structure problem. We demonstrated that IPNC trained with discontinuous initial conditions can be directly applied to Sod's shock tube problem without retraining. In the classical shock structure problem, the generalization of the IPNC with respect to the Mach number was also tested. We demonstrated that the IPNC trained on initial values with the Mach numbers within a specific interval could well generalize to initial values with the Mach number significantly beyond the interval. The extension to the 1D-3D problem shows the potential IPNC has to solve the high-dimensional problems.  

This is our first attempt to seek moment closure with deep neural networks that preserve physical invariances, where we only worked on the 1D BGK model. We will consider models with more complex collisions, such as the quadratic collision, as part of future work. We will also consider the problem in 2-3 dimensional spaces and higher order of moments.

\section*{Acknowledgements}
We thank Prof. Ruo Li from Peking University, and Dr. Jiequn Han from Flatiron Institute for their valuable suggestions. Zhengyi Li and Bin Dong are supported in part by Natural Science Foundation of Beijing Municipality (No. 180001) and National Natural Science Foundation of China (Grant No. 12090022). Yanli Wang is supported by the National Natural Science Foundation of China (Grant No. 12171026, U1930402 and 12031013).

%%% Local Variables: 
%%% mode: latex
%%% TeX-master: "article"
%%% End: 

%% file: appendix.tex
\section{Appendix}
\label{app:supp}
In the appendix, we will provide the detailed initial condition of the first sample in the wave and mix problem. 

The detailed initial condition for the first sample in Sec. \ref{sec:num_smooth_dis} is in Table \ref{tab:smooth_first}, the numerical results of which are plotted in Figure \ref{fig:ex1_solution}. The detailed initial condition for the first sample in Sec. \ref{sec:num_discon_dis} is listed in Table \ref{tab:mix_first}, the numerical results of which is plotted in Figure \ref{fig:ex2_solution}. 
\begin{table}[ht]
    \centering
    \renewcommand\arraystretch{1.5}
    \footnotesize
    \begin{tabular}{l||llll}
        $U_1$ &$a_{\rho}^1=0.24009$   & $b_{\rho}^1=0.54997$  & $\phi_{\rho}^1=1.66063$ &$k_{\rho}^1=4$ \\
        & $a_{\theta}^1=0.28894$ & $b_{\theta}^1=0.67490$ & 
         $\phi_{\theta}^1=5.14649$ &$k_{\theta}^1=1$  \\[2mm]
         \hline 
        $U_2$ & $a_{\rho}^2=0.23105$ & $b_{\rho}^2=0.69649$ &$\phi_{\rho}^2=2.30314$ & $k_{\rho}^2=1$ \\
        & $a_{\theta}^2=0.25016$ & $b_{\theta}^2=0.65446$ & $\phi_{\theta}^2=2.72434$ & $k_{\theta}^2=2$ \\
        \hline 
        & $\alpha_1=0.51229$ & $\alpha_2=0.93889$ & $\Kn= 4.33683$ 
    \end{tabular}
    \caption{The initial condition of the initial condition sample  for the wave problem.}
    \label{tab:smooth_first}
\end{table}
\begin{table}[ht]
    \centering
    \renewcommand\arraystretch{1.5}
    \footnotesize
    \begin{tabular}{l||llll}
        $U_1$ &$a_{\rho}^1=0.29053$   & $b_{\rho}^1=0.59197$  & $\phi_{\rho}^1=3.76432$ &$k_{\rho}^1=1$ \\
        & $a_{\theta}^1=0.25770$ & $b_{\theta}^1= 0.52530$ & 
         $\phi_{\theta}^1=0.09759$ &$k_{\theta}^1=3$  \\[2mm]
       %  \hline 
        $U_2$ & $a_{\rho}^2= 0.25934$ & $b_{\rho}^2=0.58673$ &$\phi_{\rho}^2=1.15617$ & $k_{\rho}^2=3$ \\
        & $a_{\theta}^2 = 0.22352$ & $b_{\theta}^2=0.66302$ & $\phi_{\theta}^2=4.95070$ & $k_{\theta}^2=2$ \\[2mm]
      %  \hline 
        & $\alpha_1 = 0.61203$ & $\alpha_2 = 0.05390$ & $\Kn = 6.19271$ \\[1mm]
        \hline  \hline 
        $U_l$ & $\rho_l=0.59584$ & $u_l = 0$ & $\theta_l=0.67065$ &  \\
        $U_r$ & $\rho_r=1.01501$ & $u_r = 0$ & $\theta_r=1.33206$ & \\[2mm]
       % \hline 
        & $x_1 = -0.11197$ & $x_2 = 0.21640$ & $\alpha = 0.36807$ 
    \end{tabular}
    \caption{The initial condition of the initial condition sample for the mix problem.}
    \label{tab:mix_first}
\end{table}